%% file: partII.tex
\documentclass[11pt,reqno]{amsart}

\usepackage{amssymb,amsmath,amsthm,amsfonts}
\usepackage{hyperref, multirow, bm, color, marginnote, graphicx, wrapfig, subcaption,esint}
\usepackage{diagbox}

\makeatletter
\renewcommand{\tocsection}[3]{%
  \indentlabel{\@ifnotempty{#2}{\bfseries\ignorespaces#1 #2\quad}}\bfseries#3}
\renewcommand{\tocsubsection}[3]{%
  \indentlabel{\@ifnotempty{#2}{\ignorespaces#1 #2\quad}}#3}

\newcommand\@dotsep{4.5}
\def\@tocline#1#2#3#4#5#6#7{\relax
  \ifnum #1>\c@tocdepth 
  \else
    \par \addpenalty\@secpenalty\addvspace{#2}%
    \begingroup \hyphenpenalty\@M
    \@ifempty{#4}{%
      \@tempdima\csname r@tocindent\number#1\endcsname\relax
    }{%
      \@tempdima#4\relax
    }%
    \parindent\z@ \leftskip#3\relax \advance\leftskip\@tempdima\relax
    \rightskip\@pnumwidth plus1em \parfillskip-\@pnumwidth
    #5\leavevmode\hskip-\@tempdima{#6}\nobreak
    \leaders\hbox{$\m@th\mkern \@dotsep mu\hbox{.}\mkern \@dotsep mu$}\hfill
    \nobreak
    \hbox to\@pnumwidth{\@tocpagenum{\ifnum#1=1\bfseries\fi#7}}\par
    \nobreak
    \endgroup
  \fi}
\AtBeginDocument{%
\expandafter\renewcommand\csname r@tocindent0\endcsname{0pt}
}
\def\l@subsection{\@tocline{2}{0pt}{2.5pc}{6pc}{}}
\makeatother

\usepackage[margin=1in]{geometry}

\newcommand{\R}{\mathbb{R}}

\newcommand{\E}{\mathcal{E}}

\newcommand{\bu}{\bm{u}}
\newcommand{\bw}{\bm{w}}
\newcommand{\bx}{\bm{x}}
\newcommand{\by}{\bm{y}}

\newcommand{\X}{\bm{X}}

\newcommand{\be}{\bm{e}}

\newcommand{\bv}{\bm{v}}
\newcommand{\p}{\partial}
\renewcommand{\div}{{\rm{div}\,}}

\newcommand{\grad}{\nabla}
\newcommand{\abs}[1]{\left\lvert #1 \right\rvert}
\newcommand{\norm}[1]{\left\lVert #1 \right\rVert}
\newcommand{\wh}[1]{\widehat{#1}}
\newcommand{\wt}[1]{\widetilde{#1}}
\newcommand{\mc}[1]{\mathcal{#1}}
\newcommand{\mb}{\bm}

\newtheorem{theorem}{Theorem}[section]
\newtheorem{lemma}[theorem]{Lemma}
\newtheorem{proposition}[theorem]{Proposition}
\newtheorem{corollary}[theorem]{Corollary}
\newtheorem{definition}[theorem]{Definition}
\theoremstyle{definition}

\allowdisplaybreaks

\numberwithin{equation}{section}

\begin{document}
\title{A hierarchy of blood vessel models, Part II: 3D-3D to 3D-1D and 1D}

\author{Laurel Ohm}
\address{Department of Mathematics, University of Wisconsin - Madison, Madison, WI 53706}
\email{lohm2@wisc.edu}

\author{Sarah Strikwerda}
\address{Department of Mathematics, University of Wisconsin - Madison, Madison, WI 53706}
\email{sstrikwerda@wisc.edu}

\begin{abstract}
We propose and analyze a hierarchy of three models of blood perfusion through a tissue surrounding a thin arteriole or venule. Our goal is to rigorously link 3D-3D Darcy--Stokes, 3D-1D Darcy--Poiseuille, and 1D Green's function methods commonly used to model this process. Here in Part II, we consider the most detailed level, a 3D-3D Darcy-Stokes system coupled across the permeable vessel surface by mass conservation and pressure/stress balance conditions. We derive a convergence result between the 3D-3D model and both the 3D-1D Darcy--Poiseuille model and 1D Green's function model proposed in Part I \cite{partI} at a rate proportional to $\epsilon^{1/6}\abs{\log\epsilon}$, where $\epsilon$ is the maximum vessel radius. The rate is limited by the inclusion of a degenerate endpoint where the vessel radius vanishes, i.e. becomes indistinguishable from a capillary. Key to our proof are \emph{a priori} estimates for the 1D integrodifferential model obtained in Part I.
\end{abstract}

\maketitle

\tableofcontents

\input{introduction}
\input{existence3D3D}
\input{asymptotics}

\input{existence3D1D}
\input{errorest}

\appendix
\numberwithin{theorem}{section}
\input{appendix}


\subsubsection*{Acknowledgments} We thank Yoichiro Mori for guidance in formulating the 3D-3D and 3D-1D models and for helpful initial discussions. LO acknowledges support from NSF grant DMS-2406003 and from the VCRGE with funding from the Wisconsin Alumni Research Foundation. SS acknowledges support from NSF grant DMS-2037851.


\bibliographystyle{abbrv}
\bibliography{BloodBib.bib}


\end{document}

%% file: introduction.tex

\section{Introduction}
We propose and analyze a hierarchy of three models of blood perfusion through a tissue surrounding a thin blood vessel. Our goal is to rigorously establish links among 3D-3D Darcy--Stokes, 3D-1D Darcy--Poiseuille, and 1D Green's function methods commonly used to model perfusion \cite{d2011robust,secomb2004green,hsu1989green,pries2008blood,qi2021control,nobile2009coupling,d2007multiscale,notaro2016mixed,blanco2009potentialities,blanco2007unified,d2008coupling}. Our models are at the level of microvasculature: we consider a single arteriole or venule embedded in the tissue which, respectively, delivers or takes blood to/from the surrounding dense capillary bed from/to a larger artery/vein connecting back to the heart. From the capillaries, vital nutrients and oxygen from the blood pass through the capillary walls into the tissue. The location and geometries of arterioles and venules are important in determining large scale properties of blood flow within tissue \cite{pittman1995influence, cassot2010branching}, with arterioles serving as ``bottlenecks" to flow \cite{nishimura2007penetrating,shih2015robust} and with geometric attributes such as tortuosity correlating with disease \cite{johnston2009tortuosity, lorthois2014tortuosity}.
To probe effects of possibly complicated vessel geometries on perfusion, it is useful to derive models which strike a balance between physical relevance and simplicity.

Here we are concerned primarily with the most detailed level of the modeling hierarchy, a 3D-3D Darcy--Stokes system, and with showing how it relates to the other two levels of the hierarchy. The relationship between the 3D-1D Darcy--Poiseuille model and the 1D Green's function model is treated in Part I \cite{partI}. 
At all levels, as in \cite{koch2020modeling, d2007multiscale, qi2021control, qi2024hemodynamic, hyde2013parameterisation, shipley2020hybrid, sweeney2024three}, we take a continuum approach to modeling the capillary network surrounding the vessel, treating it as a porous medium with isotropic permeability. We consider this porous medium to occupy the upper half space $\R^3_+$, which is a convenience that allows us access to an explicit Green's function for building the 1D model. We will rely on this explicit information to show convergence among all three models.
The blood vessel itself is taken to have one end embedded in the tissue boundary at $z=0$ and one end free in the porous medium, and may be highly curved. The maximum radius of the vessel is a small parameter $\epsilon>0$, and the non-constant radius decays to zero at the free end, i.e., the vessel becomes indistinguishable from a capillary. We model the blood flow within the vessel at two levels of detail: first, as 3D Stokes flow, and second, due to the slenderness of the vessel, as a 1D Poiseuille approximation. 
In both descriptions, the vessel wall is treated as a semi-permeable membrane, and blood flow within the vessel is coupled to the exterior pressure by mass conservation and pressure/stress balance conditions across the interface. For the 3D Darcy--3D Stokes coupling, we take no-slip conditions for the tangential component of the Stokes fluid within the vessel, as in \cite{discacciati2009navier, d2011robust, zunino2002mathematical}.   
For the 3D Darcy--1D Poiseuille coupling, we propose a geometrically constrained Robin boundary condition for the exterior pressure coupled with a novel angle-averaged Neumann boundary condition, similar to boundary conditions arising in a PDE justification of slender body theories for impermeable filaments in Stokes flow \cite{closed_loop, free_ends,rigid,inverse,laplace}.

After showing well-posedness for the 3D-3D and 3D-1D models, we prove convergence of the pressure fields from the 3D-3D Darcy--Stokes system to the pressure fields of the 3D-1D Darcy--Poiseuille system as $\epsilon\to 0$. The convergence rate of $\epsilon^{1/6}\abs{\log\epsilon}$ is limited by the free end, but holds up to the endpoint. 
To obtain this result, it is necessary to pass through the 1D model, as we require further asymptotics for the pressure exterior to the vessel which make use of the vessel slenderness. Thus, our actual convergence result is direct from 3D-3D to 1D, and we obtain the 3D-3D to 3D-1D convergence as a consequence of Part I \cite{partI}. Our 3D-3D to 1D convergence result, at the same rate of $\epsilon^{1/6}\abs{\log\epsilon}$, is more detailed, as we construct a velocity ansatz from the 1D pressure and use it to additionally obtain an approximation of the blood velocity field within the vessel. 
The 3D-1D model is still a vital component of the hierarchy: it provides a roadmap for the asymptotics of section \ref{sec:asymptotics}, and is necessary for motivating the 1D model. Furthermore, related 3D-1D models with different choices of coupling are commonly used in both analytical and computational studies of flow in and around thin blood vessels \cite{berrone2023optimization, koch2020modeling, d2007multiscale, d2012finite, d2008coupling, blanco2007unified, formaggia2007stability, koppl20203d, nobile2009coupling, formaggia2006coupling, blanco2009potentialities, notaro2016mixed, possenti2021mesoscale, laurino2019derivation, kuchta2021analysis, boulakia2024mathematical}.
We highlight in particular a related error analysis starting from a coupled 3D Darcy--3D Darcy system in \cite{laurino2019derivation}, and the general framework for analyzing 3D-3D to 3D-1D modeling error proposed in \cite{heltai2023reduced}. However, to our knowledge, our results are the first to establish convergence among all three modeling levels, and we again stress the novelty of the coupling used in our 3D-1D model. Our approach treats the vessel as a truly 3D object even in the reduced 3D-1D and 1D models, which offers a different perspective on coupling and convergence than many existing methods which conceive of the flow within the reduced vessel as a Dirac mass.

We emphasize that the upper half space is used due to our reliance on an explicit Green's function to formulate and study the 1D model, which ends up being vital to the convergence analysis. This is in many ways a modeling convenience, but it means that we need to consider a free end in order to allow the centerline geometry to be as general as possible. We could instead consider only vessels which reconnect back to the tissue boundary at $z=0$. This restrictive assumption would only simplify the arguments here. As mentioned, the free end should be understood as the vessel radius becoming so small that it is indistinguishable from the capillaries making up the surrounding porous medium. The presence of a free end with decaying radius poses analytical difficulties. Both the 3D-1D and 1D models are degenerate at the free end due to the decay in the radius and do not admit a well-defined boundary value. The convergence rate of $\epsilon^{1/6}\abs{\log\epsilon}$ is due entirely to the endpoint behavior; otherwise, we would be limited primarily by the $\epsilon^{1/2}\abs{\log\epsilon}$ convergence rate obtained in Part I \cite{partI} between the 3D-1D and 1D models. We could incorporate more refined asymptotics for the boundary layer at the vessel end, and note that related outlet boundary considerations are typical of pipe flows \cite{ghosh2021modified,canic2003effective,castineira2019rigorous}. More refined asymptotics would result in more accurate but more complicated 3D-1D and 1D models at the vessel free end. Here we opt for simplicity of the reduced models over overall accuracy.

Furthermore, a perhaps more physically realistic setting would be a bounded domain with a vessel connecting to the tissue boundary at opposite ends, with the flow within and outside the vessel driven by the pressure drop between the two ends. In this setting, the 3D-3D and 3D-1D models themselves are straightforward to formulate and analyze on their own, but convergence is less straightforward, since we rely on 1D asymptotics for the exterior pressure field. We expect an analogous convergence result between the 3D-3D and 3D-1D models to be true, but much more technical to actually show without access to an explicit 1D model for comparison. Obtaining 3D-3D to 3D-1D convergence in a bounded exterior domain is a subject for future work.

\subsection{Setup}\label{sec:setup}
\begin{figure}[!ht]
\centering
\includegraphics[scale=0.4]{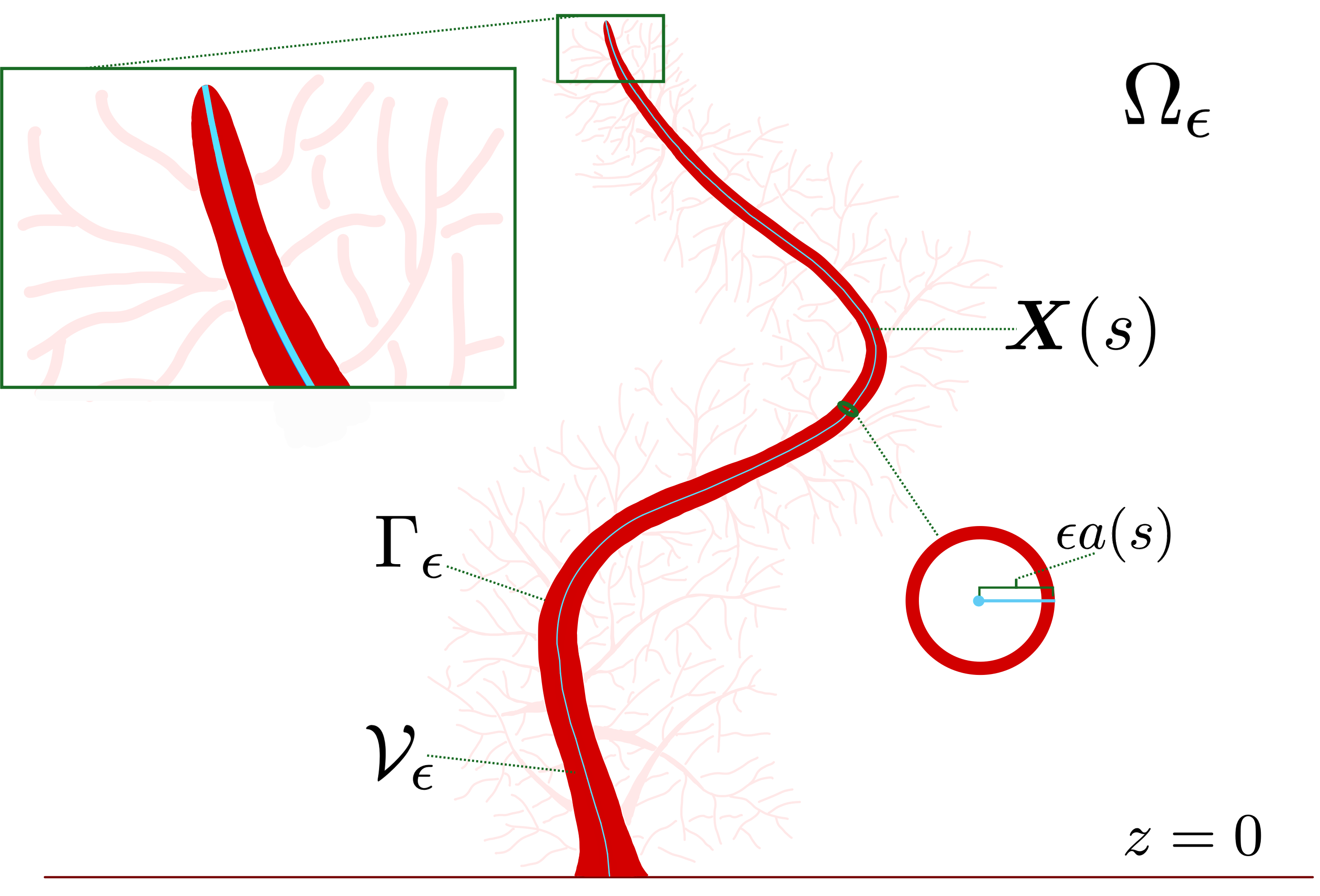}
\caption{An example geometry for the blood vessel $\mc{V}_\epsilon$.}
\label{fig:vessel1}
\end{figure}

 The dense capillary network exterior to the blood vessel, modeled as a porous medium, is assumed to occupy the entire upper half space $\R^3_+$ with a boundary at $z=0$. We consider the centerline of the vessel as a unit-length $C^3$ curve $\Gamma_0$ with one end embedded in the $z=0$ plane and the other end free in the porous medium. Let $\X:[0,1]\to \R^3_+$ denote the arclength parameterization of $\Gamma_0$ with parameter $s$ (see figure \ref{fig:vessel1}). For convenience, we will require that the tangent vector $\X_s := \frac{d\X}{ds}$ is perpendicular to the $xy$-plane when $s=0$. We further require that the curve does not self-intersect and, except for the endpoints at $s=0$, does not intersect the domain walls. In particular, 
\begin{equation}\label{eq:cGamma}
\begin{aligned}
    c_{\Gamma} = \min&\bigg\{\inf_{s_1\neq s_2}\frac{\abs{\X(s_1)-\X(s_2)}}{\abs{s_1-s_2}}\,, \; 
    \inf_{0<s\le 1}\frac{{\rm dist}(\X(s),\{z=0\})}{s}
    \bigg\} >0\,.
\end{aligned}    
\end{equation} 

We will allow for a variable radius along the length of the vessel, and define an admissible radius function as follows. 
\begin{definition}[Admissible radius]\label{def:rad}
    An admissible radius function $a:[0,1]\to [0,1]$ satisfies $a(s)\in C^2[0,1)$ with 
  \begin{equation}\label{eq:astar}
  a_\star=\sup_{s\in[0,1)}\abs{a(s)a'(s)} <\infty\,, \quad a_{\star\star}=\sup_{s\in[0,1)}\abs{a^3(s)a''(s)} <\infty\,.
  \end{equation}
  In addition, we require $\norm{a}_{L^\infty}=1$ and
  \begin{equation}\label{eq:delta}
  a(s)\ge a_0>0 \qquad \text{for }0\le s\le 1-\delta
  \end{equation}
  for some $\delta>0$, with spheroidal decay as $s\to 1$. In particular, we require 
    \begin{equation}\label{eq:spheroidal}
  \abs{a(s)-\sqrt{1-s^2}} \le C\epsilon^2\sqrt{1-s^2}\,.
  \end{equation}
  for $s>1-\delta$ where $\delta$ is as above. 
\end{definition}

Given $0<\epsilon\ll1$ and $a(s)$ as in Definition \ref{def:rad}, and letting $\X_s(s)^\perp$ denote the plane perpendicular to $\X(s)$ at cross section $s$, we may define the blood vessel $\mc{V}_\epsilon$ and its surface $\Gamma_\epsilon$ as
\begin{align}
    \mc{V}_\epsilon &= \big\{ \bx\in \X_s(s)^\perp\,:\, {\rm dist}(\bx,\X(s))<\epsilon a(s)\,, \; 0< s< 1 \big\} \,, \label{eq:Veps} \\
    \Gamma_\epsilon &= \big\{ \bx\in \X_s(s)^\perp\,:\, {\rm dist}(\bx,\X(s)) = \epsilon a(s)\,, \; 0< s< 1 \big\}\,. \label{eq:Gameps}
\end{align}
Here we take $\epsilon$ small enough relative to the vessel centerline curvature such that each $\bx\in \mc{V}_\epsilon$ belongs to a unique cross section, and aside from the $s=0$ cross section, we have $\mc{V}_\epsilon\cap \{z=0\}=\varnothing$. We denote the bottom cross section of $\mc{V}_\epsilon$ by
\begin{equation}\label{eq:Beps}
    B_\epsilon = \overline{\mc{V}_\epsilon}\cap \{z=0\}\,,
\end{equation}
and denote the domain exterior to the vessel as
\begin{equation}
    \Omega_\epsilon=\R^3_+\backslash\overline{\mc{V}_\epsilon}\,.
\end{equation}

We will often think of the surface $\Gamma_\epsilon$ as a function of arclength $s$ and the angle $\theta$ over each cross section (see figure \ref{fig:vessel1}). 
We then write the surface element $dS_\epsilon$ along $\Gamma_\epsilon$ as
\begin{equation}
    dS_\epsilon = \mc{J}_\epsilon(s,\theta)\,d\theta\,ds\,,
\end{equation}
where $\mc{J}_\epsilon$ is a Jacobian factor that we parameterize explicitly when needed.

\subsection{The models}
We begin with the 3D Darcy--3D Stokes description of perfusion outside of the vessel and flow within. 
Within the vessel, we consider velocity and pressure fields $(\bu^\epsilon,p^\epsilon):\mc{V}_\epsilon\times \mc{V}_\epsilon\to \R^3\times \R$. Letting $\E(\bv) = \frac{1}{2}(\nabla\bv+\nabla\bv^{\rm T})$ denote the symmetric gradient, we also define the stress tensor $\Sigma^\epsilon = 2\mu\,\E(\bu^\epsilon)-p^\epsilon{\bf I}$, where $\mu$ is the blood viscosity and is independent of $\epsilon$. We denote the exterior pressure field by $q^\epsilon:\Omega_\epsilon\to \R$. We will typically use the superscript $(\cdot)^\epsilon$ to refer to quantities satisfying the 3D-3D system. 
The 3D-3D Darcy--Stokes system about the vessel is given by 
\begin{subequations}
\begin{align}
-\mu\Delta\bu^\epsilon +\nabla p^\epsilon &= 0\,, \quad \div \bu^\epsilon=0 \hspace{3.7cm} \text{in } \mc{V}_\epsilon \label{eq:3D3D_1} \\
\zeta\epsilon^4\Delta q^\epsilon &= 0 \hspace{6.05cm} \text{in }\Omega_\epsilon \label{eq:3D3D_2}\\
\bu^\epsilon\cdot\bm{n} &= -\zeta\epsilon^4\frac{\p q^\epsilon}{\p\bm{n}}= -\kappa\epsilon^3\big((\Sigma^\epsilon\bm{n})\cdot\bm{n}+q^\epsilon\big) \qquad \text{on }\Gamma_\epsilon \label{eq:3D3D_3}\\
\bu^\epsilon-(\bu^\epsilon\cdot\bm{n})\bm{n}&=0 \hspace{6.05cm} \text{on }\Gamma_\epsilon \label{eq:3D3D_4}\\
\frac{\p q^\epsilon}{\p\bm{n}} &= 0 \quad \text{on }\p\Omega_\epsilon\backslash\Gamma_\epsilon\,, \quad q^\epsilon\to 0 \;\text{ as } \abs{\bx}\to\infty \label{eq:3D3D_5}\\
\Sigma^\epsilon\bm{n}\big|_{B_\epsilon} &= -p_0\bm{n} \,,\label{eq:3D3D_6}
\end{align}
\end{subequations}
where $p_0\in \R$ is incoming pressure data. Throughout, we take the normal vector $\bm{n}$ to point out of the vessel $\mc{V}_\epsilon$ and into the domain $\Omega_\epsilon$. 
Here the Darcy flow through the porous medium outside of the vessel has hydraulic conductivity $\zeta\epsilon^4$, where $\zeta$ is an $O(1)$ constant. The conductivity comes into play in the first part of boundary condition \eqref{eq:3D3D_3} through the conservation of mass across the interface $\Gamma_\epsilon$. 
We may think of the $\epsilon$-dependent scaling here as a small Darcy number scaling: in particular, the ratio between the permeability of the exterior medium and the vessel cross sectional area satisfies $\frac{(\zeta\epsilon^4)\mu}{\epsilon^2}=O(\epsilon^2)$. 
The fluid velocity across $\Gamma_\epsilon$ is additionally balanced by the pressure difference between the interior and exterior of the vessel, where $\kappa\epsilon^3$ corresponds to the conductivity of the vessel itself. Again $\kappa$ is an $O(1)$ constant. The choice of $\epsilon$-dependent scaling here arises due to energy considerations in comparing the 3D-3D system to the proposed 3D-1D reduction. In particular, the reduction holds provided that the vessel walls and surrounding capillary bed are relatively impermeable compared to the vessel interior.

We next turn to the reduced 3D Darcy--1D Poiseuille description of perfusion about the vessel. Here we denote the interior pressure by  $p:[0,1]\to \R$ and exterior pressure by $q:\Omega_\epsilon \to \R$. Throughout, we refer to quantities satisfying the 3D-1D system without superscripts. The 3D-1D Darcy--Poiseuille system is given by
\begin{subequations}
\begin{align}
\Delta q &= 0 \hspace{4cm} \text{in }\Omega_\epsilon  \label{eq:3D1D_1} \\
-\frac{\p q}{\p \bm{n}} &= \frac{\kappa}{\zeta\epsilon} (p(s)-q) \hspace{2.15cm} \text{on }\Gamma_\epsilon  \label{eq:3D1D_2}\\
\int_0^{2\pi}\frac{\p q}{\p \bm{n}}\,\mc{J}_\epsilon(s,\theta)\,d\theta &= -\frac{\pi}{8\zeta\mu}\frac{d}{ds}\bigg(a^4(s)\frac{d p}{ds}\bigg) \qquad \text{on }\Gamma_\epsilon  \label{eq:3D1D_3} \\
\frac{\p q}{\p\bm{n}} &= 0 \quad \text{on }\p\Omega_\epsilon\backslash\Gamma_\epsilon\,, \quad q\to 0 \;\text{ as }\abs{\bx}\to\infty \label{eq:3D1D_4} 
\end{align}
\end{subequations}
with prescribed incoming pressure data $p(0)=p_0\in \R$. Note that there is no boundary condition for $p$ at $s=1$ due to the degeneracy of the coefficient $a^4(s)$ in \eqref{eq:3D1D_3} as $s\to 1$. Here again $\bm{n}$ points out of the vessel $\mc{V}_\epsilon$ and into the domain $\Omega_\epsilon$.

The boundary condition \eqref{eq:3D1D_2} corresponds directly with the second part of \eqref{eq:3D3D_3}: the flux of the exterior pressure $q$ across $\Gamma_\epsilon$ is balanced by the pressure difference between the interior and exterior. The $\theta$-averaged boundary condition \eqref{eq:3D1D_4} is less obvious. The right hand side, including the $a^4$ dependence on the radius function, is Poiseuille's law for pressure-driven flow in a thin domain, and we note that when $\zeta=0$ (corresponding to an impermeable vessel), we recover usual Poiseuille flow in the vessel. The 1D nature of the interior pressure $p(s)$ due to the Poiseuille asymptotics imposes a geometric constraint on the Robin boundary value of $q$ via \eqref{eq:3D1D_2}. Because of this geometric constraint on $q$, we have only partial degrees of freedom left to determine via the boundary condition \eqref{eq:3D1D_3}, resulting in the $\theta$-averaged condition. It will also be convenient to combine the boundary conditions \eqref{eq:3D1D_2} and \eqref{eq:3D1D_3} to rewrite \eqref{eq:3D1D_3} as 
\begin{equation}\label{eq:3D1D_3alt}
    \frac{\pi}{8\kappa\mu}\frac{d}{ds}\bigg(a^4(s)\frac{dp}{ds}\bigg) = \frac{j_\epsilon(s)}{\epsilon}p(s) -\frac{1}{\epsilon}\int_0^{2\pi}q\,\mc{J}_\epsilon(s,\theta)\,d\theta \,,
\end{equation}
where $j_\epsilon(s) = \int_0^{2\pi}\mc{J}_\epsilon(s,\theta)\,d\theta\approx 2\pi \epsilon a(s)$.

We obtain the 3D-1D system \eqref{eq:3D1D_1}-\eqref{eq:3D1D_4} from the 3D-3D system \eqref{eq:3D3D_1}-\eqref{eq:3D3D_6} via formal asymptotics in section \ref{sec:asymptotics}, and ultimately prove a convergence result between the two systems in section \ref{sec:error_estimate}.

\subsection{Theorem statements}
We begin with well-posedness and quantitative energy estimates for the 3D-3D and 3D-1D systems. For both models, since the exterior domain $\Omega_\epsilon$ is unbounded, it will be natural consider the exterior pressures $q^\epsilon$ and $q$, respectively, belonging to the function space 
\begin{equation}
  D^{1,2}(\Omega_\epsilon) := \{ g \in L^6(\Omega_\epsilon)\,:\, \nabla g\in L^2(\Omega_\epsilon) \}
\end{equation}
with norm $\norm{g}_{D^{1,2}(\Omega_\epsilon)}:= \norm{\nabla g}_{L^2(\Omega_\epsilon)}$ by the 3D Sobolev inequality. 
We show the following well-posedness result for the 3D-3D system \eqref{eq:3D3D_1}-\eqref{eq:3D3D_6}. 
\begin{theorem}[Well-posedness for 3D-3D system]\label{thm:3D3D}
    Given $\mc{V}_\epsilon$, $\Omega_\epsilon$ as in section \ref{sec:setup} and incoming pressure data $p_0\in \R$, there exists a unique solution $(\bu^\epsilon,p^\epsilon,q^\epsilon)\in H^1(\mc{V}_\epsilon)\times L^2(\mc{V}_\epsilon)\times D^{1,2}(\Omega_\epsilon)$ to the system \eqref{eq:3D3D_1}-\eqref{eq:3D3D_6} satisfying the bound
    \begin{equation}\label{eq:3d3d_energy}
    \begin{aligned}
        &\frac{1}{|\mc{V}_\epsilon|^{1/2}}\bigg(\norm{\epsilon^{-2}\bu^\epsilon}_{L^2(\mc{V}_\epsilon)} + \norm{\epsilon^{-1}\nabla\bu^\epsilon}_{L^2(\mc{V}_\epsilon)} + \norm{p^\epsilon}_{L^2(\mc{V}_\epsilon)}\bigg)+ \norm{q^\epsilon}_{D^{1,2}(\Omega_\epsilon)}\\
        &+ \frac{1}{|\Gamma_\epsilon|^{1/2}} 
        \bigg(\norm{\epsilon^{-3}(\bu^\epsilon\cdot\bm{n})}_{L^2(\Gamma_\epsilon)}+ \epsilon\norm{\frac{\p q^\epsilon}{\p\bm{n}}}_{L^2(\Gamma_\epsilon)} + \norm{(\Sigma^\epsilon\bm{n})\cdot\bm{n}+q^\epsilon}_{L^2(\Gamma_\epsilon)}\bigg)
        \le C\abs{p_0}\,.
    \end{aligned}    
    \end{equation}
    Here $C$ depends on the constants $\zeta$, $\mu$, and $\kappa$ but not on $\epsilon$.
\end{theorem}
The proof of Theorem \ref{thm:3D3D} appears in section \ref{sec:3D3Dthm}, including justification for the choice of $\epsilon$-scaling of the coefficients $\zeta\epsilon^4$ and $\kappa\epsilon^3$ appearing in \eqref{eq:3D1D_3}.

We next consider well-posedness for the reduced 3D-1D system \eqref{eq:3D1D_1}-\eqref{eq:3D1D_4}. Due to the decay of the vessel radius function $a(s)$ as $s\to 1$, it is natural to consider the interior pressure $p(s)$ belonging to the weighted function space
\begin{equation}\label{eq:Ha_def}
\begin{aligned}
  \mc{H}^a(0,1) &= \big\{ v\in L^2(0,1)\,:\, a^2v_s \in L^2(0,1) \big\}\,, \\
  \norm{v}_{\mc{H}^a(0,1)}^2 &:= \norm{v}_{L^2(0,1)}^2+\norm{a^2v_s}_{L^2(0,1)}^2\,.
\end{aligned}
\end{equation}
We show the following.
\begin{theorem}[Well-posedness for 3D-1D system]\label{thm:3D1D}
     Given $\mc{V}_\epsilon$, $\Omega_\epsilon$ as in section \ref{sec:setup} and incoming pressure data $p_0\in \R$, there exists a unique solution $(p,q)\in \mc{H}^a(0,1)\times D^{1,2}(\Omega_\epsilon)$ to \eqref{eq:3D1D_1}-\eqref{eq:3D1D_4} satisfying the bound
    \begin{equation}\label{eq:energy3D1D}
\norm{p}_{\mc{H}^a(0,1)} + \frac{1}{|\Gamma_\epsilon|^{1/2}}\norm{p-q}_{L^2(\Gamma_\epsilon)} + \norm{q}_{D^{1,2}(\Omega_\epsilon)}
        \le C\abs{p_0}\,.
    \end{equation}
    Here $C$ depends on the parameters $\zeta$, $\mu$, and $\kappa$ but not on $\epsilon$.   
\end{theorem}
The proof of Theorem \ref{thm:3D1D} appears in section \ref{sec:3D1Dthm}. Note the correspondence between the energy estimates \eqref{eq:3d3d_energy} and \eqref{eq:energy3D1D} for the 3D-3D and 3D-1D systems. 
Between the well-posedness proofs, section \ref{sec:asymptotics} contains formal asymptotics which demonstrate the connection between the 3D-3D and 3D-1D models.

Our main result is a quantitative bound on the rate at which the interior and exterior blood pressures $(p^\epsilon,q^\epsilon)$ from the 3D-3D model approach the 3D-1D solution $(p,q)$ as the maximum vessel radius $\epsilon$ tends to zero. We note that the rate we obtain here is valid out to the vessel tip at $s=1$, despite the fact that we do not include a more detailed asymptotic description of the boundary layer at the free end.
We show the following.
\begin{theorem}[3D-3D to 3D-1D error estimate]\label{thm:error}
Given a vessel $\mc{V}_\epsilon$ as in section \ref{sec:setup} and pressure data $p_0\in \R$, let $(\bu^\epsilon,p^\epsilon,q^\epsilon)$ be the solution to the 3D-3D system \eqref{eq:3D3D_1}-\eqref{eq:3D3D_6} and let $(p,q)$ denote the solution to the 3D-1D system \eqref{eq:3D1D_1}-\eqref{eq:3D1D_4}. For $\epsilon$ sufficiently small, the differences $(p^\epsilon-p,q^\epsilon-q)$ satisfy
\begin{equation}
    \frac{1}{|\mc{V}_\epsilon|^{1/2}}\norm{p^\epsilon-p}_{L^2(\mc{V}_\epsilon)} + \norm{q^\epsilon-q}_{D^{1,2}(\Omega_\epsilon)}
    \le C\epsilon^{1/6}\abs{\log\epsilon}\abs{p_0}\,.
\end{equation}
\end{theorem}

We rely on a somewhat circuitous strategy to obtain this result. In particular, the $\theta$-averaged boundary condition \eqref{eq:3D1D_3} in the 3D-1D model is especially tricky to work with, as it makes it difficult to obtain necessary pointwise information about the exterior pressure $q$ on $\Gamma_\epsilon$. Luckily, the 3D-1D system itself admits a lower dimensional slender-body-type approximation, introduced below, which lends itself better to detailed analysis. In Part I \cite{partI} of our program, we analyze this slender body model in detail, and show that its solution converges to the 3D-1D solution with a rate $O(\epsilon^{1/2}\abs{\log\epsilon})$ (see Theorem \ref{thm:3D1Dto1D}). We rely on this convergence result, as well as some favorable analytical properties of the 1D model, to obtain Theorem \ref{thm:error}. We detail our strategy in the next section.

\subsection{1D approximation}
The formal asymptotics of section \ref{sec:asymptotics} suggest the 3D-1D model \eqref{eq:3D1D_1}-\eqref{eq:3D1D_4} as a reduced order approximation of the 3D-3D system \eqref{eq:3D3D_1}-\eqref{eq:3D3D_6} for small vessel radius. However, to actually obtain the quantitative error bound of Theorem \ref{thm:error}, it turns out that we will require a more detailed asymptotic characterization of the 3D-1D solution as well. For this, we turn to a fully 1D ``slender body"-type approximation of the 3D-1D system, which is introduced and studied in depth in Part I of our program \cite{partI}.

Let $\mc{G}_N(\bx,\by)$ denote the Green's function for the Neumann problem for the Laplace equation in $\R^3_+$, given by
\begin{equation}\label{eq:GN}
\mc{G}_N(\bx,\by) = \frac{1}{4\pi}\bigg(\frac{1}{\abs{\bx-\by}}+\frac{1}{\abs{\bx-\by^*}} \bigg)\,.
\end{equation}
Here $\by^*\in \R^3_{-}$ denotes the reflection of $\by\in \R^3_+$ across $y_3=0$, i.e. for $\by=(y_1,y_2,y_3)$, $\by^* = (y_1,y_2,-y_3)$.
Given $f:[0,1]\to \R$, for $\bx\in \R^3_+$, we may then define the operator
\begin{equation}\label{eq:SN_op}
\mc{S}_N[f](\bx) = \int_0^1\mc{G}_N\big(\bx,\X(\sqrt{1-\epsilon^2}\,t)\big)\,f(t)\,dt\,.
\end{equation}
Note that, by design, the kernel in \eqref{eq:SN_op} does not quite extend up to the tip of the vessel at $\X(1)$, which ensures that $\mc{S}_N$ is well defined everywhere along the blood vessel surface $\Gamma_\epsilon$, up to and including the tip.

We then define the ``slender body" approximation $q^{\rm SB}(\bx)$ to the exterior pressure $q(\bx)$ by the explicit expression 
\begin{equation}\label{eq:qSB}
q^{\rm SB}(\bx) = \frac{\pi}{8\zeta\mu}\mc{S}_N\bigg[\frac{d}{ds}\bigg(a^4\frac{dp^{\rm SB}}{ds}\bigg)\bigg](\bx)\,, \qquad \bx\in \Omega_\epsilon\,.
\end{equation}
Here the interior pressure $p^{\rm SB}(s)$ comes from solving the 1D integrodifferential equation
\begin{equation}\label{eq:pSB}
\begin{aligned}
\frac{d}{ds}\bigg(a^4\frac{dp^{\rm SB}}{ds}\bigg)
&=  \frac{8\mu\kappa}{\pi}a(s)\int_0^{2\pi} \big( p^{\rm SB}(s) - q^{\rm SB}\big|_{\Gamma_\epsilon}\big)\, d\theta \\
&= 16\mu\kappa a(s)\,p^{\rm SB}(s)- \frac{8\mu\kappa}{\pi}a(s)\int_0^{2\pi} \mc{S}_N\bigg[\frac{d}{ds}\bigg(a^4\frac{dp^{\rm SB}}{ds}\bigg)\bigg]\bigg|_{\Gamma_\epsilon}\,d\theta
\end{aligned}
\end{equation}
along the length of the vessel, with the boundary condition $p^{\rm SB}(0)=p_0$. As in the 3D-1D system, the degenerate $a^4$ coefficient precludes a well-defined boundary condition for $p^{\rm SB}$ at $s=1$.

In Part I \cite{partI}, we analyze the equation \eqref{eq:pSB} in detail and prove the following convergence result:
\begin{theorem}[3D-1D to 1D convergence \cite{partI}]\label{thm:3D1Dto1D}
Given $\mc{V}_\epsilon$ as in section \ref{sec:setup} and incoming pressure data $p_0\in\R$, consider the solution $(p(s),q)\in \mc{H}^a(0,1)\times D^{1,2}(\Omega_\epsilon)$ to the 3D-1D model \eqref{eq:3D1D_1}-\eqref{eq:3D1D_4} as well as the 1D approximation $(p^{\rm SB}(s),q^{\rm SB})$ given by \eqref{eq:qSB}, \eqref{eq:pSB}. The difference $(p^{\rm SB}-p,q^{\rm SB}-q)$ satisfies
\begin{equation}
\begin{aligned}
\norm{p^{\rm SB}-p}_{\mc{H}^a(0,1)} + &\epsilon^{-1/2}\norm{(p^{\rm SB}-p)-(q^{\rm SB}-q)}_{L^2(\Gamma_\epsilon)} + \norm{q^{\rm SB}-q}_{D^{1,2}(\Omega_\epsilon)} \\
&\le C \epsilon^{1/2}\abs{\log\epsilon}\,\abs{p_0}
\end{aligned}
\end{equation}
for $C$ independent of $\epsilon$ as $\epsilon\to 0$.
\end{theorem}

Part I of our program is devoted to the proof of Theorem \ref{thm:3D1Dto1D}, which is stated in \cite{partI} as Theorem 1.4. Here we will rely on a few important features of the 1D model that make it more convenient than the 3D-1D model for comparison with the full 3D-3D system. First, we require a near-$\theta$-independence result for the reduced order approximant of the exterior pressure $q^\epsilon$ along the vessel surface $\Gamma_\epsilon$. This type of $\theta$-independence condition is typical of slender body theories \cite{closed_loop, free_ends,rigid,inverse,laplace}. A difficulty here is that the \emph{Robin} boundary value \eqref{eq:3D1D_2} for $q$ is explicitly $\theta$-independent, but it is not straightforward to extract enough smallness (in $\epsilon$) for $q$ itself to close an error estimate using the 3D-1D description alone. However, the small $\theta$-dependence of the 1D approximant $q^{\rm SB}$ on $\Gamma_\epsilon$ can be calculated directly. 

In addition, we will be relying on the geometric smallness of the vessel free end due to the spheroidal endpoint decay \eqref{eq:spheroidal} to obtain a convergence result that holds up to the free endpoint. To leverage this smallness, we need the reduced order approximant to the interior pressure $p^\epsilon$ to be sufficiently ``well behaved" out to the endpoint. \emph{A priori}, the 1D pressure approximant $p^{\rm SB}$ is slightly better behaved at the tip than the 3D-1D approximant $p$, which is also necessary for the error estimate to close.
The improvement is due to replacing the surface element $\mc{J}_\epsilon(s,\theta)$ in the 3D-1D boundary condition \eqref{eq:3D1D_3}, which looks like $\epsilon a(s) + C\epsilon^2$ at the tip, with $\epsilon a(s)$ in the 1D approximation \eqref{eq:pSB}, which decays.

Given the 3D-1D to 1D convergence, and the fact that the 1D model is more convenient for analysis, we will obtain Theorem \ref{thm:error} by directly showing convergence of the 3D-3D system to the 1D model. The formal asymptotics of section \ref{sec:asymptotics} for deriving the 3D-1D model suggest a suitable form of velocity ansatz built out of the 1D pressure $p^{\rm SB}(s)$, given by the variable-radius Poiseuille's law. 
In particular, using curvilinear cylindrical coordinates $\bx=\bx(r,\theta,s)$ defined with respect to the vessel centerline $\X(s)$ (see section \ref{subsec:geometry}), we define
\begin{equation}\label{eq:Udef0}
\bm{U}(\bx) = -\frac{1}{16\mu}\p_s\big( \phi_\epsilon(s)(r^3 -2(\epsilon a)^2r) p^{\rm SB}_s \big)\be_r(s,\theta)+ \frac{1}{4\mu}\phi_\epsilon(s)\big(r^2-(\epsilon a)^2\big)p^{\rm SB}_s\,\be_{\rm t}(s)\,.
\end{equation} 
Here $p^{\rm SB}(s)$ is the solution to \eqref{eq:pSB}, and $\phi_\epsilon(s)$ is a smooth cutoff function in $s$ satisfying 
\begin{equation}\label{eq:phieps_def0}
	\phi_\epsilon(s) = \begin{cases}
		1\,, & 0\le s\le 1-2\epsilon^{4/3}\\
		0 \,, & 1-\epsilon^{4/3}\le s\le 1\,,
	\end{cases}
	\qquad \abs{\p_s\phi_\epsilon}\le C\epsilon^{-4/3}\,.
\end{equation}
The choice of cutoff extent, which defines a ``tip region'' $s\in[1-2\epsilon^{4/3},1]$, is optimized to balance two competing effects: (1) The ansatz is derived using the assumption that the vessel is approximately cylindrical, which ceases to be true at the tip. (2) To counteract the fact that $\bm{U}$ is not a good approximation toward the tip, we will rely on the geometric smallness of the tip. For this to work, the length of the tip region has to be small enough to provide additional factors of $\epsilon$. The balance of these two effects in section \ref{sec:error_estimate} yields the $O(\epsilon^{4/3})$ tip extent, which in turn yields the $O(\epsilon^{1/6})$ convergence rate of Theorem \ref{thm:error}. This is in constrast to the ``full" $O(\epsilon^{1/2})$ rate suggested by the 3D-1D to 1D convergence result of Theorem \ref{thm:3D1Dto1D}.

Instead of showing Theorem \ref{thm:error} directly, we will use the triple $(\bm{U},p^{\rm SB},q^{\rm SB})$ built from the 1D approximation to obtain a convergence result from 3D-3D to 1D. This result, coupled with the 3D-1D to 1D convergence of Theorem \ref{thm:3D1Dto1D}, yields Theorem \ref{thm:error} as a byproduct. We thus aim to show the following. 
\begin{theorem}[3D-3D to 1D convergence]\label{thm:3D3Dto1D}
Given a vessel $\mc{V}_\epsilon$ as in section \ref{sec:setup} and incoming pressure data $p_0\in \R$, let $(\bu^\epsilon,p^\epsilon,q^\epsilon)$ be the solution to the 3D-3D system \eqref{eq:3D3D_1}-\eqref{eq:3D3D_6} and let $(\bm{U},p^{\rm SB},q^{\rm SB})$ be the 1D approximation \eqref{eq:qSB}-\eqref{eq:pSB} with velocity ansatz $\bm{U}$ as in \eqref{eq:Udef0}. For $\epsilon$ sufficiently small, the differences $(\bu^\epsilon-\bm{U},p^\epsilon-p^{\rm SB},q^\epsilon-q^{\rm SB})$ satisfy 
\begin{equation}
\begin{aligned}
&\frac{1}{|\mc{V}_\epsilon|^{1/2}}\bigg(\norm{\epsilon^{-2}(\bu^\epsilon-\bm{U})}_{L^2(\mc{V}_\epsilon)}+\norm{\epsilon^{-1}\nabla(\bu^\epsilon-\bm{U})}_{L^2(\mc{V}_\epsilon)} + \norm{p^\epsilon-p^{\rm SB}}_{L^2(\mc{V}_\epsilon)}\bigg)  \\
&\qquad  + \frac{1}{\epsilon^3|\Gamma_\epsilon|^{1/2}}\norm{(\bu^\epsilon-\bm{U})\cdot\bm{n}}_{L^2(\Gamma_\epsilon)} + \norm{q^\epsilon-q^{\rm SB}}_{D^{1,2}(\Omega_\epsilon)} \le C\epsilon^{1/6}\abs{\log\epsilon}\abs{p_0}
\end{aligned}
\end{equation}
for $C$ independent of $\epsilon$ as $\epsilon\to 0$.
\end{theorem}

The remainder of the paper is devoted to proving the well-posedness results of Theorems \ref{thm:3D3D} and \ref{thm:3D1D} in sections \ref{sec:3D3Dthm} and \ref{sec:3D1Dthm}, respectively, and ultimately to proving the 3D-3D to 1D convergence result of Theorem \ref{thm:3D1Dto1D} in section \ref{sec:error_estimate}.

%% file: existence3D3D.tex

\section{Proof of Theorem \ref{thm:3D3D}}\label{sec:3D3Dthm}

\subsection{Preliminaries}\label{subsec:geometry}
Throughout, it will be useful to refer to an explicit parameterization for the vessel $\mc{V}_\epsilon$. We define the unit tangent vector $\be_{\rm t}(s)=\X_s(s)$ to the vessel centerline, along with an initial choice of unit normal vector $\be_1(0)\perp\be_{\rm t}(0)$. We may then define a $C^2$ orthonormal Bishop frame $(\be_{\rm t}(s),\be_1(s),\be_2(s))$ about $\X(s)$ \cite{bishop1975there}, where
\begin{equation}
\frac{d}{ds}\begin{pmatrix}
\be_{\rm t}\\
\be_1\\
\be_2\\
\end{pmatrix}
= \begin{pmatrix}
0 & \kappa_1(s) &\kappa_2(s) \\
-\kappa_1(s) & 0 & 0\\
-\kappa_2(s) & 0 & 0
\end{pmatrix}\,,
\quad \kappa_1^2+\kappa_2^2 = \abs{\X_{ss}}^2\,.
\end{equation}
We denote the maximum curvature of $\X$ by  
\begin{equation}\label{eq:kappastar}
    \kappa_\star = \norm{\X_{ss}}_{L^\infty(0,1)}\,. 
\end{equation}
Bounds on the third derivative $\X_{sss}$ will be required only for a quantitative pressure estimate in Appendix \ref{subsec:pressure_estimate}.
Using the definition \eqref{eq:kappastar}, we may make the requirement stated below \eqref{eq:Gameps} more precise: throughout, we require the maximum vessel radius $\epsilon$ to satisfy $\epsilon<\frac{1}{8\kappa_\star}$ to ensure that cross sections are well-defined. 

It will be convenient to define a curved version of cylindrical coordinates in a neighborhood of $\X(s)$. Taking
\begin{equation}
\begin{aligned}
\be_r(s,\theta) &= \cos\theta\be_1(s) + \sin\theta\be_2(s)\,,\\
\be_\theta(s,\theta) &= -\sin\theta\be_1(s) + \cos\theta\be_2(s)\,,
\end{aligned}
\end{equation}
we may parameterize points $\bx$ within $\mc{V}_\epsilon$ and near $\Gamma_\epsilon$ as 
\begin{equation}\label{eq:param}
\bx = \bx(r,\theta,s) = \X(s) + r\be_r(s,\theta)\,.
\end{equation}
Within $\mc{V}_\epsilon$, we may parameterize the volume element $d\bx=\mc{J}(r,s,\theta)\,drd\theta ds$ as
\begin{equation}\label{eq:volume_element}
\begin{aligned}
    \mc{J}(r,s,\theta)&= r(1-r\wh\kappa(s,\theta))\,, \\
    \wh\kappa(s,\theta) &= \kappa_1(s)\cos\theta+\kappa_2(s)\sin\theta\,.
\end{aligned}
\end{equation}
In addition, following \cite[Appendix A.2]{rigid}, we may write the gradient and divergence in curvilinear cylindrical coordinates as 
\begin{equation}\label{eq:grad_div}
    \begin{aligned}
        \nabla \bv &= \frac{\p\bv}{\p r}\otimes\be_r + \frac{1}{r}\frac{\p\bv}{\p\theta}\otimes \be_\theta + \frac{1}{1-r\wh\kappa}\frac{\p\bv}{\p s}\otimes\be_{\rm t}  \\
        \div\bv &=\frac{1}{1-r\wh\kappa}\bigg(\frac{1}{r}\p_r\big(r(1-r\wh\kappa) v^r\big)+\frac{1}{r}\p_\theta\big((1-r\wh\kappa) v^\theta\big) +  \p_s v^s\bigg)\,,
    \end{aligned}
\end{equation}
where $v^r=\bv\cdot\be_r$, $v^\theta=\bv\cdot\be_\theta$, and $v^s=\bv\cdot\be_{\rm t}$.
Furthermore, along the vessel surface $\bx=\X(s)+\epsilon a(s)\be_r(s,\theta)\in \Gamma_\epsilon$, we may parameterize the outward-pointing unit normal vector to the vessel as 
\begin{equation}\label{eq:normal_param}
    \bm{n}(s,\theta) = \frac{\be_r(s,\theta)}{\sqrt{1+\epsilon^2a'(s)^2}} - \frac{\epsilon a'(s)\be_{\rm t}(s)}{\sqrt{1+\epsilon^2a'(s)^2}}\,.
\end{equation}
We may additionally parameterize the surface Jacobian $dS_\epsilon = \mc{J}_\epsilon(s,\theta)\,d\theta ds$ along $\Gamma_\epsilon$ as 
\begin{equation}\label{eq:free_jacfac}
\mc{J}_\epsilon(s,\theta) = \epsilon a(s) \sqrt{(1-\epsilon a(s)\wh\kappa(s,\theta))^2+\epsilon^2(a'(s))^2}\,,
\end{equation}
where $\wh\kappa$ is as in \eqref{eq:volume_element}.


\subsection{Existence and uniqueness for the 3D-3D system}

We begin by showing existence and uniqueness for fixed $\epsilon>0$. Throughout this section, we will treat the parameters appearing in the 3D-3D system \eqref{eq:3D3D_1}-\eqref{eq:3D3D_6} as unknown functions of $\epsilon$, which we will denote by $\mu^\epsilon$, $\kappa^\epsilon$, and $\zeta^\epsilon$. The $\epsilon$-dependent scaling used in \eqref{eq:3D3D_1}-\eqref{eq:3D3D_6} will be derived in section \ref{subsec:energy} when we consider the $\epsilon$-dependence of the solution.

We first introduce the relevant definitions of function spaces and weak solutions. Here we will use $(\cdot, \cdot)_D$ to denote the inner product on the space $D$ or the $L^2(D)$ inner product if $D$ is a domain.

We define $\mb{V}=\{\mb{w} \in H^1(\mc{V}_\epsilon)~:~ \mb{w}-(\mb{w} \cdot \mb{n})\mb{n}=0~\text{on}~\Gamma_\epsilon\}$. This space 
is a Hilbert space with the inner product 
\begin{equation}(\mb{u}, \mb{w})_{\mb{V}}=2\mu^\epsilon\big(\E (\mb{u}), \E (\mb{w})\big)_{\mc{V}_\epsilon} + \frac{1}{\kappa^\epsilon}(\mb{u} \cdot \mb{n}, \mb{w} \cdot \mb{n})_{\Gamma_\epsilon}\,.
\end{equation}
We define the bilinear form $\mathcal{A}:(\mb{V} \times D^{1,2}(\Omega_\epsilon))^2 \to \mathbb{R}$ as
\begin{equation}\label{eq:Aform0}
\mathcal{A}((\mb{u}, q), (\mb{w}, w))=(\mb{u}, \mb{w})_{\mb{V}}+(q, \mb{w} \cdot \mb{n})_{\Gamma_\epsilon} + (\zeta^\epsilon \nabla q, \nabla w)_{\Omega_\epsilon} - (\mb{u} \cdot \mb{n}, w)_{\Gamma_\epsilon}\,,
\end{equation}
and we let $b: (\mb{V}\times D^{1,2}(\Omega_\epsilon)) \times L^2(\mc{V}_\epsilon) \to \mathbb{R}$ be defined as 
\begin{equation}\label{eq:bform0}
b((\mb{v},v),p)=-(p, \nabla \cdot \mb{v})_{L^2(\mc{V}_\epsilon)}\,.
\end{equation}

We define a solution to the 3D-3D system in the following manner:
\begin{definition} \label{def:weak3D}
We will call $(\mb{u}^\epsilon,p^\epsilon,q^\epsilon)\in \mb{V} \times L^2(\mc{V}_\epsilon) \times D^{1,2}(\Omega_\epsilon)$ a weak solution to the 3D-3D system \eqref{eq:3D3D_1}-\eqref{eq:3D3D_6} if 
\begin{equation}
\begin{aligned}\label{3Dweakdef}
\mathcal{A}((\mb{u}^\epsilon, q^\epsilon), (\mb{w},w)) + b((\mb{w},w),p^\epsilon)&=-(p_0\mb{n},\mb{w})_{B_\epsilon}\,,\\
b((\mb{u}^\epsilon,q^\epsilon),r)&=0 
\end{aligned}
\end{equation}
for all $(\mb{w}, r,w) \in \mb{V} \times L^2(\mc{V}_\epsilon) \times D^{1,2}(\Omega_\epsilon)$.
\end{definition}

The definition of a weak solution comes from multiplying \eqref{eq:3D1D_1} by a test function $\mb{w} \in \mb{V}$ and integrating over $\mc{V}_\epsilon$. Integrating by parts and using the properties of $\mb{V}$, \eqref{eq:3D3D_3}, and \eqref{eq:3D3D_6}, we see
\begin{equation}\begin{split}
    \int_{\mc{V}_\epsilon}2\mu\,\E(\bu^\epsilon):\E(\bw)\,d\bx &=\int_{\Gamma_\epsilon} (\Sigma^\epsilon\mb{n}) \cdot\mb{w} \,dS_\epsilon + \int_{B_\epsilon} (\Sigma^\epsilon\mb{n})\cdot\mb{w}\,d\mb{x} \\
    &=\int_{\Gamma_\epsilon} \big((\Sigma^\epsilon \mb{n})\cdot\bm{n}\big)(\mb{w}\cdot \mb{n})\,dS_\epsilon -p_0\int_{B_\epsilon} \mb{w}\cdot \mb{n} \,d\mb{x}\\
    &=-\int_{\Gamma_\epsilon} \frac{1}{\kappa^\epsilon} (\mb{u}^\epsilon \cdot \mb{n}) (\mb{w} \cdot \mb{n})\, dS_\epsilon-\int_{\Gamma_\epsilon} q^\epsilon (\mb{w} \cdot \mb{n}) \,dS_\epsilon
    -p_0\int_{B_\epsilon} \mb{w}\cdot \mb{n} \,d\mb{x}\,,
\end{split}
\end{equation}
which is equivalent to 
\begin{equation} \label{weak1} 
\mathcal{A}((\mb{u}^\epsilon, q^\epsilon),(\mb{w},0))+b((\mb{w},0),p^\epsilon)=-(p_0\mb{n},\mb{w})_{B_\epsilon}\,.
\end{equation}
To complete the definition of a weak solution, we multiply \eqref{eq:3D3D_2} by $w\in D^{1,2}(\Omega_\epsilon)$ and integrate by parts. Using the first equation of \eqref{eq:3D3D_3} and \eqref{eq:3D3D_5}, we obtain
\begin{equation} \label{weak2}
\int_{\Omega_\epsilon} \zeta^\epsilon \nabla q^\epsilon \cdot \nabla w \,d\mb{x}-\int_{\Gamma_\epsilon} (\mb{u} \cdot \mb{n}) w \,d\mb{x}= \mathcal{A}(\mb{u}^\epsilon, q^\epsilon),(0,w))+b((0,w),p^\epsilon)=0\,.
\end{equation}
Note that the minus sign in front of the second term is due to the fact that $\mb{n}$ represents the inward normal vector with respect to $\Omega_\epsilon$.
Adding \eqref{weak1} and \eqref{weak2} yields \eqref{3Dweakdef}.

We will find a weak solution using a version of the Babu\v{s}ka-Brezzi theorem, which requires us to show an inf-sup condition. We aim to show that there exists a constant $\beta>0$ such that 
\begin{equation}\sup_{(\bm{\tau},0) \in \mb{V}\times \{0\}} \frac{b((\bm{\tau},0),p)}{\|\bm{\tau}\|_{\mb{V}} }\geq \beta \|p\|_{L^2(\mc{V}_\epsilon)}\,.
\end{equation}
To do so, we first show the following result for $\mb{V}_\sigma=\{\mb{w} \in \mb{V}~:~\nabla \cdot \mb{w}=0\}$.

\begin{lemma}\label{1.2}
Let $\mb{V}_\sigma^\perp$ be the orthogonal complement of $\mb{V}_\sigma$ with respect to the norm on $\mb{V}$. 
Then the map $\mb{v} \to \nabla \cdot \mb{v}$ is an isomorphism of $\mb{V}_\sigma^\perp$ onto $L^2(\mc{V}_\epsilon)$.
\end{lemma}
\begin{proof}
Let $f \in L^2(\mc{V}_\epsilon)$. Then $f=f_0+c$ where $f_0 \in L_0^2(\mc{V}_\epsilon):=\{f\in L^2(\mc{V}_\epsilon)\,:\,\int_{\mc{V}_\epsilon}f\,d\bx=0\}$ and $c \in \mathbb{R}$. Using the curvilinear cylindrical coordinates of section \ref{subsec:geometry}, let 
\begin{equation}
\mb{v}_1=c\frac{r\theta-r^2(\kappa_1(s) \sin \theta - \kappa_2(s) \cos \theta)}{(1-r\hat{\kappa})}\mb{e}_\theta\,.
\end{equation} 
Then, using \eqref{eq:grad_div}, $\nabla \cdot \mb{v}_1=c$. Additionally, let $\mb{v}_2=-\mb{v}_1 + \frac{|\mc{V}_\epsilon|}{\int_{\mc{V}_\epsilon} (\epsilon^2 a(s)^2 -r^2) d\bx}(\epsilon^2 a(s)^2-r^2)\mb{v}_1$. Then
\begin{equation}
\int_{\mc{V}_\epsilon} \nabla \cdot \mb{v}_2d\bx=\int_{\mc{V}_\epsilon}-c d\bx +\frac{c|\mc{V}_\epsilon|}{\int_{\mc{V}_\epsilon} (\epsilon^2a(s)^2-r^2) d\bx}\int_{\mc{V}_\epsilon} (\epsilon^2a(s)^2-r^2) d\bx=0\,.
\end{equation}
Therefore, $f_0-\nabla \cdot \mb{v}_2\in L_0^2(\mc{V}_\epsilon)$. Since the divergence operator is surjective from $H_0^1(\mc{V}_\epsilon)$ to $L_0^2(\mc{V}_\epsilon)$ \cite{galdi2011introduction}, there exists $\mb{v}_3\in H_0^1(\mc{V}_\epsilon)$ such that $\nabla \cdot \mb{v}_3=f_0-\nabla \cdot \mb{v}_2$.

Putting this together, we let $\mb{v}=\mb{v}_1+\mb{v}_2 + \mb{v}_3$ and see $\nabla \cdot \mb{v}=f_0+c$. Additionally,
on $\Gamma_\epsilon$, $\mb{v}=\mb{v}_1-\mb{v}_1+\mb{v}_3$.
Since $\mb{v}_3 \in H_0^1(\mc{V}_\epsilon)$, $\mb{v}$ is zero on $\Gamma_\epsilon$, so $\mb{v}-(\mb{v} \cdot \mb{n})\mb{n}=0$ on $\Gamma_\epsilon$. Therefore, $\mb{v} \in \mb{V}$.
Since $\mb{V}_\sigma$ is closed, $\mb{V}=\mb{V}_\sigma \oplus \mb{V}_\sigma^\perp$. 
Thus, $\mb{v}=\mb{v}^\perp+\mb{v}^0$ 
where $\mb{v}^\perp\in \mb{V}_\sigma^\perp$ and 
$\mb{v}^0\in\mb{V}_\sigma$. Clearly, $\nabla \cdot \mb{v}^\perp=\nabla \cdot \mb{v}=f$.
Thus, $\nabla \cdot$ is surjective from $\mb{V}_\sigma^\perp \to L^2(\mc{V}_\epsilon)$. It is also clearly injective and continuous. Thus, by the Bounded Inverse Theorem, there exists a continuous inverse from $L^2(\mc{V}_\epsilon)$ to $\mb{V}_\sigma^\perp$. More specifically, there exists a constant $C_\epsilon$ such that 
\begin{equation}\label{409}
C_\epsilon\|p\|_{L^2(\mc{V}_\epsilon)}\geq \|\mb{v}\|_{\mb{V}} \quad \text{when } \nabla \cdot \mb{v}=p\,.
\end{equation}
\end{proof}
\begin{corollary}\label{10.3}
There exists a constant $\beta_\epsilon>0$ such that 
\begin{equation} 
\sup_{\mb{v} \in \mb{V}} \int_{\mc{V}_\epsilon} \frac{p\nabla \cdot \mb{v}}{\|\mb{v}\|_{\mb{V}}}d\mb{x}\geq \beta_\epsilon \|p\|_{L^2(\mc{V}_\epsilon)} \quad \forall p \in L^2(\mc{V}_\epsilon)\,.
\end{equation}
\end{corollary}
\begin{proof}
Let $p\in L^2(\mc{V}_\epsilon)$. We know from Lemma \ref{1.2} that there exists a $\mb{w} \in \mb{V}$ such that $\nabla \cdot \mb{w} = p$ and \eqref{409} is satisfied. Thus,
\begin{equation}
\sup_{\mb{v} \in \mb{V}} \int_{\mc{V}_\epsilon} \frac{p \nabla \cdot \mb{v}}{\|\mb{v}\|_{\mb{V}}}d\mb{x} \geq \int_{\mc{V}_\epsilon}\frac{p \nabla \cdot \mb{w}}{\|\mb{w}\|_{\mb{V}}}d\mb{x}=\frac{1}{\|\mb{w}\|_{\mb{V}}}\int_{\mc{V}_\epsilon} p^2 d\mb{x}\geq \frac{1}{C_\epsilon} \|p\|_{L^2(\mc{V}_\epsilon)}\,.
\end{equation}
Therefore, with $\beta_\epsilon=1/C_\epsilon$, the result is shown.
\end{proof}

We may now show existence and uniqueness for the 3D-3D system.
Notice that 
\begin{equation} \begin{split}
|\mathcal{A}((\mb{u}, q), &(\mb{w}, w))|
\leq \|\mb{u}\|_{\mb{V}}\|\mb{w}\|_{\mb{V}} + \|q\|_{L^2(\Gamma_\epsilon)}\|\mb{w} \cdot \mb{n}\|_{L^2(\Gamma_\epsilon)} \\
&\qquad\quad +\|q\|_{D^{1,2}(\Omega_\epsilon)}\|w\|_{D^{1,2}(\Omega_\epsilon)} +\|\mb{u} \cdot \mb{n}\|_{L^2(\Gamma_\epsilon)}\|w\|_{L^2(\Gamma_\epsilon)}\\
&\leq C_\epsilon\left(\|\mb{u}\|_{\mb{V}}\|\mb{w}\|_{\mb{V}}+\|q\|_{L^2(\Gamma_\epsilon)}\|\mb{w}\|_{\mb{V}}+\|q\|_{D^{1,2}(\Omega_\epsilon)}\|w\|_{D^{1,2}(\Omega_\epsilon)}+\|\mb{u}\|_{\mb{V}}\|w\|_{L^2(\Gamma_\epsilon)}\right)\\
&\le C_\epsilon\|(\mb{u},q)\|_{\mb{V} \times D^{1,2}(\Omega_\epsilon}\|(\mb{w}, w)\|_{\mb{V} \times D^{1,2}(\Omega_\epsilon)}\,.
\end{split} 
\end{equation}
Hence, $\mathcal{A}$ is continuous. Additionally,
\begin{equation} \begin{split}
|\mathcal{A}((\mb{u}, q), (\mb{u}, q))|&=(\mb{u},\mb{u})_{\mb{V}}+(q, \mb{u} \cdot \mb{n})_{\Gamma_\epsilon} + +(\zeta^\epsilon \nabla q, \nabla q)_{\Omega_\epsilon}-(\mb{u} \cdot \mb{n}, q)_{\Gamma_\epsilon}\\
&=(\mb{u}, \mb{u})_{\mb{V}}+(\zeta^\epsilon \nabla q, \nabla q)_{\Omega_\epsilon} =\|\mb{u}\|^2_{\mb{V}}+\zeta^\epsilon \|q\|^2_{D^{1,2}(\Omega_\epsilon)}\\
&\geq C_\epsilon\|(\mb{u}, q)\|_{\mb{V} \times D^{1,2}(\Omega_\epsilon)}\,.
\end{split} \end{equation}
Furthermore, we note that $(p_0\mb{n},\cdot)_{B_\epsilon}$ is linear and continuous, so there exists $F\in \mb{V}$ such that $(F, \mb{w})_{\mb{V}}=-(p_0\mb{n},\mb{w})_{B_\epsilon}$ for all $\mb{w}\in \mb{V}$.
Therefore, by the Lax-Milgram lemma, there exists a unique $(\mb{u}^\epsilon, q^\epsilon)\in \mb{V}_\sigma \times V$ that satisfies
\begin{equation} 
\mathcal{A}((\mb{u}^\epsilon,q^\epsilon), (\mb{w},w))= -(p_0\mb{n}, \mb{w})_{B_\epsilon}\quad \forall (\mb{w},w)\in \mb{V}_\sigma \times V\,.
\end{equation}

Given $(\mb{u}, q)$, since $\mathcal{A}((\mb{u}, q), \cdot)$ is continuous and linear, there exists $(\bm{\tau},\sigma) \in \mb{V}\times V$ such that 
\begin{equation} 
\mathcal{A}((\mb{u}, q), (\mb{w},w))=(\bm{\tau}, \mb{w})_{\mb{V}} +(\sigma, w)_V
\end{equation} 
for all $(\mb{w},w) \in \mb{V} \times V$. We will define a map $A: \mb{V}\times V\to \mb{V}\times V$ that maps $(\mb{u},q)$ to $(\bm{\tau},\sigma)$ in this manner. Then a weak solution will satisfy
\begin{equation} 
A(\mb{u}^\epsilon,q^\epsilon)+(\nabla \cdot)^*p^\epsilon=(F,0)\,, \quad \nabla \cdot \mb{u}^\epsilon=0\,,
\end{equation}
where $(\nabla \cdot)^*$ is the dual of the divergence operator.
Since 
\begin{equation}
\big(A(\mb{u}^\epsilon,q^\epsilon)-(F,0),(\mb{w},w)\big)=\mathcal{A}((\mb{u}^\epsilon,q^\epsilon),(\mb{w},w))-\mathcal{F}(\mb{w})=0 \qquad \forall (\mb{w},w) \in \mb{V}_\sigma \times V\,,
\end{equation}
we see $A(\mb{u}^\epsilon,q^\epsilon)-(F,0) \in (\mb{V}_\sigma\times V)^\perp$. Also, Range $(\nabla \cdot)^*=$ Kernel $(\nabla \cdot)^\perp=\mb{V}_\sigma^\perp$. Thus, $A(\mb{u}^\epsilon, q^\epsilon)-(F,0)\in$ Range $(\nabla \cdot)^* \times 0$, so there exists a $p^\epsilon \in L^2(\mc{V}_\epsilon)$ such that $((\nabla \cdot)^*p^\epsilon,0)=A(\mb{u}^\epsilon, q^\epsilon)-(F,0)$. This means
\begin{equation}
(p^\epsilon, \nabla \cdot \mb{w})=\mathcal{A}((\mb{u}^\epsilon, q^\epsilon), (\mb{w},w))-\mathcal{F}(\mb{w})\quad \forall (\mb{w},w) \in \mb{V}\times V\,.
\end{equation}
Thus, $(\mb{u}^\epsilon, p^\epsilon, q^\epsilon)$ is a weak solution.

We now must show that there is only one weak solution. Assume there are two solutions and let $(\mb{u}_d, p_d, q_d)$ be the difference between the two solutions. Then \begin{equation}\mathcal{A}((\mb{u}_d,q_d),(\mb{w},w))+b((\mb{w},w),p_d)=0.\end{equation} Hence, $\mathcal{A}((\mb{u}_d,q_d),(\mb{u}_d,q_d))=0$. This shows $\|\mb{u}_d\|_{\mb{V}}=0$, so $\mb{u}_d=0$. Additionally, $\|\grad q_d\|_{L^2(\Omega_\epsilon)}=0$, so $q_d=0$. We see
\begin{equation} \label{uniqhelp} 
0=\mathcal{A}((\mb{u}_d,q_d),(\mb{w},w))-(p_d,\nabla \cdot \mb{w})=-(p_d, \nabla \cdot \mb{w})\qquad \forall \mb{w} \in \mb{V}\,.\end{equation} 
Using Corollary \ref{10.3}, we see 
\begin{equation}
\|p_d\|_{L^2(\mc{V}_\epsilon)} \leq C_\epsilon\sup_{\mb{w} \in \mb{V}} \frac{p_d \nabla \cdot \mb{w}}{\|\mb{w}\|_{\mb{V}}} d\mb{x} =0
\end{equation}
Thus, $p_d=0$. We conclude that $(\mb{u}^\epsilon, p^\epsilon, q^\epsilon)$ is the unique solution, thus establishing the qualitative existence and uniqueness result of Theorem \ref{thm:3D3D}.

\subsection{Key inequalities}
To obtain the $\epsilon$-dependent scaling in the energy estimate \eqref{eq:3d3d_energy}, we will rely on a series of $\epsilon$-dependent bounds for the Stokes equations in the slender domain $\mc{V}_\epsilon$. We begin with the Korn inequality. 
\begin{lemma}[Korn inequality in $\mc{V}_\epsilon$]\label{lem:korn}
Let $\mc{V}_\epsilon$ and $\Gamma_\epsilon$ be as in section \ref{sec:setup}. Given $\bv\in H^1(\mc{V}_\epsilon)$ with $\bv-(\bv\cdot \bm{n})\bm{n}=0$ on $\Gamma_\epsilon$, for $\epsilon$ sufficiently small, the following bounds hold for the symmetric gradient $\E(\bv)=\frac{1}{2}(\nabla\bv+(\nabla\bv)^{\rm T})$:
    \begin{align}
    \norm{\bv}_{L^2(\mc{V}_\epsilon)}&\le C\epsilon\norm{\E(\bv)}_{L^2(\mc{V}_\epsilon)} \label{eq:Korn_lower}\\
     \norm{\nabla\bv}_{L^2(\mc{V}_\epsilon)}&\le C\norm{\E(\bv)}_{L^2(\mc{V}_\epsilon)}\,. \label{eq:Korn_gradient} 
    \end{align}
    Here the constants $C$ are independent of $\epsilon$.
\end{lemma}
The scaling \eqref{eq:Korn_lower}, \eqref{eq:Korn_gradient} is shown in \cite{fabricius2019pressure, miroshnikova2020pressure, manjate2024mathematical, maruvsic2000two, ghosh2021modified} in related slender geometries with homogeneous Dirichlet boundary conditions. Here the vanishing tangential boundary condition on $\Gamma_\epsilon$ is essential for the inequalities \eqref{eq:Korn_lower}, \eqref{eq:Korn_gradient} to hold without lower order terms on the right hand side. In addition, the endpoint requires particular care, so we provide a full proof of Lemma \ref{lem:korn} in Appendix \ref{subsec:korn_inequality}.

In addition to the Korn inequality, we will need the following bound to obtain $\epsilon$-dependent estimates for the pressure $p^\epsilon$.
\begin{lemma}[Pressure operator in $\mc{V}_\epsilon$]\label{lem:bogovskii}
Let $\mc{V}_\epsilon$ be as in section \ref{sec:setup}. Given $g\in L^2(\mc{V}_\epsilon)$ with $|\mc{V}_\epsilon|^{-1}\int_{\mc{V}_\epsilon}g\,d\bx=g_{\rm c}$, there exists a function $\bv\in H^1(\mc{V}_\epsilon)$ satisfying 
\begin{equation}
\begin{aligned}
    \div\bv &= g \qquad \text{in }\mc{V}_\epsilon \\
    \bv &=0 \qquad \text{on }\Gamma_\epsilon\,, \qquad \bv = g_{\rm c}\frac{|\mc{V}_\epsilon|}{|B_\epsilon|}\bm{n} \qquad \text{on }B_\epsilon
\end{aligned}    
\end{equation}
along with the $\epsilon$-dependent bound
\begin{equation}\label{eq:pressure_epsilon}
    \norm{\nabla \bv}_{L^2(\mc{V}_\epsilon)}\le C\epsilon^{-1}\norm{g}_{L^2(\mc{V}_\epsilon)}\,,
\end{equation}
where $C$ is independent of $\epsilon$.
\end{lemma}
The existence of such a function $\bv$ is due to the classical construction by Bogovskii \cite{bogovskii1980solutions} and generalizations \cite{galdi2011introduction, fabricius2019pressure, manjate2024mathematical}, as well as Lemma \ref{1.2}. The $\epsilon$-dependence is shown in closely related settings in \cite {maruvsic2000two, maruvsic2001effects, ghosh2021modified, castineira2019rigorous}, including curved pipes with constant radius and straight pipes with non-constant but non-vanishing radius. We provide a proof of the scaling \eqref{eq:pressure_epsilon} in Appendix \ref{subsec:pressure_estimate} to verify that the endpoint does not introduce additional issues.

Finally, we will rely on a series of $\epsilon$-dependent trace inequalities at different locations along the vessel surface. 
\begin{lemma}\label{lem:trace_ineqs}
  The following $\epsilon$-dependent trace inequalities hold along the vessel surface. In all cases, $C$ is some constant independent of $\epsilon$.

\emph{1. Vessel surface $\Gamma_\epsilon$ into vessel exterior $\Omega_\epsilon$.} 
Given $\mc{V}_\epsilon$, $\Omega_\epsilon$ as in section \ref{sec:setup}, and any $g\in D^{1,2}(\Omega_\epsilon)$, the trace of $g$ on the vessel surface $\Gamma_\epsilon$ satisfies
\begin{equation}\label{eq:trace_Gamma2ext}
  \norm{g}_{L^2(\Gamma_\epsilon)} \le C\epsilon^{1/2}\abs{\log\epsilon}^{1/2}\norm{\nabla g}_{L^2(\Omega_\epsilon)}\,.
\end{equation}

\emph{2. Flat endpoint $B_\epsilon$ into vessel interior.}
Given $v\in H^1(\mc{V}_\epsilon)$, at the bottom cross section $B_\epsilon$, we have 
  \begin{equation}\label{eq:trace_Beps2interior}
    \norm{v}_{L^2(B_\epsilon)} \le C\big(\epsilon\norm{\nabla v}_{L^2(\mc{V}_\epsilon)} + \norm{v}_{L^2(\mc{V}_\epsilon)}\big)\,.
  \end{equation}

\emph{3. Free end region $\Gamma_\epsilon\cap\{1-\epsilon\le s\le 1\}$ into vessel exterior.} 
For $g\in D^{1,2}(\Omega_\epsilon)$, the trace of $g$ on the surface of the vessel free end region $\Gamma_\epsilon\cap\{1-\epsilon\le s\le 1\}$ satisfies 
\begin{equation}\label{eq:trace_tip2exterior}
    \norm{g}_{L^2(\Gamma_\epsilon)} \le C\epsilon^{5/4}\abs{\log\epsilon}^{1/2}\norm{\nabla g}_{L^2(\Omega_\epsilon)}\,.
  \end{equation}

\end{lemma}
The proof of Lemma \ref{lem:trace_ineqs} appears in Appendix \ref{subsec:trace_app}.

\subsection{Energy estimate}\label{subsec:energy}
Given the unique solution $(\bu^\epsilon,p^\epsilon,q^\epsilon)$ to the 3D-3D system \eqref{eq:3D3D_1}-\eqref{eq:3D3D_6} and Lemmas \ref{lem:korn}, \ref{lem:bogovskii}, and \ref{lem:trace_ineqs}, we may proceed to show the $\epsilon$-dependent energy estimate \eqref{eq:3d3d_energy}. Along the way, we obtain the $\epsilon$-scaling of the parameters $\mu^\epsilon$, $\kappa^\epsilon$, and $\zeta^\epsilon$ used in \eqref{eq:3D3D_1}-\eqref{eq:3D3D_6}.

Recall the weak form (Definition \ref{def:weak3D}) of the 3D-3D system, which we write out here for convenience:
\begin{equation}\label{eq:weakform}
\begin{aligned}
    &\int_{\mc{V}_\epsilon}\bigg(2\mu^\epsilon\,\E(\bu^\epsilon):\E(\bw)- p^{\epsilon}\div\bw\bigg)\,d\bx + \frac{1}{\kappa^\epsilon}\int_{\Gamma_\epsilon}(\bu^\epsilon\cdot\bm{n})(\bw\cdot\bm{n})\,dS_\epsilon  + \int_{\Gamma_\epsilon}q^{\epsilon}\,(\bw\cdot\bm{n})\,dS_\epsilon\\
    &\qquad +\zeta^\epsilon\int_{\Omega_\epsilon}\nabla q^\epsilon\cdot\nabla g\,d\bx - \int_{\Gamma_\epsilon}(\bu^\epsilon\cdot\bm{n})\,g\,dS_\epsilon
     = -p_0\int_{B_\epsilon}\bw\cdot\bm{n}\,d\bx\,.
\end{aligned}
\end{equation}
Taking $(\bw,g)=(\bu^\epsilon, q^\epsilon)$ in \eqref{eq:weakform}, we obtain the energy identity 
\begin{equation}\label{eq:energy_id}
    \begin{aligned}
        & 2\mu^\epsilon\int_{\mc{V}_\epsilon} \abs{\E(\bu^\epsilon)}^2\,d\bx + \frac{1}{\kappa^\epsilon}\int_{\Gamma_\epsilon}(\bu^\epsilon\cdot\bm{n})^2\,dS_\epsilon +  \zeta^\epsilon\int_{\Omega_\epsilon}\abs{\nabla q^\epsilon}^2 \,d\bx   
        =-p_0\int_{B_\epsilon}\bu^\epsilon\cdot\bm{n}\,dS_\epsilon\,.
    \end{aligned}
\end{equation}
By Lemma \ref{lem:trace_ineqs}, we have
\begin{equation}
    \begin{aligned}
    \abs{p_0\int_{B_\epsilon}\bu^\epsilon\cdot\bm{n}\,dS_\epsilon} &\le \abs{p_0}\abs{B_\epsilon}^{1/2}\norm{\bu^\epsilon\cdot\bm{n}}_{L^2(B_\epsilon)}\\
       &\le C\abs{p_0}\abs{B_\epsilon}^{1/2}\big(\epsilon\norm{\nabla\bu^\epsilon}_{L^2(\mc{V}_\epsilon)}+ \norm{\bu^\epsilon}_{L^2(\mc{V}_\epsilon)}\big) \,.
    \end{aligned}
\end{equation}
Using Lemma \ref{lem:korn} and Young's inequality, we thus obtain 
\begin{equation}
\begin{aligned}
    &\mu^\epsilon\epsilon^{-2}\norm{\bu^\epsilon}_{L^2(\mc{V}_\epsilon)}^2+ \mu^\epsilon\norm{\nabla\bu^\epsilon}_{L^2(\mc{V}_\epsilon)}^2 +\frac{1}{\kappa^\epsilon}\norm{\bu^\epsilon\cdot\bm{n}}_{L^2(\Gamma_\epsilon)}^2  + \zeta^\epsilon\norm{\nabla q^\epsilon}_{L^2(\Omega_\epsilon)}^2\\
    &\qquad\le 
    \frac{C\epsilon^2}{\mu^\epsilon}\abs{B_\epsilon}\abs{p_0}^2\,.
\end{aligned}
\end{equation}

Next, letting $p^\epsilon_{\rm c} =|\mc{V}_\epsilon|^{-1}\int_{\mc{V}_\epsilon}p^\epsilon\,d\bx$, by Lemma \ref{lem:bogovskii}, we may consider $\bv^\epsilon\in H^1(\mc{V}_\epsilon)$ satisfying
\begin{equation}\label{eq:divv_peps}
\begin{aligned}
    \div\bv^\epsilon &= p^\epsilon \quad \text{in }\mc{V}_\epsilon\,, \qquad  \bv^\epsilon\big|_{\Gamma_\epsilon}=0\,, \qquad \bv^\epsilon\big|_{B_\epsilon} = p^\epsilon_{\rm c}\frac{|\mc{V}_\epsilon|}{|B_\epsilon|}\bm{n}\,,\\
    \norm{\nabla\bv^\epsilon}_{L^2(\mc{V}_\epsilon)}&\le C\epsilon^{-1}\norm{p^\epsilon}_{L^2(\mc{V}_\epsilon)}\,.
\end{aligned}    
\end{equation}
Using the test function pair $(\bv^\epsilon,0)$ in \eqref{eq:weakform}, we then have 
\begin{equation}
    \begin{aligned}
        \int_{\mc{V}_\epsilon}\abs{p^\epsilon}^2\,d\bx &= \abs{\int_{\mc{V}_\epsilon} 2\mu^\epsilon\,\E(\bu^\epsilon):\E(\bv^\epsilon)\,d\bx
        + p_0\,p_{\rm c}^\epsilon\,\abs{\mc{V}_\epsilon} }\\
        &\le 2\mu^\epsilon\norm{\E(\bu^\epsilon)}_{L^2(\mc{V}_\epsilon)}\norm{\nabla\bv^\epsilon}_{L^2(\mc{V}_\epsilon)} + \abs{p_0}\abs{\mc{V}_\epsilon}^{1/2}\norm{p^\epsilon}_{L^2(\mc{V}_\epsilon)}\\
        &\le C\big(\mu^\epsilon \epsilon^{-1}\norm{\E(\bu^\epsilon)}_{L^2(\mc{V}_\epsilon)}+ \abs{p_0}\abs{B_\epsilon}^{1/2}\big)\norm{p^\epsilon}_{L^2(\mc{V}_\epsilon)}\,,
    \end{aligned}
\end{equation}
where we note that $|\mc{V}_\epsilon|=C|B_\epsilon|$.
Using Young's inequality and Lemma \ref{lem:korn}, we thus obtain the pressure estimate
\begin{equation}
    \epsilon^2\norm{p^\epsilon}_{L^2(\mc{V}_\epsilon)}^2\le C(\mu^\epsilon)^2\norm{\nabla \bu^\epsilon}_{L^2(\mc{V}_\epsilon)}^2 + C\epsilon^2\abs{B_\epsilon}\abs{p_0}^2\,.
\end{equation}
In total, we have
\begin{equation}
\begin{aligned}
    &(\mu^\epsilon)^2\epsilon^{-2}\norm{\bu^\epsilon}_{L^2(\mc{V}_\epsilon)}^2 + (\mu^\epsilon)^2\norm{\nabla\bu^\epsilon}_{L^2(\mc{V}_\epsilon)}^2 +\epsilon^2\norm{p^\epsilon}_{L^2(\mc{V}_\epsilon)}^2\\
    &\qquad +\frac{\mu^\epsilon}{\kappa^\epsilon}\norm{\bu^\epsilon\cdot\bm{n}}_{L^2(\Gamma_\epsilon)}^2  + \mu^\epsilon\zeta^\epsilon\norm{\nabla q^\epsilon}_{L^2(\Omega_\epsilon)}^2
    \le 
    C\epsilon^2\abs{B_\epsilon}\abs{p_0}^2\,.
\end{aligned} 
\end{equation}
Dividing through by the volume $|\mc{V}_\epsilon|$ of the vessel and using that $|\mc{V}_\epsilon|= C|B_\epsilon|\le C\epsilon^2$, we finally obtain
\begin{equation}\label{eq:energy_epsilon}
\begin{aligned}
    &\frac{1}{|\mc{V}_\epsilon|}\bigg((\mu^\epsilon)^2\epsilon^{-4}\norm{\bu^\epsilon}_{L^2(\mc{V}_\epsilon)}^2 + (\mu^\epsilon)^2\epsilon^{-2}\norm{\nabla\bu^\epsilon}_{L^2(\mc{V}_\epsilon)}^2 +\norm{p^\epsilon}_{L^2(\mc{V}_\epsilon)}^2\bigg)\\
    &\qquad+\frac{\mu^\epsilon}{\epsilon^2\kappa^\epsilon |\mc{V}_\epsilon|}\norm{\bu^\epsilon\cdot\bm{n}}_{L^2(\Gamma_\epsilon)}^2  + \frac{\mu^\epsilon\zeta^\epsilon}{\epsilon^2|\mc{V}_\epsilon|}\norm{\nabla q^\epsilon}_{L^2(\Omega_\epsilon)}^2
    \le 
    C\abs{p_0}^2\,.
\end{aligned}  
\end{equation}
Using the boundary conditions \eqref{eq:3D3D_3} on $\Gamma_\epsilon$, from the bound on $\norm{\bu^\epsilon\cdot\bm{n}}_{L^2(\Gamma_\epsilon)}$ we also obtain
\begin{equation}\label{eq:BC_epsilon}
    \frac{\mu^\epsilon(\zeta^\epsilon)^2}{\epsilon^2\kappa^\epsilon |\mc{V}_\epsilon|}\norm{\frac{\p q^\epsilon}{\p\bm{n}}}_{L^2(\Gamma_\epsilon)}^2 + \frac{\mu^\epsilon\kappa^\epsilon}{\epsilon^2 |\mc{V}_\epsilon|}\norm{(\Sigma^\epsilon\bm{n})\cdot\bm{n}+q^\epsilon}_{L^2(\Gamma_\epsilon)}^2
    \le 
    C\abs{p_0}^2\,.
\end{equation}
From \eqref{eq:energy_epsilon}, we see that $\mu^\epsilon\zeta^\epsilon\sim \epsilon^4$ is needed for $\nabla q^\epsilon$ to remain $O(1)$ as $\epsilon\to 0$, and $\mu^\epsilon\sim 1$ is needed to obtain the classical Poiseuille law scaling \cite {fabricius2019pressure, maruvsic2000two,canic2003effective} for $(\bu^\epsilon,p^\epsilon)$. Furthermore, from \eqref{eq:BC_epsilon} we see that $\kappa^\epsilon\sim\epsilon^3$ is needed for $|\Gamma_\epsilon|^{-1}\norm{(\Sigma^\epsilon\bm{n})\cdot\bm{n}+q^\epsilon}_{L^2(\Gamma_\epsilon)}^2$ to be $O(1)$ as $\epsilon\to 0$.




%% file: asymptotics.tex

\section{Formal asymptotics}\label{sec:asymptotics}
Here we present a formal derivation of the 3D-1D system \eqref{eq:3D1D_1}-\eqref{eq:3D1D_4} from the 3D-3D system \eqref{eq:3D3D_1}-\eqref{eq:3D3D_6} in the case of a straight, semi-infinite vessel with general non-constant, non-degenerate radius function:
\begin{equation}
  \X(z) = z\be_z\,, \qquad 0\le z<\infty  \,,
\end{equation}
with $a(z)\in C^2(0,\infty)$ satisfying $0<a_0\le a(z)\le 1$ everywhere. In this way, we avoid asymptotics for the free end, instead relying on geometric smallness of the endpoint in the error analysis. 

The straight centerline allows us to work in usual cylindrical coordinates $(r,\theta,z)$ centered about the $z$-axis, with corresponding unit vectors $\be_r(\theta)$, $\be_\theta(\theta)$, $\be_z$. For convenience, we rescale $R=r/\epsilon$ so that the interior of the vessel is given by $0\leq R\leq a(z)$.

We first consider the Stokes equations within the vessel. Due to axisymmetry, we consider 
\begin{equation}
\bu^\epsilon= u^{\epsilon,R}(R,z)\be_r(\theta) + u^{\epsilon,z}(R,z)\be_z\,, \qquad p^\epsilon = p^\epsilon(R,z)\,.
\end{equation}
Note that even for a curved centerline, this axisymmetry within the vessel is preserved to leading order in $\epsilon$.
Due to the boundary condition $\bu^\epsilon-(\bu^\epsilon\cdot\bm{n})\bm{n}\big|_{\Gamma_\epsilon}=0$, we note that $\bu^\epsilon\cdot\bm{n}^\perp=0$ on $\Gamma_\epsilon$, where $\bm{n}^\perp=\frac{1}{\sqrt{1+(\epsilon a')^2}}\big(\epsilon a'\be_r+\be_z \big)$. In particular, 
\begin{equation}\label{eq:nontang_BC}
  u^{\epsilon,z}=-\epsilon a'u^{\epsilon,R} \qquad \text{at }R=a\,.
\end{equation}

The energy inequality \eqref{eq:3d3d_energy} motivates the following expansions: 
\begin{equation}\label{eq:expand}
\begin{aligned}
  p^{\epsilon}(R,z) = \sum_{j=0}^\infty \epsilon^jp^{\epsilon}_j\,, \qquad
  u^{\epsilon,z}(R,z) = \epsilon^2\sum_{j=0}^\infty \epsilon^ju^{\epsilon,z}_j\,, \qquad 
  u^{\epsilon,R}(R,z) = \epsilon^3\sum_{j=0}^\infty \epsilon^ju^{\epsilon,R}_j(R,z) \,,
\end{aligned}
\end{equation}
with $u^{\epsilon,z}_0(a(z),z)=0$ due to \eqref{eq:nontang_BC}.
Writing the Stokes equations as
\begin{align}
\frac{1}{\epsilon R} \frac{\partial}{\partial R}(Ru^{\epsilon,R}) + \frac{\partial u^{\epsilon,z}}{\partial z}=0  \label{eq:fullcylindrical1}\\
-\frac{1}{\epsilon} \frac{\partial p^\epsilon}{\partial R} + \mu\left(\frac{1}{\epsilon^2R} \frac{\partial}{\partial R}\left(R\frac{\partial u^{\epsilon,R}}{\partial R}\right) + \frac{\partial^2u^{\epsilon,R}}{\partial z^2} -\frac{ u^{\epsilon,R}}{\epsilon^2R^2}\right)=0  \label{eq:fullcylindrical2}\\
-\frac{\partial p^\epsilon}{\partial z} + \mu \left(\frac{1}{\epsilon^2R}\frac{\partial}{\partial R}\left(R \frac{\partial u^{\epsilon,z}}{\partial R}\right) + \frac{\partial^2 u^{\epsilon,z}}{\partial z^2}\right)=0\,,  \label{eq:fullcylindrical3}
\end{align}
we may plug in the expansions \eqref{eq:expand} and collect matching orders in $\epsilon$. From \eqref{eq:fullcylindrical2}, we obtain
\begin{equation}
  \frac{\partial p^\epsilon_0}{\partial R} =0\,,
\end{equation}
i.e. $p^\epsilon_0=p^\epsilon_0(z)$.
Next, from \eqref{eq:fullcylindrical3}, we obtain
\begin{equation}\label{eq:899}
\frac{\partial p^\epsilon_0}{\partial z} =\mu \frac{1}{R} \frac{\partial}{\partial R} \left(R \frac{\partial u^{\epsilon,z}_0}{\partial R}\right).
\end{equation}
Using the fact that $p_0^{\epsilon,z}$ is not a function of $R$ and that $u^{\epsilon,z}_0=0$ when $R=a(z)$, integrating \eqref{eq:899} in $R$ yields 
\begin{align}
 -\frac{a^2(z)-R^2}{4\mu}\frac{\partial p^\epsilon_0}{\partial z}&=u^{\epsilon,z}_0(R,z)\,. \label{eq:uzdef}
\end{align}
Finally, combining \eqref{eq:uzdef} with the incompressibility constraint \eqref{eq:fullcylindrical1}, we see 
\begin{align}
\frac{1}{R} \frac{\partial}{\partial R} (Ru^{\epsilon,R}_0)&=\frac{a(z)a'(z)}{2\mu}\frac{\partial p^\epsilon_0}{\partial z}+\frac{a^2(z)-R^2}{4\mu}\frac{\partial^2 p^\epsilon_0}{\partial z^2}\,.
\end{align}
Upon integrating in $R$, we then obtain
\begin{equation}\label{eq:uRdef}
  u^{\epsilon,R}_0=\frac{Ra(z)a'(z)}{4\mu}\frac{\partial p^\epsilon_0}{\partial z}+\frac{2a^2(z)R-R^3}{16\mu} \frac{\partial^2 p^\epsilon_0}{\partial z^2}\,,
\end{equation}
and, on the vessel surface,
\begin{equation}\label{eq:uRBCdef} 
u^{\epsilon,R}_0\big|_{R=a(z)}= \frac{a^2(z)a'(z)}{4\mu}\frac{\partial p^\epsilon_0}{\partial z}+\frac{a^3(z)}{16\mu}\frac{\partial^2 p^\epsilon_0}{\partial z^2}\,.
\end{equation}

Before turning to the boundary conditions for $(\bu^\epsilon,p^\epsilon)$, we first consider the asymptotic behavior of the exterior pressure $q^\epsilon$ in tandem with its boundary condition on the vessel surface. For this, we will make use of the weak formulation \eqref{eq:weakform} of the 3D-3D system.
We write, in rescaled variables,
\begin{equation}
  q^\epsilon(R,\theta,z) = \sum_{j=0}^\infty \epsilon^jq^\epsilon_j(R,\theta,z)\,.
\end{equation}
Here we do not make use of axisymmetry since this is not (approximately) preserved far away from the vessel when the vessel centerline is curved. For any $g\in D^{1,2}$ exterior to the vessel, we consider the test function pair $(0,g)$ in the weak form \eqref{eq:weakform}. Using the boundary condition \eqref{eq:3D3D_3} to replace $\bu^\epsilon\cdot\bm{n}=\epsilon^3\kappa((\Sigma^\epsilon\bm{n})\cdot\bm{n}+q^\epsilon)$, we obtain 
\begin{equation}
\begin{aligned}
  \epsilon^4\zeta\int_0^\infty\int_0^{2\pi}\int_a^\infty\nabla q^\epsilon\cdot\nabla g\,\epsilon^2RdRd\theta dz+\epsilon^3\kappa \int_0^\infty\int_0^{2\pi}\big((\Sigma^\epsilon\bm{n})\cdot\bm{n}+q^\epsilon\big)\,g\,\epsilon a\,d\theta dz
     = 0\,.
\end{aligned}
\end{equation}
To leading order in $\epsilon$ (noting that the volume factor $\epsilon^2RdR$ does not remain small exterior to the vessel), we have 
\begin{equation}
  \zeta\int_0^\infty\int_0^{2\pi}\int_a^\infty\nabla q^\epsilon_0\cdot\nabla g\,\epsilon^2RdRd\theta dz+\kappa \int_0^\infty\int_0^{2\pi}\big(q^\epsilon_0-p^\epsilon_0\big)\,g\, a(z)\,d\theta dz
     = 0\,.
\end{equation}
From this, we obtain the leading order exterior equation
\begin{equation}\label{eq:qeps0}
\begin{aligned}
  \Delta q^\epsilon_0 &= 0\,,  \hspace{3cm} R>a \\
  \frac{\p q^\epsilon_0}{\p R} &= \frac{\kappa}{\zeta}\big(q^\epsilon_0-p^\epsilon_0(z)\big)\,, \qquad R=a\,.
\end{aligned}
\end{equation}
Rescaling to $r=\epsilon R$ yields the expected form \eqref{eq:3D1D_1}-\eqref{eq:3D1D_2}.

Finally, we consider the asymptotic behavior of the boundary conditions for $(\bu^\epsilon,p^\epsilon)$ on the vessel surface, for which we again rely on the weak form \eqref{eq:weakform}.
We consider interior test functions $\bw\in \bm{V}$ which inherit the asymptotic structure of $\bu^\epsilon$, in particular, the (near-) axisymmetry of $\bu^\epsilon$ within the vessel. 
Thus we consider $\bw$ of the form 
\begin{equation}\label{eq:wform}
  \bw = w^R(R,z)\be_r(\theta) + w^z(R,z)\be_z
\end{equation}
satisfying the following boundary condition at $R=a(z)$:
\begin{equation}
  w^z(a(z),z) = -\epsilon a' w^R(a(z),z)\,.
\end{equation}
Using $(\bw,0)$ in the weak form \eqref{eq:weakform}, we have
\begin{equation}
\begin{aligned}
    &\int_0^\infty\int_0^{2\pi}\int_0^a\bigg(2\mu\,\E(\bu^\epsilon):\E(\bw)- p^{\epsilon}\div\bw\bigg)\epsilon^2R\, dRd\theta dz + \int_0^\infty\int_0^{2\pi}q^{\epsilon}\,(\bw\cdot\bm{n})\,\epsilon a\,d\theta dz\\
    &\quad + \frac{1}{\epsilon^3\kappa}\int_0^\infty\int_0^{2\pi}(\bu^\epsilon\cdot\bm{n})(\bw\cdot\bm{n})\,\epsilon a\,d\theta dz 
     = p_0\int_0^{2\pi}\int_0^a\bw\cdot\be_z\,\epsilon^2R\,dRd\theta\,.
\end{aligned}
\end{equation}
First note that if we take $w^R=0$, the leading order terms in $\epsilon$ yield
\begin{equation}
\begin{aligned}
    &\int_0^\infty\int_0^{2\pi}\int_0^a\bigg(\mu\frac{\p u^{\epsilon,z}}{\p R}\frac{\p w^z}{\p R} -p^{\epsilon}_0\,\p_z w^z\bigg)\epsilon^2R\, dRd\theta dz 
     = p_0\int_0^{2\pi}\int_0^aw^z\,\epsilon^2R\,dRd\theta\,,
\end{aligned}
\end{equation}  
or, using \eqref{eq:899}, $p^{\epsilon}_0\big|_{z=0}=p_0$.

For more general $\bw$ of the form \eqref{eq:wform}, at leading order in $\epsilon$, we obtain
\begin{equation}
\begin{aligned}
  &-\int_0^\infty\int_0^{2\pi}\int_0^a p^\epsilon_0(z)\,\frac{\p}{\p R}\big(R w^R\big)\,\epsilon dRd\theta dz + \int_0^\infty\int_0^{2\pi}q^{\epsilon}_0\,w^R(a(z),z)\,\epsilon a\,d\theta dz \\
  &\qquad + \frac{1}{\kappa}\int_0^\infty\int_0^{2\pi}u^{\epsilon,R}_0\,w^R(a(z),z)\,\epsilon a\,d\theta dz = 0\,.
\end{aligned}
\end{equation}
Integrating the first term in $R$ and replacing $u^{\epsilon,R}_0$ with \eqref{eq:uRBCdef} in the third, we obtain 
\begin{equation}
\begin{aligned}
  &-\int_0^\infty\int_0^{2\pi} p^\epsilon_0(z)\,w^R(a(z),z)\,\epsilon a \,d\theta dz + \int_0^\infty\int_0^{2\pi}q^{\epsilon}_0\,w^R(a(z),z)\,\epsilon a\,d\theta dz\\
  &\quad + \frac{1}{\kappa}\int_0^\infty\int_0^{2\pi}\bigg(\frac{a^2(z)a'(z)}{4\mu}\frac{\partial p^\epsilon_0}{\partial z}+\frac{a^3(z)}{16\mu}\frac{\partial^2 p^\epsilon_0}{\partial z^2} \bigg)w^R(a(z),z)\,\epsilon a\,d\theta dz   = 0\,.
\end{aligned}
\end{equation}
Rearranging, we thus have
\begin{equation}
\begin{aligned}
  & \int_0^\infty \bigg(\int_0^{2\pi}q^{\epsilon}_0\,d\theta-2\pi p^\epsilon_0(z)\bigg)\,w^R(a(z),z)a(z)\, dz  \\
  &\qquad + \frac{\pi}{8\mu\kappa}\int_0^\infty\frac{\p}{\p z}\bigg(a^4(z)\frac{\partial p^\epsilon_0}{\partial z}\bigg)w^R(a(z),z) \,dz = 0\,.
\end{aligned}
\end{equation}
Since $w^R(a(z),z)$ is arbitrary, we obtain the boundary condition 
\begin{equation}\label{eq:peps0}
\begin{aligned}
  & \frac{\pi}{8\mu\kappa}\frac{\p}{\p z}\bigg(a^4(z)\frac{\partial p^\epsilon_0}{\partial z}\bigg) -2\pi a(z) p^\epsilon_0(z) + \int_0^{2\pi}q^{\epsilon}_0\,a(z)d\theta = 0\,,
\end{aligned}
\end{equation}
which we recognize as the straight form of the boundary condition \eqref{eq:3D1D_3alt}.

Thus, combining \eqref{eq:qeps0} and \eqref{eq:peps0}, we formally obtain the 3D-1D system \eqref{eq:3D1D_1}-\eqref{eq:3D1D_4}, ignoring a treatment of the endpoint behavior. In particular, we emphasize that the $\theta$-integrated boundary condition \eqref{eq:3D1D_3} or \eqref{eq:3D1D_3alt} arises due to the geometric constraint on $\bu^\epsilon$ on $\Gamma_\epsilon$ in the thin limit. Due to the asymptotic $\theta$-independence of $\bu^\epsilon\cdot\bm{n}$ on $\Gamma_\epsilon$, the limiting boundary value problem has fewer degrees of freedom than the original 3D-3D system, leading to a corresponding weaker, $\theta$-averaged boundary condition for $q^\epsilon$.

%% file: existence3D1D.tex
\section{Proof of Theorem \ref{thm:3D1D}}\label{sec:3D1Dthm}
We next show well-posedness of the 3D-1D model \eqref{eq:3D1D_1}-\eqref{eq:3D1D_4}.
Recalling the definition \eqref{eq:Ha_def} of the function space $\mc{H}^a$, we introduce the following notion of weak solution:
\begin{definition}\label{weak1Ddef}
We will call $(p,q)\in \mathcal{H}^a(0,1) \times D^{1,2}(\Omega_\epsilon)$ a weak solution to the 3D-1D system \eqref{eq:3D1D_1}-\eqref{eq:3D1D_4} if 
\begin{equation}\label{weak1D}
0=\int_{\Omega_\epsilon}\zeta \nabla q \cdot \nabla g d \mb{x} +\int_{\Gamma_\epsilon}\frac{\kappa}{\epsilon}(p-q)(\gamma-g)\,dS_\epsilon+\int_0^1 \frac{\pi}{8\mu}a^4(s)\frac{dp}{ds} \frac{d\gamma}{ds}ds \quad
\end{equation}
for all $(\gamma, g) \in \mathcal{H}_0^a(0,1) \times D^{1,2}(\Omega_\epsilon)$ and
 $p(0)=p_0$.
\end{definition}
Note that $\mathcal{H}^a(0,1)\subset C[0,1)$, so the boundary condition for $p$ can be understood pointwise. Due to this property, we can define the test function space $\mc{H}_0^a(0,1)$ as 
\begin{equation} 
\mc{H}_0^a(0,1)=\{\gamma \in \mathcal{H}^a(0,1)~:~\gamma(0)=0\}\,.
\end{equation}
Because of the Dirichlet boundary condition for $p$, we seek $\overline{p}(s)=p(s)-p_0$, so that $\overline{p} \in \mathcal{H}_0^a(0,1)$.

We define $d:(D^{1,2}(\Omega_\epsilon) \times \mc{H}^a(0,1))\times (D^{1,2}(\Omega_\epsilon) \times \mc{H}^a(0,1))\to \mathbb{R}$ as
\begin{equation}
d((q,\overline{p}), (g, \gamma)):=\int_{\Omega_\epsilon} \nabla q\cdot \nabla g \,d\mb{x} + \int_{\Gamma_\epsilon} \frac{\kappa}{\zeta\epsilon} (\overline{p}-q)(\gamma-g) dS_\epsilon + \int_0^1 \frac{\pi}{8\zeta\mu}a^4(s) \frac{d\overline{p}}{ds}\frac{d\gamma}{ds}ds\,.
\end{equation}
To proceed, we will need the following weighted Poincar\'e inequality for $\mc{H}_0^a(0,1)$:
\begin{lemma}\label{lem:HaPoincare}
For all $u \in \mc{H}_{0}^a(0,1)$, 
\begin{equation} \label{HbP} 
\|u\|_{L^2(0,1)}\leq C\|a^2u_s\|_{L^2(0,1)}\,.
\end{equation}
\end{lemma}
\begin{proof}
Let $u\in \mc{H}_0^a(0,1)$. Define $w:[0,\infty)\to \R$ as
\begin{equation}
w(s)=\begin{cases} u(1-s) \quad&\text{when }0 \leq s\leq 1\\
0\quad &\text{when } s>1\,.
\end{cases}
\end{equation}
We note that
\begin{equation}
w_s(s)=\begin{cases}
-u_s(1-s) \quad &\text{when }0 \leq s\leq 1\\
0 \quad &\text{when }s>1\,.
\end{cases}
\end{equation}
Using a change of variables and \cite[Lemma 1, ``Hardy's inequality at infinity"]{mironescu2018role} with $r=1$ and $q=2$, we see
\begin{equation}\begin{split}
\frac{4}{(1-C\epsilon^2)^4}&\int_0^1(a^2(t) u_t(t))^2 dt \geq 
4\int_0^1 (1-t^2)^2 (u_t(t))^2dt\geq 4\int_0^1 (1-t)^2 (u_t(t))^2 dt\\
&=4\int_0^1 s^2(u_s(1-s))^2 ds=4\int_0^\infty s^2(w_s(s))^2 ds \geq \int_0^\infty w^2(s) ds=\int_0^1 u^2(t) dt\,.
\end{split} \end{equation}
Taking the square root of both sides shows the result.
\end{proof}
Due to Lemma \ref{lem:HaPoincare}, for $v\in \mc{H}^a_0$, we may use the norm
\begin{equation}
\|v\|_{\mc{H}^a_0(0,1)}:=\|a^2v_s\|_{L^2(0,1)}\,.
\end{equation}
We see that the bilinear form $d$ satisfies
\begin{equation}\label{eq:dcoercive}
d((q,\overline{p}),(q, \overline{p})) \geq\|\nabla q\|_{L^2(\Omega_\epsilon)}^2 + \frac{\pi}{8\zeta\mu}\left\|a^2\overline{p}_s\right\|_{L^2(0,1)}^2\geq \min\left\{1,\frac{\pi}{8\zeta\mu}\right\}\|(q,\overline{p})\|^2_{D^{1,2}(\Omega_\epsilon)\times \mc{H}^a(0,1)}\,.
\end{equation}
Additionally,
\begin{equation} 
\begin{split}
|d((q, \overline{p}), (g, \gamma))| &\leq \| \nabla q \|_{L^2(\Omega_\epsilon)} \|\nabla g\|_{L^2(\Omega_\epsilon)}+\frac{\kappa}{\zeta \epsilon}\|q-\overline{p}\|_{L^2(\Gamma_\epsilon)}\|\gamma-g\|_{L^2(\Gamma_\epsilon)} \\
&\quad+ \frac{\pi}{8\zeta\mu}\|\overline{p}\|_{\mc{H}_0^a(0,1)}\|\gamma\|_{\mc{H}_0^a(0,1)}\\
&\leq C_\epsilon \left(\|\nabla q\|_{L^2(\Omega_\epsilon)}\|\nabla g\|_{L^2(\Omega_\epsilon)}+ \|q\|_{L^2(\Gamma_\epsilon)}\|\gamma\|_{L^2(\Gamma_\epsilon)} + \|q\|_{L^2(\Gamma_\epsilon)}\|g\|_{L^2(\Gamma_\epsilon)}\right.\\
&\left.\quad +  \|\gamma\|_{L^2(\Gamma_\epsilon)}\|\overline{p}\|_{L^2(\Gamma_\epsilon)} +\|\overline{p}\|_{L^2(\Gamma_\epsilon)}\|g\|_{L^2(\Gamma_\epsilon)}+\|\overline{p}\|_{\mc{H}_0^a(0,1)}\|\gamma\|_{\mc{H}_0^a(0,1)}\right)\\
&\le C_\epsilon\|(q,\overline{p})\|_{D^{1,2}(\Omega_\epsilon)\times\mc{H}^a_0(0,1)}\|(g,\gamma)\|_{D^{1,2}(\Omega_\epsilon)\times\mc{H}^a_0(0,1)}\,.
\end{split} \end{equation}
Finally, the term arising due to lifting the boundary condition $p(0)=p_0$ satisfies
\begin{equation}\begin{split} \label{contright}
\int_{\Gamma_\epsilon}\frac{\kappa}{\zeta\epsilon}p_0 (g-\gamma) \,dS_\epsilon
&\leq C\epsilon^{-1/2} \abs{p_0}\|g-\gamma\|_{L^2(\Gamma_\epsilon)}\,.
\end{split} \end{equation}

Combining \eqref{contright} with the fact that $d$ is bounded and coercive, by the Lax-Milgram lemma, there is a unique solution $(q,\overline{p})\in D^{1,2}(\Omega_\epsilon) \times \mc{H}_0^a(0,1)$ satisfying
\begin{equation}\label{eq:eqn}
\begin{split}
\int_{\Omega_\epsilon}  \nabla q \cdot\nabla g dx &+ \int_{\Gamma_\epsilon}\frac{\kappa}{\zeta\epsilon} (\overline{p}-q)(\gamma-g) dS_\epsilon + \int_0^1 \frac{\pi}{8\zeta\mu}a^4(s) \frac{d\overline{p}}{ds}\frac{d\gamma}{ds}ds\\
&=\int_{\Gamma_\epsilon} \frac{\kappa}{\zeta\epsilon} p_0(g-\gamma)\, dS_\epsilon\,.
\end{split}\end{equation}
Furthermore, using $(g,\gamma)=(q,\overline{p})$ in \eqref{eq:eqn} and using Young's inequality to split \eqref{contright}, we have
\begin{equation} \label{bound}
\norm{q}_{D^{1,2}(\Omega_\epsilon)}^2 + \frac{1}{\epsilon}\norm{\overline{p}-q}_{L^2(\Gamma_\epsilon)}^2 + \norm{\overline{p}}_{\mc{H}^a(0,1)}^2 \leq C\abs{p_0}^2
\end{equation}
for $C$ independent of $\epsilon$.
Thus there is a unique solution satisfying Definition \ref{weak1Ddef} and the energy estimate \eqref{eq:energy3D1D}.

%% file: errorest.tex

\section{Error estimate}\label{sec:error_estimate}
Given the solution theory of sections \ref{sec:3D3Dthm} and \ref{sec:3D1Dthm} for the 3D-3D system and 3D-1D system, respectively, along with the asymptotics of section \ref{sec:asymptotics}, our goal is now to show that the 3D-1D solution $(p(s),q)$ approximates the 3D-3D solution $(p^\epsilon,q^\epsilon)$ in the sense of Theorem \ref{thm:error}. As mentioned, we will not show Theorem \ref{thm:error} directly, as we will need details about the lower dimensional approximant that are not straightforward to obtain. Instead, we will rely on the 1D model \eqref{eq:qSB}-\eqref{eq:pSB} and the convergence result of Theorem \ref{thm:3D1Dto1D} to obtain Theorem \ref{thm:error} as a consequence of the 3D-3D to 1D convergence result of Theorem \ref{thm:3D3Dto1D}.

We begin by stating the properties of the 1D solution $(p^{\rm SB},q^{\rm SB})$ that will be critical for the error analysis. These are taken directly from Part I \cite{partI}, where the 1D model is explored in greater detail. We then proceed to the proof of Theorem \ref{thm:3D3Dto1D} in sections \ref{subsec:stokes} and \ref{subsec:error}.

\subsection{Properties of the 1D model}
We will rely crucially on some of the additional information we can extract from the 1D model \eqref{eq:qSB}-\eqref{eq:pSB} that is not as easy to see for the 3D-1D model. Here we outline the analytical properties of the 1D solution that will be especially useful. 
Note that it will occasionally be useful to write the $p^{\rm SB}$ equation \eqref{eq:pSB} as
\begin{equation}\label{eq:pSB_eqn}
	\frac{d}{ds}\bigg(a^4\frac{dp^{\rm SB}}{ds}\bigg) = \frac{8\mu\kappa}{\pi}a(s)\int_0^{2\pi} \big( p^{\rm SB}(s) - q^{\rm SB}\big|_{\Gamma_\epsilon}\big)\, d\theta \,.
\end{equation}

We begin with bounds for the interior pressure $p^{\rm SB}$. Here it is advantageous to use the 1D model rather than the 3D-1D model due to a slight improvement in behavior at the tip $s=1$. In particular, we may obtain $L^\infty$ bounds for $p^{\rm SB}$ and its derivatives with improved weights at the vessel endpoint: we gain one factor of $a(s)$ for $p^{\rm SB}$ and two factors of $a$ for higher derivatives over na\"ive \emph{a priori} estimates for the 3D-1D model. Again, this improvement is due to replacing the surface element $\mc{J}_\epsilon\sim \epsilon a(s) + C\epsilon^2$ in the 3D-1D boundary condition \eqref{eq:3D1D_3} with $\epsilon a(s)$ in the 1D approximation. The improvement in weights comes at the expense of a factor of $\epsilon^{-1/2}$. However, we can take advantage of this tradeoff very close to the tip (i.e. $o(\epsilon)$ away from $s=1$). 
The statement from Part I, Theorem 1.3 \cite{partI} is as follows:
\begin{lemma}[Estimates for the $p^{\rm SB}$ equation \cite{partI}]\label{lem:press_derivs}
Given a vessel $\mc{V}_\epsilon$ as in section \ref{sec:setup} and pressure data $p_0\in\R$, for $\epsilon$ sufficiently small, the unique solution $p^{\rm SB}(s)\in \mc{H}^a(0,1)$ to the integrodifferential equation \eqref{eq:pSB} additionally satisfies the weighted $L^2$ bounds 
\begin{equation}
\begin{aligned}
	\norm{p^{\rm SB}}_{L^2(0,1)}\le C\abs{p_0}\,, &\qquad
	\norm{a^2p^{\rm SB}_s}_{L^2(0,1)}\le C\abs{p_0}\,,\\
	\norm{a^{-1/2}\big(a^4p^{\rm SB}_s\big)_s}_{L^2(0,1)} \le C\abs{p_0}\,,
	&\qquad 
	\norm{a^{3/2}(a^4p^{\rm SB}_s)_{ss}}_{L^2(0,1)} \le C\abs{\log\epsilon}\abs{p_0}\,.
\end{aligned}
\end{equation}
Furthermore, $p^{\rm SB}$ satisfies the following weighted $L^\infty$ bounds along the extent of the vessel: 
\begin{equation}
\begin{aligned}
	\norm{p^{\rm SB}}_{L^\infty(0,1)}\le C\epsilon^{-1/2}\abs{p_0}\,,
	&\qquad
	\norm{a\,p^{\rm SB}_s}_{L^\infty(0,1)}\le C\epsilon^{-1/2}\abs{p_0}\,,\\
	\norm{a^3p^{\rm SB}_{ss}}_{L^\infty(0,1)}\le C\epsilon^{-1/2}\abs{p_0}\,,
	&\qquad
	\norm{a^5p^{\rm SB}_{sss}}_{L^\infty(0,1)}\le C\epsilon^{-1/2}\abs{\log\epsilon}\abs{p_0}\,.
\end{aligned}
\end{equation}
In all cases, $C$ is bounded independent of $\epsilon$ as $\epsilon\to 0$.
\end{lemma}

The next important feature of the 1D model is the near-$\theta$-independence of the exterior pressure approximation $q^{\rm SB}$ along $\Gamma_\epsilon$. As mentioned, this type of $\theta$-independence condition is typical of slender body theories and plays a crucial role in the error analysis. However, it is not straightforward to show this directly for the 3D-1D exterior pressure $q$. Even without the endpoint, some kind of 1D asymptotic approximation of the 3D-1D model seems necessary to obtain this near-$\theta$-independence condition in general. For the 1D model, it is straightforward to show the following.
\begin{lemma}[Near $\theta$-independence of exterior pressure]\label{lem:qtheta}
Given a vessel $\mc{V}_\epsilon$ as in section \ref{sec:setup} and pressure data $p_0\in \R$, let $q^{\rm SB}(\bx)$ be given by \eqref{eq:qSB} for $p^{\rm SB}(s)$ satisfying \eqref{eq:pSB}. For $\epsilon$ sufficiently small, along $\Gamma_\epsilon$, $q^{\rm SB}$ satisfies
	\begin{equation}
		\norm{q^{\rm SB}-\frac{1}{2\pi}\int_0^{2\pi}q^{\rm SB}\,d\theta}_{L^2(\Gamma_\epsilon)} \le C\epsilon \abs{\log\epsilon}\abs{p_0}\,.
	\end{equation}
\end{lemma}

\begin{proof}
From the proof of Lemma 5.1 in Part I, we have the following pointwise bound along $\Gamma_\epsilon$:
\begin{equation}
\begin{aligned}
	\abs{\frac{\p q^{\rm SB}}{\p\theta}}&\le C\epsilon \abs{\log\epsilon}a(s)\norm{(a^4p^{\rm SB}_s)_s}_{L^\infty(0,1)} \le C\epsilon^{1/2}\abs{\log\epsilon}a(s)\abs{p_0}\,.
\end{aligned}
\end{equation}
We may thus calculate
\begin{equation}
\begin{aligned}
	\norm{q^{\rm SB}-\frac{1}{2\pi}\int_0^{2\pi}q^{\rm SB}\,d\theta}_{L^2(\Gamma_\epsilon)} &\le C\norm{\p_\theta q^{\rm SB}}_{L^2(\Gamma_\epsilon)}
	\le C\abs{\Gamma_\epsilon}^{1/2}\norm{\p_\theta q^{\rm SB}}_{L^\infty(\Gamma_\epsilon)}\,.
\end{aligned}
\end{equation}
Using that $\abs{\Gamma_\epsilon}\sim\epsilon$, we obtain Lemma \ref{lem:qtheta}.
\end{proof}

Finally, we will make use of the 3D-1D to 1D convergence result of Theorem \ref{thm:3D1Dto1D} to note that the 1D approximation $(p^{\rm SB},q^{\rm SB})$ approximately satisfies the weak form of the 3D-1D system (Definition \ref{weak1Ddef}). In particular, plugging $(p^{\rm SB},q^{\rm SB})$ into the 3D-1D model, we obtain the following error bound. 
\begin{lemma}[1D solution in 3D-1D weak form]\label{lem:weakdiff}
Given $\mc{V}_\epsilon$ as in section \ref{sec:setup} and pressure data $p_0\in \R$, let $(p^{\rm SB},q^{\rm SB})$ be the solution to the 1D approximation given by \eqref{eq:pSB}, \eqref{eq:qSB}, respectively. For any $\wt g\in D^{1,2}(\Omega_\epsilon)$, the exterior pressure $q^{\rm SB}$ satisfies 
\begin{equation}
\begin{aligned}
	&\int_{\Omega_\epsilon}\nabla q^{\rm SB}\cdot\nabla \wt g\,d\bx - \frac{\kappa}{\zeta\epsilon}\int_{\Gamma_\epsilon}(p^{\rm SB}(s)-q^{\rm SB})\wt g\,dS_\epsilon = \wt E_1(\wt g)\,, \\
	&\qquad |\wt E_1(\wt g)| \le C\epsilon^{1/2}\abs{\log\epsilon}^{3/2}\abs{p_0}\norm{\nabla \wt g}_{L^2(\Omega_\epsilon)}\,.
\end{aligned}
\end{equation}
Furthermore, for any $\gamma(s)\in L^2(0,1)$, the interior pressure $p^{\rm SB}(s)$ satisfies
\begin{equation}
\begin{aligned}
	&\int_{\Gamma_\epsilon}(p^{\rm SB}(s)-q^{\rm SB})\gamma(s)\,dS_\epsilon - \frac{\pi\epsilon}{8\kappa\mu}\int_0^1\big(a^4p^{\rm SB}_s\big)_s\gamma(s)\,ds = \wt E_2(\gamma)\,, \\
	&\qquad |\wt E_2(\gamma)| \le C\epsilon^{3/2}\abs{p_0}\norm{\gamma}_{L^2(0,1)} \,.
\end{aligned}
\end{equation}
\end{lemma}

\begin{proof}
By Lemma 5.1 in Part I \cite{partI}, we have that the exterior pressure approximation $q^{\rm SB}(\bx)$ satisfies the following PDE in $\Omega_\epsilon$:
\begin{equation}\label{eq:qSB_PDE}
\begin{aligned}
	\Delta q^{\rm SB}&= 0 \hspace{5.5cm} \text{in }\Omega_\epsilon\\
	\frac{\p q^{\rm SB}}{\p\bm{n}}&= -\frac{\kappa}{\zeta\epsilon}(p^{\rm SB}(s)-q^{\rm SB}) + \mc{R}_\epsilon(s,\theta) \qquad \text{on }\Gamma_\epsilon\\
	\frac{\p q^{\rm SB}}{\p\bm{n}}&=0 \hspace{5.5cm} \text{on } \p\Omega_\epsilon\backslash\Gamma_\epsilon\,.
\end{aligned}
\end{equation}
Note that here the sign of the normal vector $\bm{n}$ on $\Gamma_\epsilon$ is reversed compared to Part I: here $\bm{n}$ is the outward normal vector to $\mc{V}_\epsilon$, pointing into $\Omega_\epsilon$. The remainder term $\mc{R}_\epsilon$ satisfies the pointwise bound
\begin{equation}
\begin{aligned}
	\abs{\mc{R}_\epsilon(s,\theta)}&\le C\min\{a^{-1}(s),\epsilon^{-1}\}\big(\epsilon^{-1/2}\abs{\log\epsilon}\abs{p_0} \big) \,.
\end{aligned}
\end{equation}
Using the form \eqref{eq:free_jacfac} of the Jacobian factor $\mc{J}_\epsilon$ and the bound $\abs{\mc{J}_\epsilon-\epsilon a(s)}\le C\epsilon^2$, we may also obtain the following $L^2$ bound for $\mc{R}_\epsilon$ along $\Gamma_\epsilon$:
\begin{equation}
\begin{aligned}
	\norm{\mc{R}_\epsilon}_{L^2(\Gamma_\epsilon)}^2 &\le C\bigg(\int_0^1\int_0^{2\pi} \min\{a^{-1}(s),\epsilon^{-1}\}^2\,\epsilon a(s)\,d\theta ds\\
	&\qquad + \int_0^1\int_0^{2\pi} \min\{a^{-1}(s),\epsilon^{-1}\}^2\,\epsilon^2\,d\theta ds\bigg)\epsilon^{-1}\abs{\log\epsilon}^2\abs{p_0}^2\\
	&\le C\bigg(\int_0^1a^{-1}(s)\,ds\bigg)\abs{\log\epsilon}^2\abs{p_0}^2
	\le C\abs{\log\epsilon}^2\abs{p_0}^2\,,
\end{aligned}
\end{equation}
where we have used that $a^{-1}(s)$ is integrable on $(0,1)$ due to the spheroidal endpoint condition \eqref{eq:spheroidal}, and the value is independent of $\epsilon$. 
From \eqref{eq:qSB_PDE}, we may then calculate that for any $\wt g\in D^{1,2}(\Omega_\epsilon)$, we have 
\begin{equation}
\begin{aligned}
	&\int_{\Omega_\epsilon}\nabla q^{\rm SB}\cdot\nabla \wt g\,d\bx - \frac{\kappa}{\zeta\epsilon}\int_{\Gamma_\epsilon}(p^{\rm SB}(s)-q^{\rm SB})\wt g\,dS_\epsilon 
	 = \int_{\Gamma_\epsilon}\mc{R}_\epsilon(\wt g-\gamma(s))\,dS_\epsilon \,,
\end{aligned}
\end{equation}
where the right hand side satisfies the estimate 
\begin{equation}
\begin{aligned}
	\abs{\int_{\Gamma_\epsilon}\mc{R}_\epsilon(\wt g-\gamma(s))\,dS_\epsilon}&\le C\abs{\log\epsilon}\abs{p_0}\norm{g-\gamma(s)}_{L^2(\Gamma_\epsilon)}\\
	&\le C\epsilon^{1/2}\abs{\log\epsilon}^{3/2}\abs{p_0}\norm{\nabla \wt g}_{L^2(\Omega_\epsilon)} + C\epsilon^{1/2}\abs{\log\epsilon}\abs{p_0}\norm{\gamma}_{L^2(0,1)}\,.
\end{aligned}
\end{equation}

Next, again using the form \eqref{eq:free_jacfac} of $\mc{J}_\epsilon$, we note that
\begin{equation}
	\abs{\mc{J}_\epsilon(s,\theta)-\epsilon a\sqrt{1+(\epsilon a')^2}} \le C\epsilon^2a^2
\end{equation}
for some $C$ depending only on the curvature of the vessel centerline. Then, using the equation \eqref{eq:pSB_eqn} satisfied by $p^{\rm SB}$, we may write
\begin{equation}
\begin{aligned}
	\int_{\Gamma_\epsilon}(p^{\rm SB}-q^{\rm SB})\gamma(s)\,dS_\epsilon &=
	\int_0^1\bigg(\int_0^{2\pi}(p^{\rm SB}-q^{\rm SB}) \,d\theta\bigg)\gamma(s) \epsilon a\sqrt{1+(\epsilon a')^2}\,ds + Q_0\\
	&= \frac{\pi}{8\mu\kappa}\int_0^1a^{-1}(s)\p_s(a^4\p_sp^{\rm SB})\gamma(s) \epsilon a\sqrt{1+(\epsilon a')^2}\,ds + Q_0\\
	&= \frac{\pi\epsilon}{8\mu\kappa}\int_0^1\p_s(a^4\p_sp^{\rm SB})\gamma(s) \,ds + Q_0 + Q_1\,.
\end{aligned}
\end{equation}
Here the remainder terms $Q_0$ and $Q_1$ satisfy
\begin{equation}
\begin{aligned}
	\abs{Q_0}&\le C\epsilon\norm{p^{\rm SB}-q^{\rm SB}}_{L^2(\Gamma_\epsilon)}\norm{\gamma}_{L^2(\Gamma_\epsilon)} \le C\epsilon^{3/2}\norm{p^{\rm SB}-q^{\rm SB}}_{L^2(\Gamma_\epsilon)}\norm{\gamma}_{L^2(0,1)}\\
	&\le C\epsilon^2\abs{p_0}\norm{\gamma}_{L^2(0,1)}\,,\\
	\abs{Q_1}&\le C\abs{\int_0^1a^{-1}(s)\p_s(a^4\p_sp^{\rm SB})\gamma(s) \epsilon^2\abs{aa'} ds} \\
	&\le C\epsilon^2\norm{a^{-1}\p_s(a^4p^{\rm SB}_s)}_{L^\infty(0,1)}\norm{\gamma}_{L^2(0,1)}
	\le C\epsilon^{3/2}\abs{p_0}\norm{\gamma}_{L^2(0,1)}\,,
\end{aligned}
\end{equation}
where for $Q_0$ we have used the 3D-1D energy estimate \eqref{eq:energy3D1D} along with the 3D-1D to 1D error bound of Theorem \ref{thm:3D1Dto1D} to bound $p^{\rm SB}-q^{\rm SB}$ on $\Gamma_\epsilon$, and for $Q_1$ we have applied Lemma \ref{lem:press_derivs}.
Altogether, we obtain Lemma \ref{lem:weakdiff}. 
\end{proof}


\subsection{Stokes residuals}\label{subsec:stokes}
Equipped with the relevant analytical properties of the 1D model, we may proceed to show the 3D-3D to 1D error bound of Theorem \ref{thm:3D3Dto1D}. We will begin by using the asymptotics of section \ref{sec:asymptotics} to construct a velocity field ansatz for the 1D approximation which, along with $p^{\rm SB}(s)$, nearly satisfies the Stokes equations in $\mc{V}_\epsilon$.

We consider $p^{\rm SB}(s)$ satisfying the 1D integrodifferential equation \eqref{eq:pSB} and, following the asymptotics of section \ref{sec:asymptotics}, for $\bx=\bx(r,\theta,s)\in \mc{V}_\epsilon$, we build the following velocity field ansatz:
\begin{equation}\label{eq:Udef}
\bm{U}(\bx) = -\frac{1}{16\mu}\p_s\big( \phi_\epsilon(s)(r^3 -2(\epsilon a)^2r) p^{\rm SB}_s \big)\be_r(s,\theta)+ \frac{1}{4\mu}\phi_\epsilon(s)\big(r^2-(\epsilon a)^2\big)p^{\rm SB}_s\,\be_{\rm t}(s)\,.
\end{equation} 
Here $(r,\theta,s)$ are the curvilinear cylindrical coordinates defined in section \ref{subsec:geometry}, and $\phi_\epsilon(s)$ is a smooth cutoff function in $s$ satisfying 
\begin{equation}\label{eq:phieps_def}
	\phi_\epsilon(s) = \begin{cases}
		1\,, & 0\le s\le 1-2\epsilon^{4/3}\\
		0 \,, & 1-\epsilon^{4/3}\le s\le 1\,,
	\end{cases}
	\qquad \abs{\p_s\phi_\epsilon}\le C\epsilon^{-4/3}\,.
\end{equation}
As mentioned in the introduction, the choice of the cutoff extent is optimized to balance two competing effects: (1) The ansatz \eqref{eq:Udef} is not a good approximation of $\bu^\epsilon$ at the tip: the gradient $\nabla\bm{U}$ is (\emph{a priori}) unbounded as $s\to 1$, and, more importantly for the error analysis, $\bm{U}$ does a worse and worse job of approximately satisfying the boundary condition $\bm{U}-(\bm{U}\cdot\bm{n})\bm{n}=0$ on $\Gamma_\epsilon$ as $s\to 1$.
(2) To counteract the effect (1), we will rely on the geometric smallness of the tip, which works to our advantage due to the $L^\infty$ bounds of Lemma \ref{lem:press_derivs} for $p^{\rm SB}$ up to the endpoint. However, to make use of this, the tip region has to be small enough. The choice of $\epsilon^{4/3}$ in \eqref{eq:phieps_def} balances effects (1) and (2) to yield the $\epsilon^{1/6}\abs{\log\epsilon}$ convergence rate of Theorem \ref{thm:3D3Dto1D}. 

Now, by construction, the pair $(\bm{U},p^{\rm SB}(s))$ nearly satisfies the Stokes equations within $\mc{V}_\epsilon$. In particular, we may show the following.
\begin{proposition}[1D velocity ansatz]\label{prop:almostStokes}
Given incoming pressure data $p_0\in \R$, let $p^{\rm SB}(s)$ be given by the solution to the 1D approximation \eqref{eq:pSB}.
Given any $\bw\in H^1(\mc{V}_\epsilon)$ with $\bw-(\bw\cdot\bm{n})\bm{n}=0$ on $\Gamma_\epsilon$, the velocity field ansatz $\bm{U}$ given by \eqref{eq:Udef} satisfies
	 \begin{equation}\label{eq:almostStokes}
 	\begin{aligned}
 	\int_{\mc{V}_\epsilon}\bigg(2\mu\,\E(\bm{U}):\E(\bw) &-\big(p^{\rm SB}(s)-p_0\big)\,\div\bw \bigg)\,d\bx \\
 	&= -\int_0^1\big(p^{\rm SB}(s)-p_0\big)\int_0^{2\pi}(\bw\cdot\bm{n}) \,\mc{J}_\epsilon\,d\theta\,ds + E_0\,,\\
 	|E_0|&\le C\abs{p_0}\big(\epsilon^{13/6}\norm{\nabla\bw}_{L^2(\mc{V}_\epsilon)} + \epsilon^{5/3}\norm{\bw\cdot\bm{n}}_{L^2(\Gamma_\epsilon)}\big)\,.
 \end{aligned}
 \end{equation}
\end{proposition}

\begin{proof}[Proof of Proposition \ref{prop:almostStokes}]
First, using the form \eqref{eq:grad_div} of the divergence in curvilinear coordinates, within $\mc{V}_\epsilon$ we may calculate 
\begin{equation}\label{eq:divU}
\begin{aligned}
	\div\bm{U} &= \frac{-1}{4\mu(1-r\wh\kappa)}\bigg(\frac{1}{4r}\p_r\bigg((1-r\wh\kappa)\p_s\big[\phi_\epsilon\,(r^4-2(\epsilon a)^2r^2)p^{\rm SB}_s\big] \bigg) - \p_s\big[\phi_\epsilon\,(r^2-(\epsilon a)^2)p^{\rm SB}_s \big] \bigg)\\
	&= \frac{\wh\kappa}{16\mu(1-r\wh\kappa)}\p_s\big[\phi_\epsilon\,(5r^3-6(\epsilon a)^2r)p^{\rm SB}_s\big] \,.
\end{aligned}
\end{equation}
Using the form \eqref{eq:volume_element} of the volume element $\mc{J}$, integrating in $r$, and applying Lemma \ref{lem:press_derivs}, we may then bound
\begin{equation}\label{eq:divUfbd}
\begin{aligned}
	\|\div\bm{U}\|_{L^2(\mc{V}_\epsilon)}^2 &\le 
	C\int_0^{1-\epsilon^{4/3}}\int_0^{2\pi}\int_0^{\epsilon a}\bigg((r^3+(\epsilon a)^2r)^2(p^{\rm SB}_{ss})^2+\epsilon^4r^2(p^{\rm SB}_s)^2\bigg)\,r\,drd\theta ds \\
	&\qquad + C\int_{1-2\epsilon^{4/3}}^{1-\epsilon^{4/3}}\int_0^{2\pi}\int_0^{\epsilon a}\abs{\p_s\phi_\epsilon}^2\epsilon^6a^4(a\,p^{\rm SB}_s)^2\,r\,drd\theta ds\\
	&\le C\epsilon^8a^2\big(\norm{a^3p^{\rm SB}_{ss}}_{L^\infty(0,1)}^2 +\norm{a\,p^{\rm SB}_s}_{L^\infty(0,1)}^2\big)\\
	&\qquad + C\sup_{s\in[1-2\epsilon^{4/3},1-\epsilon^{4/3}]}\epsilon^{4/3} \epsilon^{-8/3}\epsilon^8a^6\norm{a\,p^{\rm SB}_s}_{L^\infty(0,1)}^2
	\le C\epsilon^7\abs{p_0}^2\,.
\end{aligned}
\end{equation}
Here we have used that $\p_s\phi_\epsilon$ is supported only for $s\in[1-2\epsilon^{4/3},1-\epsilon^{4/3}]$, and, by the spheroidal endpoint condition \eqref{eq:spheroidal}, we may bound $\sup_{s\in[1-2\epsilon^{4/3},1-\epsilon^{4/3}]}a\le C\epsilon^{2/3}$. 
We further note that along the vessel surface $\Gamma_\epsilon$, using the form \eqref{eq:free_jacfac} of the surface Jacobian $\mc{J}_\epsilon$ and the $L^\infty$ bounds of Lemma \ref{lem:press_derivs}, we may estimate
\begin{equation}
\begin{aligned}
	\|\div\bm{U}\|_{L^2(\Gamma_\epsilon)}^2 
	&\le C\int_0^{1-\epsilon^{4/3}}\int_0^{2\pi}\bigg(\epsilon^6(a^3p^{\rm SB}_{ss})^2+\epsilon^6(a\,p^{\rm SB}_s)^2\bigg)\,\epsilon a\,d\theta ds \\
	&\qquad + C\int_{1-2\epsilon^{4/3}}^{1-\epsilon^{4/3}}\int_0^{2\pi}\abs{\p_s\phi_\epsilon}^2\epsilon^6a^4(a\,p^{\rm SB}_s)^2\,\epsilon a\,d\theta ds\\
	&\le 
	C\big(\epsilon^7 + \sup_{s\in[1-2\epsilon^{4/3},1-\epsilon^{4/3}]}\epsilon^{4/3}\epsilon^{-8/3}\epsilon^7a^5\big)\epsilon^{-1}\abs{p_0}^2 
	\le C\epsilon^6\abs{p_0}^2\,,
\end{aligned}
\end{equation}
and along the bottom cross section $B_\epsilon$, we may calculate 
\begin{equation}
\begin{aligned}
	\|\div\bm{U}\|_{L^2(B_\epsilon)}^2 &\le C\int_0^{2\pi}\int_0^{\epsilon a}\bigg(\epsilon^6(a^3p^{\rm SB}_{ss})^2+\epsilon^6(a\,p^{\rm SB}_s)^2\bigg)\bigg|_{s=0}\,r\,drd\theta   \\
	&\le C\epsilon^8\big(\norm{a^3p^{\rm SB}_{ss}}_{L^\infty(0,1)}^2+\norm{a\,p^{\rm SB}_s}_{L^\infty(0,1)}^2 \big)
	\le C\epsilon^7\abs{p_0}^2\,.
\end{aligned}
\end{equation}
Here we have also used that $\p_s\phi_\epsilon$ is supported far away from the $s=0$ cross section.

Within the vessel $\mc{V}_\epsilon$, we may thus write 
\begin{equation}\label{eq:divEU}
	2\,\div(\E(\bm{U})) = \Delta\bm{U} + \nabla\mc{U}_0 \,, \qquad \mc{U}_0 := \div\bm{U}\,,
\end{equation}
where, due to the above bounds, $\mc{U}_0$ may be treated as a remainder term. 
We may further decompose
\begin{equation}\label{eq:DeltaU}
	\Delta\bm{U} = \Delta(\bm{U}\cdot\be_{\rm t})\be_{\rm t} - \underbrace{ (\bm{U}\cdot\be_{\rm t})\Delta\be_{\rm t}}_{\mc{U}_1} + \underbrace{2\,\div\big((\bm{U}\cdot\be_{\rm t})\nabla\be_{\rm t}\big)+\Delta\big((\bm{U}\cdot\be_r)\be_r\big)}_{\div(\mc{U}_2)}\,,
\end{equation}
where the vector-valued remainder $\mc{U}_1$ satisfies 
\begin{equation}
\begin{aligned}
	\|\mc{U}_1\|_{L^2(\mc{V}_\epsilon)}^2 &\le C\int_0^{1-\epsilon^{4/3}}\int_0^{2\pi}\int_0^{\epsilon a}\epsilon^4 (a^2\,p^{\rm SB}_s)^2 \,r\,drd\theta ds 
	\le C\epsilon^6\norm{a^2p^{\rm SB}_s}_{L^2(0,1)}^2 \le \epsilon^6\abs{p_0}^2\,, 
\end{aligned}
\end{equation}
and the matrix-valued remainder term $\mc{U}_2$ satisfies  
\begin{equation}
\begin{aligned}
	\|\mc{U}_2\|_{L^2(\mc{V}_\epsilon)}^2 &\le C\|\mc{U}_1\|_{L^2(\mc{V}_\epsilon)}^2 + C\int_0^{1-\epsilon^{4/3}}\int_0^{2\pi}\int_0^{\epsilon a}\bigg(\epsilon^6a^{-4}(a^5p^{\rm SB}_s)_{ss}^2 +\epsilon^4a^{-2}(a^3p^{\rm SB}_s)_s^2 \bigg)\,r\,drd\theta ds\\
	&\qquad + \frac{C}{\epsilon^{8/3}}\int_{1-2\epsilon^{4/3}}^{1-\epsilon^{4/3}}\int_0^{2\pi}\int_0^{\epsilon a}\bigg(\epsilon^6(a^3p^{\rm SB}_s)_s^2 + \epsilon^4a^2(a\,p^{\rm SB}_s)^2\bigg)\,r\,drd\theta ds \\
	&\qquad + \frac{C}{\epsilon^{16/3}}\int_{1-2\epsilon^{4/3}}^{1-\epsilon^{4/3}}\int_0^{2\pi}\int_0^{\epsilon a}\epsilon^6a^4(a\,p^{\rm SB}_s)^2\,r\,drd\theta ds\\
	&\le C\epsilon^{-1}\abs{p_0}^2\bigg(\sup_{s\in[0, 1-\epsilon^{4/3}]}\big(\epsilon^8\abs{\log\epsilon}^2a^{-2}+ \epsilon^6\big) \\
	&\qquad +\sup_{s\in[1-2\epsilon^{4/3}, 1-\epsilon^{4/3}]}\big(\epsilon^{-4/3}\epsilon^8a^2+ \epsilon^{-4/3}\epsilon^6a^4+ \epsilon^{-4}\epsilon^8a^6\big)\bigg)
	\le C\epsilon^5\abs{p_0}^2\,.
\end{aligned}
\end{equation}
Here we have integrated in $r$, used the spheroidal endpoint condition \eqref{eq:spheroidal}, and applied the $L^\infty$ bounds of Lemma \ref{lem:press_derivs}. 
We will also require the following bounds for $\mc{U}_2$ on the vessel surface $\Gamma_\epsilon$:
\begin{equation}
\begin{aligned}
	\|\mc{U}_2\|_{L^2(\Gamma_\epsilon)}^2 
	&\le C\int_0^{1-\epsilon^{4/3}}\int_0^{2\pi}\bigg(\epsilon^6a^{-4}(a^5p^{\rm SB}_s)_{ss}^2 +\epsilon^4a^{-2}(a^3p^{\rm SB}_s)_s^2 \bigg)\,\epsilon a\,d\theta ds\\
	&\qquad + \frac{C}{\epsilon^{8/3}}\int_{1-2\epsilon^{4/3}}^{1-\epsilon^{4/3}}\int_0^{2\pi}\bigg(\epsilon^6(a^3p^{\rm SB}_s)_s^2 + \epsilon^4(a^2+a^4)(a\,p^{\rm SB}_s)^2\bigg)\,\epsilon a\,d\theta ds \\
	&\le C\epsilon^{-1}\abs{p_0}^2\bigg(\sup_{s\in[0, 1-\epsilon^{4/3}]}\big(\epsilon^7\abs{\log\epsilon}^2a^{-3}+ \epsilon^5a^{-1}\big)\\
	&\qquad +\sup_{s\in[1-2\epsilon^{4/3}, 1-\epsilon^{4/3}]}\epsilon^{-4/3}\big(\epsilon^7a+\epsilon^5 a^3 \big)\bigg)
	\le C\epsilon^{10/3}\abs{p_0}^2\,,
\end{aligned}
\end{equation}
and along the bottom cross section $B_\epsilon$:
\begin{equation}
\begin{aligned}
		\|\mc{U}_2\|_{L^2(B_\epsilon)}^2 
		&\le C\int_0^{2\pi}\int_0^{\epsilon a}\bigg(\epsilon^6(a^5p^{\rm SB}_s)_{ss}^2 +\epsilon^4(a^3p^{\rm SB}_s)_s^2 \bigg)\bigg|_{s=0}\,r\,drd\theta \\
		&\le C\big(\epsilon^8\abs{\log\epsilon}^2+ \epsilon^6\big)\epsilon^{-1}\abs{p_0}^2
		\le C\epsilon^5\abs{p_0}^2 \,.
\end{aligned}
\end{equation}
Here we have used Lemma \ref{lem:press_derivs}, the lower bound $a(s)\ge a_0>0$ away from $s=1$, and the fact that $\p_s\phi_\epsilon$ is supported far away from $s=0$.

Next, using the curvilinear cylindrical coordinates of section \ref{sec:setup}, we may calculate 
\begin{equation}\label{eq:DeltaUet_f}
\begin{aligned}
	\Delta(\bm{U}\cdot\be_{\rm t}) &= \frac{1}{4\mu}\bigg[\frac{1}{1-r\wh\kappa}\bigg(\frac{1}{r}\p_r\big[(1-r\wh\kappa)\,2r^2\,\phi_\epsilon\,p^{\rm SB}_s\big] \bigg) \\
	&\qquad + \div\bigg(\frac{1}{1-r\wh\kappa}\p_s\big[\phi_\epsilon\,(r^2-(\epsilon a)^2)p^{\rm SB}_s\big]\be_{\rm t} \bigg)\bigg]\\
	&= \frac{\phi_\epsilon\,p^{\rm SB}_s}{\mu(1-r\wh\kappa)} - \div\bigg(\frac{1}{(1-r\wh\kappa)}\bigg[\frac{3}{2\mu}r^2\wh\kappa\,\phi_\epsilon\,p^{\rm SB}_s\,\be_r-  \p_s\big[\phi_\epsilon(r^2-(\epsilon a)^2)p^{\rm SB}_s\big]\be_{\rm t}\bigg] \bigg)\\
	&= \frac{\phi_\epsilon\,p^{\rm SB}_s}{\mu(1-r\wh\kappa)} + \div(\mc{U}_3)\,,
\end{aligned}
\end{equation}
where the vector-valued remainder term $\mc{U}_3$ satisfies 
\begin{equation}
\begin{aligned}
	\|\mc{U}_3\|_{L^2(\mc{V}_\epsilon)}^2 &\le C\int_0^{1-\epsilon^{4/3}}\int_0^{2\pi}\int_0^{\epsilon a}\epsilon^4a^{-2}\bigg(a^6(p^{\rm SB}_{ss})^2 +a^2(p^{\rm SB}_s)^2 \bigg)\,r\,drd\theta ds \\
	&\quad +  \frac{C}{\epsilon^{8/3}}\int_{1-2\epsilon^{4/3}}^{1-\epsilon^{4/3}}\int_0^{2\pi}\int_0^{\epsilon a}\epsilon^4a^2(a\,p^{\rm SB}_s)^2 \,r\,drd\theta ds \\
	&\le C\big(\epsilon^6+ \sup_{s\in[1-2\epsilon^{4/3},1-\epsilon^{4/3}]}\epsilon^{-4/3}\epsilon^6a^4\big)\epsilon^{-1}\abs{p_0}^2
	\le C\epsilon^5\abs{p_0}^2\,.
\end{aligned}
\end{equation}
Here we have used the $L^\infty$ bounds of Lemma \ref{lem:press_derivs} and the spheroidal endpoint condition \eqref{eq:spheroidal}.
We will also need the following estimates for $\mc{U}_3$ on $\Gamma_\epsilon$: 
\begin{equation}
\begin{aligned}
	\|\mc{U}_3\|_{L^2(\Gamma_\epsilon)}^2 &\le C\int_0^{1-\epsilon^{4/3}}\int_0^{2\pi}\epsilon^4a^{-2}(a\,p^{\rm SB}_s)^2 \, \epsilon a\,d\theta ds    \\
	&\le C\big(\sup_{s\in[0,1-\epsilon^{4/3}]}\epsilon^5a^{-1} \big)\epsilon^{-1}\abs{p_0}^2
	\le C\epsilon^{10/3}\abs{p_0}^2\,,
\end{aligned}
\end{equation}
and along the bottom cross section $B_\epsilon$: 
\begin{equation}
\begin{aligned}
	\|\mc{U}_3\|_{L^2(B_\epsilon)}^2 &\le C\int_0^{2\pi}\int_0^{\epsilon a}\epsilon^4\bigg((a^3p^{\rm SB}_{ss})^2 +(a\,p^{\rm SB}_s)^2 \bigg)\bigg|_{s=0}\,r\,drd\theta  \\
	&\le C\epsilon^6\,\epsilon^{-1}\abs{p_0}^2 
	\le C\epsilon^5\abs{p_0}^2\,.
\end{aligned}
\end{equation}

Combining \eqref{eq:divEU}, \eqref{eq:DeltaU}, and \eqref{eq:DeltaUet_f}, we may thus write 
\begin{equation}\label{eq:bulk_expr_f}
\begin{aligned}
	-2\mu\,\div\big(\E(\bm{U})\big) + \nabla p^{\rm SB} 
	&= -\frac{\phi_\epsilon \,p^{\rm SB}_s}{1-r\wh\kappa}\,\be_{\rm t} + \frac{p^{\rm SB}_s}{1-r\wh\kappa}\,\be_{\rm t} -\mu\big(\nabla\mc{U}_0 + \mc{U}_1 + \div(\mc{U}_2) + \div(\mc{U}_3)\,\be_{\rm t}\big) \\
	&= \frac{(1-\phi_\epsilon)\,p^{\rm SB}_s}{1-r\wh\kappa}\be_{\rm t} -\mu\big(\nabla\mc{U}_0 + \mc{U}_1 + \div(\mc{U}_2) + \div(\mc{U}_3)\,\be_{\rm t}\big)\,,
\end{aligned}
\end{equation}
where we have used the form \eqref{eq:grad_div} of the curvilinear gradient to calculate $\nabla p^{\rm SB}$.

Given the representation \eqref{eq:bulk_expr_f}, we may consider any $\bw\in H^1(\mc{V}_\epsilon)$ with $\bw-(\bw\cdot\bm{n})\bm{n}=0$ on $\Gamma_\epsilon$ and write  
\begin{equation}
\begin{aligned}
	&\int_{\mc{V}_\epsilon}\bigg(2\mu\,\E(\bm{U}):\E(\bw) -\big(p^{\rm SB}-p_0\big)\,\div\bw\bigg)\,d\bx = -\int_{\Gamma_\epsilon}\big(p^{\rm SB}- p_0\big)(\bw\cdot\bm{n}) \,dS_\epsilon \\
	& 
	+\underbrace{\int_{\mc{V}_\epsilon}\bigg(\frac{(1-\phi_\epsilon)p^{\rm SB}_s}{1-r\wh\kappa}\be_{\rm t}- \mu\,\mc{U}_1\bigg)\cdot\bw\,d\bx}_{W_1}
	+\underbrace{\mu\int_{\mc{V}_\epsilon} \big(\mc{U}_0\,\div\bw + \mc{U}_2:\nabla\bw + \mc{U}_3\cdot\nabla(\bw\cdot\be_{\rm t}) \big)\,d\bx}_{W_2}  \\
	& + \underbrace{\int_{\Gamma_\epsilon}2\mu(\E(\bm{U})\bm{n})\cdot\bm{n}\,(\bw\cdot\bm{n}) \,dS_\epsilon}_{W_3} 
	-\underbrace{\mu \int_{\Gamma_\epsilon} \bigg(\big(\mc{U}_0+(\mc{U}_2\bm{n})\cdot\bm{n}\big)(\bw\cdot\bm{n})+\mc{U}_3\cdot\bm{n}(\bw\cdot\be_{\rm t}) \bigg)\,dS_\epsilon}_{W_4}
	 \\
	&  
	-\underbrace{\mu\int_{B_\epsilon}\bigg((2\E(\bm{U})+\mc{U}_2)\bm{n} + \big(\mc{U}_0-\frac{p^{\rm SB}}{\mu}\big)\bm{n} + (\mc{U}_3\cdot\bm{n})\be_{\rm t}\bigg)\cdot\bw \,d\bx}_{W_5}\,.
\end{aligned}
\end{equation}
Here each of the remainder terms $W_i$ are small in $\epsilon$ and may be bounded in turn as follows. 

We may first estimate
\begin{equation}
\begin{aligned}
	\abs{W_1} &\le C\bigg[\bigg(\int_{1-2\epsilon^{4/3}}^{1-\epsilon^{4/3}} \int_0^{2\pi}\int_0^{\epsilon a}\abs{1-\phi_\epsilon}^2\abs{p^{\rm SB}_s}^2\,r\,drd\theta ds\bigg)^{1/2}+ \norm{\mc{U}_1}_{L^2(\mc{V}_\epsilon)}\bigg]\norm{\bw}_{L^2(\mc{V}_\epsilon)} \\
	&\le C\big(\epsilon^{5/3}\norm{a\,p^{\rm SB}_s}_{L^\infty(0,1)} +\epsilon^3\abs{p_0}\big)\norm{\bw}_{L^2(\mc{V}_\epsilon)}\\
	&\le C\epsilon^{7/6}\abs{p_0}\norm{\bw}_{L^2(\mc{V}_\epsilon)}
	\le C\epsilon^{13/6}\abs{p_0}\norm{\nabla\bw}_{L^2(\mc{V}_\epsilon)}\,,
\end{aligned}
\end{equation}
where we have used the Poincar\'e inequality implied by Lemma \ref{lem:korn} for $\bw$, since $\bw-(\bw\cdot\bm{n})\bm{n}=0$ on the vessel surface $\Gamma_\epsilon$. 
Note that this is one of two terms in the error analysis which limit the extent of ${\rm supp}(1-\phi_\epsilon)$: essentially, this term tells us how \emph{small} the tip region must be in order to take advantage of the geometric smallness of the tip.

Next, using the above bounds for each $\mc{U}_i$ within the vessel $\mc{V}_\epsilon$, we may estimate
 \begin{equation}
 \begin{aligned}
 	|W_2| &\le C\epsilon^{5/2}\abs{p_0}\norm{\nabla\bw}_{L^2(\mc{V}_\epsilon)}\,.
 \end{aligned}
 \end{equation}
In addition, using the expression \eqref{eq:Udef} for $\bm{U}$ and the form \eqref{eq:normal_param} of the normal vector $\bm{n}$ along the vessel surface $\Gamma_\epsilon$, we may bound 
\begin{equation}
\begin{aligned}
	\|(\E(\bm{U})\bm{n})\cdot\bm{n}\|_{L^2(\Gamma_\epsilon)}^2
	&\le C\|\p_r(\bm{U}\cdot\bm{n})\|_{L^2(\Gamma_\epsilon)}^2 + C\|\epsilon a'\p_s(\bm{U}\cdot\bm{n})\|_{L^2(\Gamma_\epsilon)}^2\\
	&\le C\int_0^{1-\epsilon^{4/3}}\int_0^{2\pi}\bigg(\epsilon^4a^{-2}(a^3p^{\rm SB}_s)_s^2 + \epsilon^8|a'|^2a^{-4}(a^5p^{\rm SB}_s)_{ss}^2 \bigg)\,\epsilon a\,d\theta ds\\
	&\qquad + \frac{C}{\epsilon^{8/3}}\int_{1-2\epsilon^{4/3}}^{1-\epsilon^{4/3}}\int_0^{2\pi}
	\bigg(\epsilon^4a^2(a\,p^{\rm SB}_s)^2 + \epsilon^8|a'|^2(a^3p^{\rm SB}_s)_{s}^2 \bigg)\,\epsilon a\,d\theta ds \\
	&\qquad + \frac{C}{\epsilon^{16/3}}\int_{1-2\epsilon^{4/3}}^{1-\epsilon^{4/3}}\int_0^{2\pi}\epsilon^8|a'|^2a^4(a\,p^{\rm SB}_s)^2 \,\epsilon a\,d\theta ds \\
	&\le C\epsilon^{-1}\abs{p_0}^2\bigg(\sup_{s\in[0,1-\epsilon^{4/3}]} (\epsilon^5a^{-1} +\epsilon^9\abs{\log\epsilon}^2a^{-5})\\
	&\quad  + \sup_{s\in[1-2\epsilon^{4/3},1-\epsilon^{4/3}]}\big(\epsilon^{-4/3}(\epsilon^5a^3+\epsilon^9a^{-1})+\epsilon^{-4}\epsilon^9a^3\big)\bigg)\\
	&\le C\epsilon^{10/3}\abs{p_0}^2\,.
\end{aligned}
\end{equation}
We may thus bound the third remainder term $W_3$ as
 \begin{equation}
 \begin{aligned}
 	|W_3| &\le C\|(\E(\bm{U})\bm{n})\cdot\bm{n}\|_{L^2(\Gamma_\epsilon)}\norm{\bw\cdot\bm{n}}_{L^2(\Gamma_\epsilon)}
 	\le  C\epsilon^{5/3}\abs{p_0}\norm{\bw\cdot\bm{n}}_{L^2(\Gamma_\epsilon)}\,.
 \end{aligned}
 \end{equation}
Similarly, we may bound the remainder term $W_4$ along $\Gamma_\epsilon$ as 
 \begin{equation}
 \begin{aligned}
 	|W_4|&\le C\epsilon^{5/3}\abs{p_0}\norm{\bw\cdot\bm{n}}_{L^2(\Gamma_\epsilon)}\,.
 \end{aligned}
 \end{equation}

Finally, using the expression \eqref{eq:Udef} to calculate $\nabla\bm{U}$ on the bottom cross section $B_\epsilon$, we may bound 
\begin{equation}
\begin{aligned}
	\|\nabla\bm{U}\|_{L^2(B_\epsilon)}^2 &\le C\int_0^{2\pi}\int_0^{\epsilon a}\bigg(\epsilon^6(p^{\rm SB}_{sss})^2+\epsilon^4(p^{\rm SB}_{ss})^2+\epsilon^2(p^{\rm SB}_s)^2 \bigg)\bigg|_{s=0}\,r\,drd\theta
	\le C\epsilon^4\abs{p_0}^2\,,
\end{aligned}
\end{equation}
by a combination of the $L^2$ and $L^\infty$ bounds of Lemma \ref{lem:press_derivs}. Here we have again used that the radius function $a$ is bounded away from zero near $B_\epsilon$. 
Then, using the above bounds for $\mc{U}_0$, $\mc{U}_2$, and $\mc{U}_3$, we may estimate the remainder term $W_5$ coming from the bottom of the vessel as 
 \begin{equation}
 \begin{aligned}
 	|W_5|&\le C\epsilon^2\abs{p_0}\norm{\bw}_{L^2(B_\epsilon)}
 	\le C\epsilon^3\abs{p_0}\norm{\nabla\bw}_{L^2(\mc{V}_\epsilon)}\,.
 \end{aligned}
 \end{equation}
 Here we have used the trace inequality (Lemma \ref{lem:trace_ineqs}) on $B_\epsilon$ along with the Poincar\'e inequality implied by Lemma \ref{lem:korn} for $\bw$.

Combining the estimates for each $W_i$, we obtain Proposition \ref{prop:almostStokes}.
\end{proof}


\subsection{Boundary residuals}
Given that the velocity ansatz $\bm{U}(\bx)$ given by \eqref{eq:Udef} satisfies Proposition \ref{prop:almostStokes}, it remains to compare the triple $(\bm{U},p^{\rm SB}(s),q^{\rm SB})$ with the full 3D-3D solution $(\bu^\epsilon,p^\epsilon,q^\epsilon)$. To do so, we will insert $(\bm{U},p^{\rm SB},q^{\rm SB})$ into the bilinear form $\mc{B}:=\mc{A}+b$ arising in the weak form (Definition \ref{def:weak3D}) of the 3D-3D system. 
We show the following proposition.
\begin{proposition}\label{prop:Adiff}
	Given pressure data $p_0\in \R$, let $(\bu^\epsilon,p^\epsilon,q^\epsilon)$ be the solution to the 3D-3D system \eqref{eq:3D3D_1}-\eqref{eq:3D3D_6}, and let $(\bm{U},p^{\rm SB}(s),q^{\rm SB})$ be given by the 1D model \eqref{eq:qSB}-\eqref{eq:pSB} with velocity ansatz $\bm{U}$ as in \eqref{eq:Udef}.
	Let $\mc{B}:=\mc{A}+b$ be the bilinear form corresponding to the 3D-3D system in Definition \ref{def:weak3D}. Consider any $g\in D^{1,2}(\Omega_\epsilon)$ and any $\bw \in H^1(\mc{V}_\epsilon)$ with $\bw-(\bw\cdot\bm{n})\bm{n}=0$ on $\Gamma_\epsilon$. 

	The differences $(\wh{\bu},\wh p,\wh q):=(\bu^\epsilon-\bm{U},p^\epsilon-p^{\rm SB}(s),q^\epsilon-q^{\rm SB})$ satisfy
\begin{equation}\label{eq:propAdiff_est}
\begin{aligned}
	\mc{B}\big((\wh{\bu},\wh q;\wh p),(\bw,g) \big)&:=
	\int_{\mc{V}_\epsilon}\bigg(2\mu\,\E(\wh{\bu}):\E(\bw)-\wh p\,\div\bw\bigg)\,d\bx + \frac{1}{\epsilon^3\kappa}\int_{\Gamma_\epsilon}(\wh{\bu}\cdot\bm{n})(\bw\cdot\bm{n})\,dS_\epsilon\\
    &+ \int_{\Gamma_\epsilon}\wh q\,(\bw\cdot\bm{n})\,dS_\epsilon +\zeta\epsilon^4\int_{\Omega_\epsilon}\nabla \wh q\cdot\nabla g\,d\bx - \int_{\Gamma_\epsilon}(\wh{\bu}\cdot\bm{n})\,g\,dS_\epsilon = E_B(\bw,g)\,,
\end{aligned}
\end{equation}
where the right hand side satisfies 
\begin{equation}
\begin{aligned}
	|E_B| &\le C\abs{p_0}\big(\epsilon^{13/6}\norm{\nabla\bw}_{L^2(\mc{V}_\epsilon)}+\epsilon^{2/3}\norm{\bw\cdot\bm{n}}_{L^2(\Gamma_\epsilon)} +
	\epsilon^{9/2}\abs{\log\epsilon}^{3/2}\norm{\nabla g}_{L^2(\Omega_\epsilon)}\big)\,.
\end{aligned}
\end{equation}
\end{proposition}

\begin{proof}[Proof of Proposition \ref{prop:Adiff}]
Given any pair $(\bw,g)$ as in the proposition statement, we may write out $\mc{B}\big((\bm{U},q^{\rm SB};p^{\rm SB}-p_0),(\bw,g) \big)$ as 
\begin{equation}\label{eq:Aform}
\begin{aligned}
    \mc{B}\big((\bm{U},q^{\rm SB};p^{\rm SB}-p_0),(\bw,g) \big)&=\int_{\mc{V}_\epsilon}\bigg(2\mu\,\E(\bm{U}):\E(\bw)- \big(p^{\rm SB}(s)-p_0\big)\,\div\bw\bigg)\,d\bx \\
    &\quad + \frac{1}{\epsilon^3\kappa}\int_{\Gamma_\epsilon}(\bm{U}\cdot\bm{n})(\bw\cdot\bm{n})\,dS_\epsilon + \int_{\Gamma_\epsilon}q^{\rm SB}\,(\bw\cdot\bm{n})\,dS_\epsilon\\
    &\quad  +\zeta\epsilon^4\int_{\Omega_\epsilon}\nabla q^{\rm SB}\cdot\nabla g\,d\bx - \int_{\Gamma_\epsilon}(\bm{U}\cdot\bm{n})\,g\,dS_\epsilon \,.
\end{aligned}
\end{equation} 
We will consider the form of each term in $\mc{B}\big((\bm{U},q^{\rm SB};p^{\rm SB}-p_0),(\bw,g) \big)$ individually. 

Immediately by Proposition \ref{prop:almostStokes}, we may write
\begin{equation}\label{eq:E0terms}
\begin{aligned}
 	\int_{\mc{V}_\epsilon}&\bigg(2\mu\,\E(\bm{U}):\E(\bw) -\big(p^{\rm SB}(s)-p_0\big)\,\div\bw \bigg)\,d\bx \\
 	&= -\int_0^1p^{\rm SB}(s)\int_0^{2\pi}(\bw\cdot\bm{n}) \,\mc{J}_\epsilon(s,\theta')\,d\theta'\,ds +\int_{\Gamma_\epsilon}p_0\,(\bw\cdot\bm{n}) \,dS_\epsilon + E_0\\
 	&= -\frac{1}{2\pi}\int_0^{2\pi}\int_0^1p^{\rm SB}(s)\bigg(\int_0^{2\pi}(\bw\cdot\bm{n}) \,d\theta'\bigg)\,\mc{J}_\epsilon(s,\theta)\,ds\,d\theta \\
 	&\qquad +\int_{\Gamma_\epsilon}p_0\,(\bw\cdot\bm{n}) \,dS_\epsilon  +E_0 + E_1\,,\\
 	|E_0|&\le C\abs{p_0}\big(\epsilon^{13/6}\norm{\nabla\bw}_{L^2(\mc{V}_\epsilon)}+ \epsilon^{5/3}\norm{\bw\cdot\bm{n}}_{L^2(\Gamma_\epsilon)} \big)\,,
\end{aligned}
\end{equation}
where the additional error term $E_1$ from moving the Jacobian factor $\mc{J}_\epsilon$ satisfies
\begin{equation}
	|E_1| \le C\epsilon|\Gamma_\epsilon|^{1/2}\norm{p^{\rm SB}}_{L^2(0,1)}\norm{\bw\cdot\bm{n}}_{L^2(\Gamma_\epsilon)}
	\le C\epsilon^{3/2}\abs{p_0}\norm{\bw\cdot\bm{n}}_{L^2(\Gamma_\epsilon)}\,.
\end{equation}

Furthermore, recalling the expression \eqref{eq:Udef} for $\bm{U}$ and using the form \eqref{eq:free_jacfac} of $\mc{J}_\epsilon$, we may evaluate 
\begin{equation}\label{eq:UdotN}
\begin{aligned}
	(\bm{U}\cdot\bm{n})\big|_{\Gamma_\epsilon}\,\mc{J}_\epsilon(s,\theta) &= \bigg(\frac{\phi_\epsilon\epsilon^4}{16\mu}\big(a^4p^{\rm SB}_{ss} + 4a^3a'p^{\rm SB}_s\big)
	+ \frac{\epsilon^4a^3}{16\mu}(\p_s\phi_\epsilon) a\,p^{\rm SB}_s\bigg) \frac{\sqrt{(1-\epsilon a\wh\kappa)^2+(\epsilon a')^2}}{\sqrt{1+(\epsilon a')^2}}\\
	&= \phi_\epsilon\frac{\epsilon^4}{16\mu} \p_s\big(a^4\p_sp^{\rm SB}\big) + \mc{P}_1(s,\theta)\,, \\
	\norm{\mc{P}_1}_{L^2((0,1)\times(0,2\pi))}&\le C\big(\epsilon^5\norm{\p_s(a^4\p_sp^{\rm SB})}_{L^2(0,1)} \\
	&\qquad + \sup_{s\in[1-2\epsilon^{4/3},1-\epsilon^{4/3}]}\epsilon^{-2/3}\epsilon^4a^3\norm{a\,p^{\rm SB}_s}_{L^\infty(0,1)}\big)
	\le C\epsilon^{29/6}\abs{p_0}\,.
\end{aligned}
\end{equation}
Here we have again used Lemma \ref{lem:press_derivs} and that $\p_s\phi_\epsilon$ is supported only for $s\in[1-2\epsilon^{4/3},1-\epsilon^{4/3}]$.
We may then calculate 
\begin{equation}\label{eq:E2terms}
\begin{aligned}
\frac{1}{\epsilon^3\kappa}\int_{\Gamma_\epsilon}(\bm{U}\cdot\bm{n})(\bw\cdot\bm{n})\,dS_\epsilon &=  \frac{\epsilon}{16\kappa\mu}\int_0^1\p_s\big(a^4\p_s p^{\rm SB}\big)\,\phi_\epsilon\int_0^{2\pi}(\bw\cdot\bm{n})\,d\theta ds  + E_2\,, \\
%
|E_2|&\le C\epsilon^{-3}\norm{\mc{P}_1}_{L^2}\norm{\bw\cdot\bm{n}}_{L^2(\Gamma_\epsilon)} 
\le C\epsilon^{11/6}\abs{p_0}\norm{\bw\cdot\bm{n}}_{L^2(\Gamma_\epsilon)}\,.
\end{aligned}
\end{equation}
%
%
Using the equation \eqref{eq:pSB_eqn} satisfied by $p^{\rm SB}(s)$, we may next calculate 
\begin{equation}\label{eq:E3terms}
\begin{aligned}
	\int_{\Gamma_\epsilon}(\bm{U}\cdot\bm{n})\,g\,dS_\epsilon &= \frac{\epsilon^4}{16\mu}\int_0^1\p_s(a^4\p_s p^{\rm SB})\,\phi_\epsilon\int_0^{2\pi}g\,d\theta\,ds  - E_3 \\
	&= \frac{\epsilon^4\kappa}{2\pi} \int_0^1a(s)\bigg(\int_0^{2\pi}(p^{\rm SB}(s)-q^{\rm SB})\,d\theta\bigg)\bigg(\phi_\epsilon\int_0^{2\pi}g\,d\theta\bigg)\,ds - E_3 \\
	&= \epsilon^3\kappa \int_0^1\int_0^{2\pi}\bigg(\frac{1}{2\pi}\int_0^{2\pi}(p^{\rm SB}(s)-q^{\rm SB})\,d\theta\bigg)\,\phi_\epsilon\,g\,\mc{J}_\epsilon(s,\theta)\,d\theta ds -E_3- E_4\\
	&= \epsilon^3\kappa \int_0^1\int_0^{2\pi}\big(p^{\rm SB}(s)-q^{\rm SB}\big)\,\phi_\epsilon\,g\,\mc{J}_\epsilon(s,\theta)\,d\theta ds -E_3- E_4-E_5\\
	&= \epsilon^3\kappa \int_{\Gamma_\epsilon}\big(p^{\rm SB}(s)-q^{\rm SB}\big)\,g\,dS_\epsilon -E_3- E_4-E_5 -E_6\,,
\end{aligned}
\end{equation}
where in the final line the cutoff $\phi_\epsilon$ has been absorbed into $E_6$. Here the remainder terms $E_3$ and $E_4$ satisfy 
\begin{equation}
\begin{aligned}
	|E_3|&\le C\norm{\mc{P}_1}_{L^2}\norm{g}_{L^2(\Gamma_\epsilon)} \le C\epsilon^{29/6}\abs{p_0}\norm{g}_{L^2(\Gamma_\epsilon)}
	\le C\epsilon^{16/3}\abs{\log\epsilon}^{1/2}\abs{p_0}\norm{\nabla g}_{L^2(\Omega_\epsilon)}\,,\\
	|E_4|&\le C\epsilon^4\norm{p^{\rm SB}-q^{\rm SB}}_{L^2(\Gamma_\epsilon)}\norm{g}_{L^2(\Gamma_\epsilon)} \le C\epsilon^{9/2}\abs{p_0}\norm{g}_{L^2(\Gamma_\epsilon)}\\
	&\le 
	C\epsilon^5\abs{\log\epsilon}^{1/2}\abs{p_0}\norm{\nabla g}_{L^2(\Omega_\epsilon)}\,,
\end{aligned}
\end{equation}
where we have used the 3D-1D energy estimate \eqref{eq:energy3D1D} for $p-q$ and the 3D-1D to 1D error estimate (Theorem \ref{thm:3D1Dto1D}) to bound $p^{\rm SB}-q^{\rm SB}$ on $\Gamma_\epsilon$. We have then applied the exterior trace inequality of Lemma \ref{lem:trace_ineqs} for $g$. 

The remainder term $E_5$ in \eqref{eq:E3terms} comes from replacing the $\theta$-integral of the exterior pressure $q^{\rm SB}$ by $q^{\rm SB}$ itself, and may be estimated using Lemma \ref{lem:qtheta} as 
\begin{equation}
\begin{aligned}
	|E_5| &\le C\epsilon^3\norm{q^{\rm SB}-\frac{1}{2\pi}\int_0^{2\pi}q^{\rm SB}\,d\theta}_{L^2(\Gamma_\epsilon)}\norm{g}_{L^2(\Gamma_\epsilon)}
	\le C\epsilon^4\abs{\log\epsilon}\abs{p_0}\norm{g}_{L^2(\Gamma_\epsilon)} \\
	&\le C\epsilon^{9/2}\abs{\log\epsilon}^{3/2}\abs{p_0}\norm{\nabla g}_{L^2(\Omega_\epsilon)}\,.
\end{aligned}
\end{equation}
Here we have again used the exterior trace inequality of Lemma \ref{lem:trace_ineqs} for $g$. 
The final term $E_6$ is due to absorbing the cutoff $\phi_\epsilon$ and may be estimated as 
\begin{equation}
\begin{aligned}
	|E_6| &\le C\epsilon^3\int_{1-2\epsilon^{4/3}}^1\int_0^{2\pi}\abs{1-\phi_\epsilon}\abs{p^{\rm SB}(s)-q^{\rm SB}} \,\abs{g}\,dS_\epsilon \\
	&\le C\epsilon^{7/2}\abs{p_0}\norm{g}_{L^2(\Gamma_\epsilon\cap\,{\rm supp}(1-\phi_\epsilon))} 
	\le C\epsilon^{19/4}\abs{\log\epsilon}^{1/2}\abs{p_0}\norm{\nabla g}_{L^2(\Omega_\epsilon)}\,,
\end{aligned}
\end{equation}
where here we have used the free end region trace inequality of Lemma \ref{lem:trace_ineqs}.

Next, again using the near-$\theta$-independence of $q^{\rm SB}$ on $\Gamma_\epsilon$ due to Lemma \ref{lem:qtheta}, we note that 
\begin{equation}\label{eq:E7terms}
\begin{aligned}
	\int_{\Gamma_\epsilon}q^{\rm SB}\,(\bw\cdot\bm{n})\,dS_\epsilon &= 
	\frac{1}{2\pi}\int_0^1\int_0^{2\pi}\bigg(\int_0^{2\pi}q^{\rm SB}\,d\theta'\bigg)\,(\bw\cdot\bm{n})\,\mc{J}_\epsilon(s,\theta)\,d\theta ds\\
	&\qquad +\int_{\Gamma_\epsilon}\bigg(q^{\rm SB}-\frac{1}{2\pi}\int_0^{2\pi}q^{\rm SB}\,d\theta'\bigg)\,(\bw\cdot\bm{n})\,dS_\epsilon \\
	&= \frac{1}{2\pi}\int_0^1\int_0^{2\pi}q^{\rm SB}\,\bigg(\int_0^{2\pi}(\bw\cdot\bm{n})\,d\theta\bigg)\,\mc{J}_\epsilon(s,\theta')\,d\theta' ds +E_7\,,
\end{aligned}
\end{equation}
where the remainder term $E_7$ satisfies
\begin{equation}
\begin{aligned}
	|E_7|&\le C\bigg(\norm{q^{\rm SB}-\frac{1}{2\pi}\int_0^{2\pi}q^{\rm SB}\,d\theta}_{L^2(\Gamma_\epsilon)}+\epsilon\norm{q^{\rm SB}}_{L^2(\Gamma_\epsilon)}\bigg)\norm{\bw\cdot\bm{n}}_{L^2(\Gamma_\epsilon)} \\
	&\le C\big(\epsilon\abs{\log\epsilon}\abs{p_0} + \epsilon\norm{p^{\rm SB}-q^{\rm SB}}_{L^2(\Gamma_\epsilon)}+\epsilon\abs{\Gamma_\epsilon}^{1/2}\norm{p^{\rm SB}}_{L^2(0,1)}\big)\norm{\bw\cdot\bm{n}}_{L^2(\Gamma_\epsilon)}\\
	&\le C\epsilon\abs{\log\epsilon}\abs{p_0}\norm{\bw\cdot\bm{n}}_{L^2(\Gamma_\epsilon)}\,.
\end{aligned}
\end{equation} 
Here we have also used the 3D-1D energy estimate \eqref{eq:energy3D1D} along with the 3D-1D to 1D error estimate of Theorem \ref{thm:3D1Dto1D} to bound $p^{\rm SB}-q^{\rm SB}$ on $\Gamma_\epsilon$.

Finally, we may reinsert the cutoff $\phi_\epsilon$ into the terms coming from \eqref{eq:E0terms} and \eqref{eq:E7terms} as
\begin{equation}\label{eq:E8terms}
\begin{aligned}
	&\frac{1}{2\pi}\int_{\Gamma_\epsilon}(q^{\rm SB}-p^{\rm SB})\int_0^{2\pi}(\bw\cdot\bm{n})\,d\theta\,dS_\epsilon = \frac{1}{2\pi}\int_{\Gamma_\epsilon}(q^{\rm SB}-p^{\rm SB}(s))\,\phi_\epsilon\int_0^{2\pi}(\bw\cdot\bm{n})\,d\theta\,dS_\epsilon + E_8\,,\\
	&\quad \abs{E_8}\le C\norm{\bw\cdot\bm{n}}_{L^2(\Gamma_\epsilon)}\bigg(\int_{1-2\epsilon^{4/3}}^1\int_0^{2\pi}\abs{1-\phi_\epsilon}^2\abs{q^{\rm SB}-p^{\rm SB}}^2\,dS_\epsilon \bigg)^{1/2}\\
	&\quad\qquad\le C\norm{\bw\cdot\bm{n}}_{L^2(\Gamma_\epsilon)}\bigg(\epsilon^{5/4}\abs{\log\epsilon}\norm{\nabla q^{\rm SB}}_{L^2(\Omega_\epsilon)} + \epsilon^{7/6}\sup_{s\in[1-2\epsilon^{4/3},1]}a(s) \norm{p^{\rm SB}}_{L^\infty(0,1)}\bigg)\\
	&\quad\qquad\le C\norm{\bw\cdot\bm{n}}_{L^2(\Gamma_\epsilon)}\big(\epsilon^{5/4}\abs{\log\epsilon} + \epsilon^{11/6}\epsilon^{-1/2}\big)\abs{p_0}
	\le C\epsilon^{5/4}\abs{\log\epsilon}\abs{p_0}\norm{\bw\cdot\bm{n}}_{L^2(\Gamma_\epsilon)}\,.
\end{aligned}
\end{equation}
Here we have used the free endpoint trace inequality of Lemma \ref{lem:trace_ineqs}, the 3D-1D energy estimate \eqref{eq:energy3D1D}, and the convergence result of Theorem \ref{thm:3D1Dto1D} to bound $q^{\rm SB}$, and we have used the spheroidal endpoint condition \eqref{eq:spheroidal} and Lemma \ref{lem:press_derivs} to bound $p^{\rm SB}$.

Inserting each of \eqref{eq:E0terms}, \eqref{eq:E2terms}, \eqref{eq:E3terms}, \eqref{eq:E7terms}, and \eqref{eq:E8terms} into \eqref{eq:Aform}, for any $(\bw,g)\in H^1(\mc{V}_\epsilon)\times D^{1,2}(\Omega_\epsilon)$ with $\bw-(\bw\cdot\bm{n})\bm{n}=0$ on $\Gamma_\epsilon$, we thus obtain 
\begin{equation}\label{eq:AUqph}
\begin{aligned}
    &\mc{B}\big((\bm{U},q^{\rm SB};p^{\rm SB}-p_0),(\bw,g) \big)
     = \frac{\epsilon}{16\kappa\mu}\int_0^1\p_s(a^4\p_s p^{\rm SB})\int_0^{2\pi}(\bw\cdot\bm{n})\,d\theta\,ds  \\
    &\qquad + \frac{1}{2\pi}\int_{\Gamma_\epsilon}(q^{\rm SB}-p^{\rm SB}(s))\int_0^{2\pi}(\bw\cdot\bm{n})\,d\theta\,dS_\epsilon 
    + \int_{\Gamma_\epsilon}p_0(\bw\cdot\bm{n})\,dS_\epsilon  \\
    &\qquad +\zeta\epsilon^4\int_{\Omega_\epsilon}\nabla q^{\rm SB}\cdot\nabla g\,d\bx - \epsilon^3\kappa \int_{\Gamma_\epsilon}(p^{\rm SB}(s)-q^{\rm SB})\,g\,dS_\epsilon
    + \sum_{i=0}^8E_i\\
  	&\hspace{4.7cm} =  \int_{\Gamma_\epsilon}p_0\,(\bw\cdot\bm{n})\,dS_\epsilon+ \sum_{i=0}^8E_i + E_9\,.
\end{aligned}
\end{equation}  
Here we have used Lemma \ref{lem:weakdiff} with the pair $(\gamma(s),\wt g)=(-\frac{\phi_\epsilon}{2\pi}\int_0^{2\pi}\bw\cdot\bm{n}\,d\theta,\epsilon^4\zeta g)$, and the additional error term $E_9$ satisfies the bound 
\begin{equation}
\begin{aligned}
	\abs{E_9} &\le \abs{\wt E_1(\epsilon^4\zeta g)} + \abs{\wt E_2\bigg(-\frac{\phi_\epsilon}{2\pi}\int_0^{2\pi}\bw\cdot\bm{n}\,d\theta\bigg)}\\
	&\le C\epsilon^{9/2}\abs{\log\epsilon}^{3/2}\abs{p_0}\norm{\nabla g}_{L^2(\Omega_\epsilon)}
	+ C\epsilon^{3/2}\abs{p_0}\norm{\phi_\epsilon\int_0^{2\pi}\bw\cdot\bm{n}\,d\theta}_{L^2(0,1)} \\
	&\le C\epsilon^{9/2}\abs{\log\epsilon}^{3/2}\abs{p_0}\norm{\nabla g}_{L^2(\Omega_\epsilon)} \\
	&\qquad+C\epsilon\abs{p_0}\bigg(\sup_{s\in[0,1-\epsilon^{4/3}]}a^{-1}(s)\int_0^{1-\epsilon^{4/3}}\int_0^{2\pi}(\bw\cdot\bm{n})^2\,\epsilon a(s)\,d\theta ds \bigg)^{1/2}\\
	&\le C\epsilon^{9/2}\abs{\log\epsilon}^{3/2}\abs{p_0}\norm{\nabla g}_{L^2(\Omega_\epsilon)} + C\epsilon^{2/3}\norm{\bw\cdot\bm{n}}_{L^2(\Gamma_\epsilon)}\,.
\end{aligned}
\end{equation}
Here we have used that $a^{-1}(s)\ge C\epsilon^{-2/3}$ for $s\in[0,1-\epsilon^{4/3}]$, by the spheroidal endpoint condition \eqref{eq:spheroidal}. Note that the error term $E_9$ is a second limiting term, in addition to $W_1$, which dictates how small the tip region ${\rm supp}(1-\phi_\epsilon)$ must be. These two terms will be balanced later on by an error term which requires a lower bound on the extent of the tip region as well.

We now consider $(\bu^\epsilon,p^\epsilon,q^\epsilon)$ satisfying the 3D-3D system \eqref{eq:3D3D_1}-\eqref{eq:3D3D_6}, and define the differences $(\wh{\bu},\wh p,\wh q):=(\bu^\epsilon-\bm{U},p^\epsilon-p^{\rm SB}, q^\epsilon -q^{\rm SB})$. It will be convenient to consider $\wh p$ as 
\begin{equation}
\wh p= p^\epsilon - p_0 - \big(p^{\rm SB}(s) - p_0\big)\,.
\end{equation}
We may use $p^\epsilon - p_0$ in the 3D-3D system (Definition \ref{def:weak3D}) to lift the boundary condition at the vessel base to the vessel interior. Using $(\bu^\epsilon,p^\epsilon-p_0, q^\epsilon)$ in Definition \ref{def:weak3D}, we have 
\begin{equation}\label{eq:Auepsqepspeps}
\begin{aligned}
	\mc{B}\big((\bu^\epsilon,q^\epsilon;p^\epsilon-p_0),(\bw,g) \big)
	&= \int_{\Gamma_\epsilon}p_0(\bw\cdot\bm{n})\,dS_\epsilon
\end{aligned}
\end{equation}
for any $(\bw,g)$ as in \eqref{eq:AUqph}. 
Subtracting \eqref{eq:AUqph} from \eqref{eq:Auepsqepspeps}, we thus have that the differences $(\wh{\bu},\wh p,\wh q)$ satisfy
\begin{equation}\label{eq:Aform2}
\begin{aligned}
    \mc{B}\big((\wh{\bu},\wh q;\wh p),(\bw,g) \big)
    &=\int_{\mc{V}_\epsilon}\bigg(2\mu\,\E(\wh{\bu}):\E(\bw)- \wh p(s)\,\div\bw\bigg)\,d\bx + \frac{1}{\epsilon^3\kappa}\int_{\Gamma_\epsilon}(\wh{\bu}\cdot\bm{n})(\bw\cdot\bm{n})\,dS_\epsilon\\
    &\qquad  + \int_{\Gamma_\epsilon}\wh q\,(\bw\cdot\bm{n})\,dS_\epsilon +\zeta\epsilon^4\int_{\Omega_\epsilon}\nabla \wh q\cdot\nabla g\,d\bx - \int_{\Gamma_\epsilon}(\wh{\bu}\cdot\bm{n})\,g\,dS_\epsilon \\
    & = -\sum_{i=0}^9 E_i\,.
\end{aligned}
\end{equation} 
Using the bounds satisfied by the remainder terms $E_i$ shown above yields Proposition \ref{prop:Adiff}.
\end{proof}


\subsection{Error estimate}\label{subsec:error}
Equipped with Proposition \ref{prop:Adiff}, we may now proceed to the proof of Theorem \ref{thm:3D3Dto1D}. 
We begin by considering $\bm{d}\in H^1(\mc{V}_\epsilon)$ satisfying 
\begin{equation}\label{eq:d_eqn_f}
\begin{aligned}
	\div\bm{d} &= \wh p \qquad \text{in }\mc{V}_\epsilon\,, \\
	\bm{d}&=0 \;\text{ on } \Gamma_\epsilon\,,\qquad 
	\bm{d} = \frac{1}{|B_\epsilon|}\bigg(\int_{\mc{V}_\epsilon}\wh p\,d\bx\bigg)\bm{n}\; \text{ on }B_\epsilon\,,\\
	\norm{\nabla\bm{d}}_{L^2(\mc{V}_\epsilon)}&\le C\epsilon^{-1}\norm{\wh p}_{L^2(\mc{V}_\epsilon)}\,.
\end{aligned}
\end{equation}
Recall that such a $\bm{d}$ exists by Lemma \ref{lem:bogovskii}. Inserting $(\bm{d},0)$ into \eqref{eq:propAdiff_est}, we may then calculate
\begin{equation}
\begin{aligned}
	\int_{\mc{V}_\epsilon}\wh p^2\,d\bx &= \abs{-\int_{\mc{V}_\epsilon}2\mu\,\E(\wh{\bu}):\E(\bm{d})\,d\bx - E_B(\bm{d},0)}\\
	&\le C\big(\norm{\E(\wh{\bu})}_{L^2(\mc{V}_\epsilon)}+ \epsilon^{13/6}\abs{p_0}\big)\norm{\nabla\bm{d}}_{L^2(\mc{V}_\epsilon)} \\
	&\le C\big(\epsilon^{-1}\norm{\E(\wh{\bu})}_{L^2(\mc{V}_\epsilon)}+ \epsilon^{7/6}\abs{p_0}\big)\norm{\wh p}_{L^2(\mc{V}_\epsilon)}\,.
\end{aligned}
\end{equation}
Applying Young's inequality, we may thus bound
\begin{equation}\label{eq:phat_bd}
	\norm{\wh p}_{L^2(\mc{V}_\epsilon)} \le C\big(\epsilon^{-1}\norm{\E(\wh{\bu})}_{L^2(\mc{V}_\epsilon)}+ \epsilon^{7/6}\abs{p_0}\big)\,.
\end{equation}

We next aim to insert $(\wh\bu,\wh q)$ as a test function pair in \eqref{eq:propAdiff_est}, but we must first correct for the fact that $\bm{U}-(\bm{U}\cdot\bm{n})\bm{n}$ does not quite vanish on $\Gamma_\epsilon$.
This will be the factor which provides a lower bound for the extent of the tip region ${\rm supp}(1-\phi_\epsilon)$. The boundary condition satisfied by the ansatz $\bm{U}$ becomes progressively more erroneous as $s\to 1$, since the radial vector $\be_r$ becomes an increasingly bad approximation of the normal vector $\bm{n}$ to $\Gamma_\epsilon$.
 In particular, using the expression \eqref{eq:Udef} for $\bm{U}$ and the form \eqref{eq:normal_param} of the normal vector $\bm{n}$, along $\Gamma_\epsilon$ we may calculate 
\begin{equation}
\begin{aligned}
	\bm{U} -(\bm{U}\cdot\bm{n})\bm{n} &=
	-\frac{\epsilon a'}{\sqrt{1+(\epsilon a')^2}}(\bm{U}\cdot\be_r)\bm{n}^\perp \\
	&= -\frac{\epsilon a'}{16\mu\sqrt{1+(\epsilon a')^2}}\big((\epsilon a)^3(\phi_\epsilon p^{\rm SB}_s)_s+4\phi_\epsilon\epsilon^3a^2a'p^{\rm SB}_s \big)\bm{n}^\perp\,,
\end{aligned}
\end{equation}
where $\bm{n}^\perp=\frac{1}{\sqrt{1+(\epsilon a')^2}}(\epsilon a'\be_r+\be_{\rm t})$ and we have used that $\bm{U}\cdot\be_{\rm t}=0$ on $\Gamma_\epsilon$. Throughout $\mc{V}_\epsilon$, we thus define 
\begin{equation}\label{eq:vn_def}
	\bv_{\rm n}(r,\theta,s) = -\frac{\epsilon a'}{16\mu\sqrt{1+(\epsilon a')^2}}\big( r^3(\phi_\epsilon p^{\rm SB}_s)_s+4 r\epsilon^2aa'\phi_\epsilon\,p^{\rm SB}_s \big)\bm{n}^\perp\,.
\end{equation}
Note that $\bv_{\rm n}\cdot\bm{n}=0$ on $\Gamma_\epsilon$ and that, by Lemma \ref{lem:press_derivs}, $\bv_{\rm n}$ satisfies the bounds 
\begin{equation}\label{eq:vn_free_bds}
\begin{aligned}
	\norm{\bv_{\rm n}}_{L^2(\mc{V}_\epsilon)}^2 &\le C\int_0^{1-\epsilon^{4/3}}\int_0^{2\pi}\int_0^{\epsilon a}\epsilon^8a^{-2}\bigg((a^3p^{\rm SB}_{ss})^2+(a\,p^{\rm SB}_s)^2\bigg)\,r\,drd\theta ds\\
	&\qquad + \frac{C}{\epsilon^{8/3}}\int_{1-2\epsilon^{4/3}}^{1-\epsilon^{4/3}}\int_0^{2\pi}\int_0^{\epsilon a}\epsilon^8a^2(a\,p^{\rm SB}_s)^2\,r\,drd\theta ds 
	\le C\epsilon^9\abs{p_0}^2\,,\\
	\norm{\nabla\bv_{\rm n}}_{L^2(\mc{V}_\epsilon)}^2&\le C\int_0^{1-\epsilon^{4/3}}\int_0^{2\pi}\int_0^{\epsilon a}\epsilon^8a^{-6}\bigg((a^5p^{\rm SB}_{sss})^2+ (a^3p^{\rm SB}_{ss})^2+ (a\,p^{\rm SB}_s)^2\bigg)\,r\,drd\theta ds\\
	&\qquad + C\int_0^{1-\epsilon^{4/3}}\int_0^{2\pi}\int_0^{\epsilon a}\epsilon^6a^{-4}\bigg((a^3p^{\rm SB}_{ss})^2+ (a\,p^{\rm SB}_s)^2\bigg)\,r\,drd\theta ds\\
	&\qquad + \frac{C}{\epsilon^{8/3}}\int_{1-2\epsilon^{4/3}}^{1-\epsilon^{4/3}}\int_0^{2\pi}\int_0^{\epsilon a}\bigg(\epsilon^8a^{-2}(a^3p^{\rm SB}_{ss})^2 + \epsilon^6(a\,p^{\rm SB}_s)^2 \bigg)\,r\,drd\theta ds\\
	&\qquad + \frac{C}{\epsilon^{16/3}}\int_{1-2\epsilon^{4/3}}^{1-\epsilon^{4/3}}\int_0^{2\pi}\int_0^{\epsilon a}\epsilon^8a^2(a\,p^{\rm SB}_s)^2\,r\,drd\theta ds
	\le C\epsilon^{17/3}\abs{p_0}^2\,.
\end{aligned}
\end{equation}
Furthermore, using crucially the form of $\bm{n}^\perp$ and the curvilinear divergence \ref{eq:grad_div}, we may calculate 
\begin{equation}
\begin{aligned}
	\norm{\div\bv_{\rm n}}_{L^2(\mc{V}_\epsilon)}^2&\le 
	 C\int_0^{1-\epsilon^{4/3}}\int_0^{2\pi}\int_0^{\epsilon a}\epsilon^8a^{-6}\bigg((a^5p^{\rm SB}_{sss})^2+(a^3p^{\rm SB}_{ss})^2 + (a\,p^{\rm SB}_s)^2 \bigg)\,r\,drd\theta ds \\
	 &\qquad + \frac{C}{\epsilon^{8/3}}\int_{1-2\epsilon^{4/3}}^{1-\epsilon^{4/3}}\int_0^{2\pi}\int_0^{\epsilon a}\epsilon^8a^{-2}\bigg((a\,p^{\rm SB}_s)^2+ (a^3p^{\rm SB}_{ss})^2 \bigg)\,r\,drd\theta ds \\
	&\qquad + \frac{C}{\epsilon^{16/3}}\int_{1-2\epsilon^{4/3}}^{1-\epsilon^{4/3}}\int_0^{2\pi}\int_0^{\epsilon a} \epsilon^8a^2(a\,p^{\rm SB}_s)^2 \,r\,drd\theta ds \\
	 &\le C\epsilon^{-1}\abs{p_0}^2\bigg(\sup_{s\in[0,1-\epsilon^{4/3}]}\epsilon^{10}\abs{\log\epsilon}^2a^{-4}\\
	 &\qquad  + \sup_{s\in[1-2\epsilon^{4/3},1-\epsilon^{4/3}]}\big(\epsilon^{-4/3}\epsilon^{10} + \epsilon^{-4}\epsilon^{10}a^4 \big)\bigg) 
	  \le C\epsilon^{19/3}\abs{\log\epsilon}^2\abs{p_0}^2\,.
\end{aligned}
\end{equation}

Now we may use the pair $(\wh{\bu}-\bv_{\rm n},\wh q)$ as test functions $(\bw,g)$ in the bilinear form $\mc{B}\big((\wh{\bu},\wh q;\wh p),(\bw,g)\big)$. By Proposition \ref{prop:Adiff}, we obtain
\begin{equation}
\begin{aligned}
	&\int_{\mc{V}_\epsilon}2\mu\,\abs{\E(\wh{\bu})}^2\,d\bx + \frac{1}{\epsilon^3\kappa}\int_{\Gamma_\epsilon}(\wh{\bu}\cdot\bm{n})^2\,dS_\epsilon +\zeta\epsilon^4\int_{\Omega_\epsilon}\abs{\nabla \wh q}^2\,d\bx \\
	&\qquad = E_B(\wh{\bu}-\bv_{\rm n},\wh q)+ \int_{\mc{V}_\epsilon}\bigg(2\mu\,\E(\wh{\bu}):\E(\bv_{\rm n})+\wh p\,\div(\wh{\bu}-\bv_{\rm n})\bigg)\,d\bx\,,
\end{aligned}
\end{equation}
where we have used that $\bv_{\rm n}\cdot\bm{n}=0$ on $\Gamma_\epsilon$.
Using the estimates of Proposition \ref{prop:Adiff} for $E_B$ along with the Korn inequality (Lemma \ref{lem:korn}) for $\wh{\bu}-\bv_{\rm n}$, we may bound
\begin{equation}\label{eq:err1}
\begin{aligned}
	|E_B(\wh{\bu}-\bv_{\rm n},\wh q)| 
	&\le C\abs{p_0}\big(\epsilon^{13/6}\norm{\E(\wh{\bu}-\bv_{\rm n})}_{L^2(\mc{V}_\epsilon)}+\epsilon^{2/3}\norm{\wh{\bu}\cdot\bm{n}}_{L^2(\Gamma_\epsilon)} \\
	&\qquad +
	\epsilon^{9/2}\abs{\log\epsilon}^{3/2}\norm{\nabla \wh q}_{L^2(\Omega_\epsilon)}\big)\\
	&\le C\abs{p_0}\big(\epsilon^{13/6}\norm{\E(\wh{\bu})}_{L^2(\mc{V}_\epsilon)}+ \epsilon^{2/3}\norm{\wh{\bu}\cdot\bm{n}}_{L^2(\Gamma_\epsilon)} \\
	&\qquad  +
	\epsilon^{9/2}\abs{\log\epsilon}^{3/2}\norm{\nabla \wh q}_{L^2(\Omega_\epsilon)}\big) +\epsilon^5\abs{p_0}^2\,,
\end{aligned}
\end{equation}
where we have also used the bound \eqref{eq:vn_free_bds} for $\nabla\bv_{\rm n}$.
Furthermore, using the above estimates for $\wh p$, $\nabla\bv_{\rm n}$, $\div\bv_{\rm n}$, and the bound \eqref{eq:divUfbd} for $\div\wh{\bu}=-\div\bm{U}$, we have
\begin{equation}\label{eq:err2}
\begin{aligned}
	&\int_{\mc{V}_\epsilon}\bigg(2\mu\,\E(\wh{\bu}):\E(\bv_{\rm n})+\wh p\,\div(\wh{\bu}-\bv_{\rm n})\bigg)\,d\bx \\
	&\qquad \le C\norm{\E(\wh{\bu})}_{L^2(\mc{V}_\epsilon)}\norm{\nabla\bv_{\rm n}}_{L^2(\mc{V}_\epsilon)} + \norm{\wh p}_{L^2(\mc{V}_\epsilon)}\big(\norm{\div\wh{\bu}}_{L^2(\mc{V}_\epsilon)} + \norm{\div\bv_{\rm n}}_{L^2(\mc{V}_\epsilon)} \big) \\
	&\qquad\le C\epsilon^{17/6}\abs{p_0}\norm{\E(\wh{\bu})}_{L^2(\mc{V}_\epsilon)} + C\big(\epsilon^{-1}\norm{\E(\wh{\bu})}_{L^2(\mc{V}_\epsilon)} + \epsilon^{7/6}\abs{p_0} \big)\epsilon^{19/6}\abs{\log\epsilon}\abs{p_0}\\
	&\qquad \le C\epsilon^{13/6}\abs{\log\epsilon}\abs{p_0}\norm{\E(\wh{\bu})}_{L^2(\mc{V}_\epsilon)} + C\epsilon^{13/3}\abs{\log\epsilon}\abs{p_0}^2\,.
\end{aligned}
\end{equation}
Combining \eqref{eq:err1} and \eqref{eq:err2}, we obtain the bound
\begin{equation}
\begin{aligned}
	2\mu&\norm{\E(\wh{\bu})}_{L^2(\mc{V}_\epsilon)}^2 + \frac{1}{\epsilon^3\kappa}\norm{\wh{\bu}\cdot\bm{n}}_{L^2(\Gamma_\epsilon)}^2 + \zeta\epsilon^4\norm{\nabla\wh q}_{L^2(\Omega_\epsilon)}^2\\
	&\le C\bigg( \epsilon^{13/6}\abs{\log\epsilon}\abs{p_0}\norm{\E(\wh{\bu})}_{L^2(\mc{V}_\epsilon)}  + \epsilon^{2/3}\norm{\wh{\bu}\cdot\bm{n}}_{L^2(\Gamma_\epsilon)} + \epsilon^{9/2}\abs{\log\epsilon}^{3/2}\norm{\nabla \wh q}_{L^2(\Omega_\epsilon)}\bigg) \\
	&\qquad + C\epsilon^{13/3}\abs{\log\epsilon}\abs{p_0}^2\,.
\end{aligned}
\end{equation}
Using Young's inequality to split each of the first three terms, we obtain a bound for the differences $(\wh{\bu},\wh q)$:
\begin{equation}\label{eq:err_est1}
\begin{aligned}
	\norm{\E(\wh{\bu})}_{L^2(\mc{V}_\epsilon)}^2 + \frac{1}{\epsilon^3}\norm{\wh{\bu}\cdot\bm{n}}_{L^2(\Gamma_\epsilon)}^2 + \epsilon^4\norm{\nabla\wh q}_{L^2(\Omega_\epsilon)}^2 
	\le C\epsilon^{13/3}\abs{\log\epsilon}^2\abs{p_0}^2\,.
\end{aligned}
\end{equation}
By applying the Korn inequality of Lemma \ref{lem:korn} to $\wh{\bu}-\bv_{\rm n}$, we have 
\begin{equation}\label{eq:whbu_korn}
\begin{aligned}
	\norm{\E(\wh{\bu})}_{L^2(\mc{V}_\epsilon)} &\ge \norm{\E(\wh{\bu}-\bv_{\rm n})}_{L^2(\mc{V}_\epsilon)} - \norm{\E(\bv_{\rm n})}_{L^2(\mc{V}_\epsilon)}\\
	&\ge \norm{\nabla(\wh{\bu}-\bv_{\rm n})}_{L^2(\mc{V}_\epsilon)} - C\norm{\nabla\bv_{\rm n}}_{L^2(\mc{V}_\epsilon)}
	\ge \norm{\nabla \wh{\bu}}_{L^2(\mc{V}_\epsilon)} - C\epsilon^{17/6}\abs{p_0}\,.
\end{aligned}
\end{equation}
Combining the initial error bound \eqref{eq:err_est1} with the interior pressure estimate \eqref{eq:phat_bd} and the above bound \eqref{eq:whbu_korn}, upon dividing through by the volume factor $\abs{\mc{V}_\epsilon}\sim\epsilon^2$, we finally obtain the error estimate
\begin{equation}
\begin{aligned}
&\frac{1}{|\mc{V}_\epsilon|}\bigg(\epsilon^{-4}\norm{\wh{\bu}}_{L^2(\mc{V}_\epsilon)}^2+\epsilon^{-2}\norm{\nabla\wh{\bu}}_{L^2(\mc{V}_\epsilon)}^2 + \norm{\wh p}_{L^2(\mc{V}_\epsilon)}^2\bigg) + \frac{\epsilon^{-6}}{|\Gamma_\epsilon|}\norm{\wh{\bu}\cdot\bm{n}}_{L^2(\Gamma_\epsilon)}^2 \\
&\qquad  + \norm{\nabla\wh q}_{L^2(\Omega_\epsilon)}^2 \le C\epsilon^{1/3}\abs{\log\epsilon}^2\abs{p_0}^2\,,
\end{aligned}
\end{equation}
which yields Theorem \ref{thm:3D3Dto1D}.
\hfill\qedsymbol

%% file: appendix.tex

\appendix
\section{Inequalities}\label{sec:app}

\subsection{Korn inequality}\label{subsec:korn_inequality}
Here we verify the $\epsilon$-dependent scaling of Lemma \ref{lem:korn}. Throughout, we refer to the construction in \cite[Appendix A.1]{free_ends} for the more complicated setting of a homogeneous Korn inequality in the exterior of a thin domain. 
We begin by noting that for $\bw\in H^1(\R^3_+)$, we have 
\begin{equation}\label{eq:kornR}
  \norm{\nabla\bw}_{L^2(\R^3_+)} \le C\norm{\E(\bw)}_{L^2(\R^3_+)}\,,
\end{equation}
where $C$ clearly has no dependence on $\epsilon$. 
We will make use of \eqref{eq:kornR} by defining an extension operator
\begin{equation}\label{eq:extension_ops}
\begin{aligned}
  &Z_{\rm ex}: H^1(\mc{V}_\epsilon)\cap\big\{\bv-(\bv\cdot\bm{n})\bm{n}\big|_{\Gamma_\epsilon}=0\big\} \to H^1(\R^3_+)
\end{aligned}
\end{equation}
for which the symmetric gradient satisfies 
\begin{equation}\label{eq:desired_ext}
\|\E(Z_{\rm ex}\bv)\|_{L^2(\R^3_+)}\le C\norm{\E(\bv)}_{L^2(\mc{V}_\epsilon)}
\end{equation}
for $C$ independent of $\epsilon$. 
Note that the vanishing tangential boundary condition in \eqref{eq:extension_ops} is needed to avoid an additional $L^2$ term on the right hand side of \eqref{eq:desired_ext}.

To construct such an extension operator, it will be convenient to first subdivide the vessel $\mc{V}_\epsilon$ into segments which scale uniformly with the vessel radius $\epsilon a$. This requires some care due to the decay at the vessel endpoint. For $0\le s\le 1-\delta$, where $\delta$ is as in the definition \ref{def:rad} of an admissible radius function, we let
\begin{equation}\label{eq:sk_def}
  s_k = C\epsilon k\,, \quad k=0,\dots,\frac{1-\delta}{C\epsilon}\,,
\end{equation}
where $C$ is such that $\frac{1-\delta}{C\epsilon}$ is an integer. For the remaining interval $1-\delta\le s \le 1$, instead of the subdivision \eqref{eq:sk_def}, we will use a non-uniform partition of the interval with spacing chosen such that the resulting segments scale as uniformly as possible in $\epsilon a$. Let $K_\delta=\frac{1-\delta}{C\epsilon}$.
Starting from $s_{K_\delta} = 1-\delta$, we define the partition recursively as 
\begin{equation}\label{eq:nonuniform_partition}
  s_{k+1} = s_k + C_k \, \epsilon a(s_k)\,, \quad k=K_\delta,\dots, N_\epsilon-1\,,
\end{equation}
where the constants $C_k$ satisfy $C^{-1}\le C_k\le C$ independent of $\epsilon$, and the final point is placed at 
\begin{equation}
s_{N_\epsilon}=1-\epsilon^2\,.
\end{equation} 

We then consider a partition of unity $\{\varphi_k(s)\}_{k=1}^{N_\epsilon}$ on $[0,1]$ satisfying 
\begin{equation}\label{eq:partition}
  \varphi_k(s_k) = 1 \,, \quad {\rm supp}(\varphi_k)\subset [s_{k-1},s_{k+1}] 
\end{equation}
for the points $s_1$ to $s_{N_\epsilon-1}$. We define the endpoint segment more broadly to satisfy 
\begin{equation}
  \varphi_{N_\epsilon}(\bx) = \begin{cases}
  1 & \text{ if } \text{\rm dist}(\bx,\X(1-\epsilon^2))\le 2\epsilon^2\\
  0 & \text{ if } \text{\rm dist}(\bx,\X(1-\epsilon^2))\ge 3\epsilon^2\,.
  \end{cases}
\end{equation}
Note that we may bound $\abs{\nabla\varphi_k}\le \frac{C}{\epsilon a(s_{k+1})}$ for each $k=1,\dots,N_\epsilon-1$, with $\abs{\nabla\varphi_{N_\epsilon}}\le \frac{C}{\epsilon^2}$.

For $0\le s\le 1-\epsilon^2$, we additionally consider a cutoff function $\phi_\epsilon(r,s)$ in the radial direction about $\mc{V}_\epsilon$ satisfying
\begin{equation}\label{eq:radial_phieps}
  \phi_\epsilon(r,s) = \begin{cases}
    1\,,  & \epsilon a(s) \le r \le \frac{3\epsilon a(s)}{2}\\
    0\,,  & r \ge 2\epsilon a(s)\,,
  \end{cases}
\end{equation}
with smooth decay between such that $\abs{\nabla \phi_\epsilon}\le C(\epsilon a(s))^{-1}$.

For $k=1$ to $N_\epsilon-1$, we then define the following subdivision of $\mc{V}_\epsilon$ and the surrounding annular region into (curved) segments:
\begin{equation}
\mc{V}_k=\mc{V}_\epsilon \cap {\rm supp}(\varphi_k)\,, \quad
\mc{O}_k = {\rm supp}(\phi_\epsilon)\cap{\rm supp}(\varphi_k) \,.
\end{equation}
For the tip, we define
\begin{equation}
\mc{V}_{N_\epsilon} = {\rm supp}(\varphi_{N_\epsilon})\cap \mc{V}_\epsilon\,,\qquad 
  \mc{O}_{N_\epsilon} = {\rm supp}(\varphi_{N_\epsilon})\cap\Omega_\epsilon \,.
\end{equation}
Given $\bv\in H^1(\mc{V}_\epsilon)$, we let $\bv_k$ denote the restriction of $\bv$ to the segment $\mc{V}_k$.

On each segment $k=1,\dots,N_\epsilon-1$, we may define $E^{\rm ref}:H^1(\mc{V}_k)\to H^1(\mc{O}_k)$ as the following extension-by-reflection across $\Gamma_\epsilon$ in the purely radial direction:
\begin{equation}\label{eq:Erefv}
\begin{aligned}
  &E^{\rm ref}\bv_k(\X(s)+r\be_r(s,\theta)) \\
  &\quad= 
  \begin{cases}
  \bv_k(\X(s)+r\be_r(s,\theta))\,, & 0\le r \le \epsilon a(s) \\
  \phi_\epsilon(r,s)\,\bv_k(\X(s)+(2\epsilon a(s)-r)\be_r(s,\theta))\,,  & \epsilon a(s)\le r \le 2\epsilon a(s)\,, \\
  0\,, & r> 2\epsilon a(s)\,.
  \end{cases}
\end{aligned}
\end{equation}
Noting that the radial extension \eqref{eq:Erefv} is not suitable on the endpoint segment $\mc{V}_{N_\epsilon}$, we define an endpoint extension operator $E^{\rm end}:H^1(\mc{V}_{N_\epsilon})\to H^1(\mc{O}_{N_\epsilon})$ as follows. For any $\bx\in \mc{O}_{N_\epsilon}$, let $\bx^*\in \Gamma_\epsilon$ denote the projection of $\bx$ onto the vessel surface $\Gamma_\epsilon$ along the straight-line distance between $\bx$ and the vessel centerline segment $\{\X(s)\,:\, s\in[0,1-\epsilon^2]\}$.
Letting $\bm{n}^*=-\frac{\bx-\bx^*}{\abs{\bx-\bx^*}}$, we may write $\bx\in \mc{O}_{N_\epsilon}$ as $\bx=(\bx\cdot\bm{n}^*)\bm{n}^* + (\bx- (\bx\cdot\bm{n}^*)\bm{n}^*)$. We then define $E^{\rm end}\bv_{N_\epsilon}$ as
\begin{equation}
  E^{\rm end}\bv_{N_\epsilon}(\bx) = \begin{cases}
    \bv_{N_\epsilon}(\bx)\,, &  \bx\in \mc{V}_{N_\epsilon}\,,\\
    \bv_{N_\epsilon}\big(\bx- 2(\bx\cdot\bm{n}^*)\bm{n}^*)\big)\,, & \bx \in \mc{O}_{N_\epsilon}\,,
  \end{cases}
\end{equation}
i.e. the extension-by-reflection of $\bv_{N_\epsilon}$ into $\mc{O}_{N_\epsilon}$ across the straight-line distance to $\{\X(s)\,:\, s\in[0,1-\epsilon^2]\}$.\footnote{Given that the maximum centerline curvature $\kappa_\star$ is independent of $\epsilon$, for $\epsilon$ sufficiently small, the centerline is effectively straight over the $O(\epsilon^2)$-length endpoint segment $\mc{V}_{N_\epsilon}$, and we may ensure that the reflection $-(\bx\cdot\bm{n}^*)\bm{n}^*+(\bx- (\bx\cdot\bm{n}^*)\bm{n}^*)$ actually lies within $\mc{V}_{N_\epsilon}$.} 
Note that $E^{\rm end}$ coincides with $E^{\rm ref}$ below the cross section $s=1-\epsilon^2$.
By an analogous computation to \cite[Proposition A.5]{free_ends}, we may show 
\begin{equation}
  |\nabla(E^{\rm end}\bv_{N_\epsilon})(\bx)| \le C|\nabla\bv_{N_\epsilon}|
\end{equation}
for $C$ independent of $\epsilon$.

We may now define our extension operator $Z_{\rm ex}$ as
\begin{equation}\label{eq:f_extension}
\begin{aligned}
  Z_{\rm ex}\bv &= \sum_{k=1}^{N_\epsilon-1}\varphi_k\,E^{\rm ref}\bv_k
  + \varphi_{N_\epsilon}(\bx) \,E^{\rm end}\bv_{N_\epsilon} \,.
\end{aligned}  
\end{equation}

We will estimate $Z_{\rm ex}\bv$ using an $\epsilon$-dependent Korn-Poincar\'e inequality on each vessel segment $\mc{V}_k$. First, given any Lipschitz domain $\mc{D}$, we denote the space of rigid motions on $\mc{D}$ by 
 \begin{equation}\label{eq:rigidmotions}
   \bm{R}(\mc{D}) = \{ \bw\in H^1(\mc{D})\,:\, \bw = \bm{A}\bx+\bm{b} \text{ for }\bm{A}=-\bm{A}^{\rm T}\in \R^{3\times 3}\,, \, \bm{b}\in \R^3 \}\,,
 \end{equation}
 and note that the symmetric gradient of any $H^1(\mc{D})$ function vanishes precisely on the set $\bm{R}(\mc{D})$.
By the tangential boundary condition $\bv-(\bv\cdot\bm{n})\bm{n}=0$ on $\Gamma_\epsilon$, we have that for each $k$, if $\bv_k\in \bm{R}(\mc{V}_k)$, then $\bv_k=0$. In particular, on each $\mc{V}_k$, including the endpoint segment $k=N_\epsilon$, we may bound 
\begin{equation}\label{eq:korn_poincare}
  \norm{\bv_k}_{L^2(\mc{V}_k)} \le C\epsilon a(s_k)\norm{\E(\bv_k)}_{L^2(\mc{V}_k)}\,, \quad 
  \norm{\nabla\bv_k}_{L^2(\mc{V}_k)} \le C\norm{\E(\bv_k)}_{L^2(\mc{V}_k)}
\end{equation}
for $C$ independent of $\epsilon$. For the endpoint segment, we recall that $s_{N_\epsilon}=1-\epsilon^2$, and $\epsilon a(1-\epsilon^2)=O(\epsilon^2)$ by the spheroidal endpoint condition \eqref{eq:spheroidal}.
The bounds \eqref{eq:korn_poincare} with undetermined $\epsilon$-dependence may be obtained using a compactness argument (cf. \cite[Corollary A.9 and Lemma A.10]{free_ends}), and the $\epsilon$-dependence follows from scaling (cf. \cite[Corollary A.11]{free_ends}). In particular, using the spheroidal endpoint condition \eqref{eq:spheroidal} and the non-uniform partition \eqref{eq:nonuniform_partition}, for $s\ge1-\delta$, the ratio $a(s_{k+1})/a(s_k)$ is bounded independent of $\epsilon$ (see \cite[A.33]{free_ends}), and each vessel segment $\mc{V}_k$ scales nearly uniformly as $\epsilon a(s_k)$.

Then, for each $k$, using the form \eqref{eq:Erefv} of $E^{\rm ref}\bv$ (or $E^{\rm end}$ on the endpoint segment) and the bound \eqref{eq:korn_poincare}, we may calculate 
\begin{equation}\label{eq:Z_estk}
\begin{aligned}
  \|Z_{\rm ex}\bv\|_{L^2(\mc{O}_k\cup \mc{V}_k)} &\le C\norm{\bv_k}_{L^2(\mc{V}_k)} \le C\epsilon a(s_k)\norm{\E(\bv_k)}_{L^2(\mc{V}_k)}\,,\\
  \|\E(Z_{\rm ex}\bv)\|_{L^2(\mc{O}_k\cup \mc{V}_k)} &\le C\big(\norm{\nabla\bv_k}_{L^2(\mc{V}_k)} + (\epsilon a(s_k))^{-1}\norm{\bv_k}_{L^2(\mc{V}_k)}\big) 
  \le C\norm{\E(\bv_k)}_{L^2(\mc{V}_k)}\,.
\end{aligned}  
\end{equation}
Summing over $k$ and using the nearly-disjoint supports of $\varphi_k$, we have 
\begin{equation}\label{eq:Z1f_ests}
\begin{aligned}
  \|Z_{\rm ex}\bv\|_{L^2(\R^3_+)} \le C\epsilon\norm{\E(\bv)}_{L^2(\mc{V}_\epsilon)}\,, \quad
  \|\E(Z_{\rm ex}\bv)\|_{L^2(\R^3_+)} \le C\norm{\E(\bv)}_{L^2(\mc{V}_\epsilon)}\,.
\end{aligned}  
\end{equation}
We may then use the Korn inequality \eqref{eq:kornR} on $\R^3_+$ to obtain
\begin{equation}\label{eq:freeKorn}
\begin{aligned}
  \norm{\bv}_{L^2(\mc{V}_\epsilon)}&\le \|Z_{\rm ex}\bv\|_{L^2(\R^3_+)} \le C\epsilon\norm{\E(\bv)}_{L^2(\mc{V}_\epsilon)} \,,\\
  \norm{\nabla\bv}_{L^2(\mc{V}_\epsilon)}&\le \|\nabla (Z_{\rm ex}\bv)\|_{L^2(\R^3_+)} \le C\|\E(Z_{\rm ex}\bv)\|_{L^2(\R^3_+)} \le C\norm{\E(\bv)}_{L^2(\mc{V}_\epsilon)} \,.
\end{aligned}  
\end{equation}
\hfill \qedsymbol

\subsection{Straightening}\label{subsec:straightening}
The following construction will be used in each of the subsequent sections.

Given a vessel $\mc{V}_\epsilon$ as in section \ref{sec:setup}, let $\underline{\mc{V}_\epsilon}$ denote a corresponding vessel with the same radius function $a(s)$ on a straight centerline $\X(s)=s\be_z$, where $\be_z$ is the usual Cartesian basis vector.
We define a $C^2$ change of variables $\Psi(\bx)$ mapping the curved vessel $\mc{V}_\epsilon$ to $\underline{\mc{V}_\epsilon}$. Writing $\bx\in \mc{V}_\epsilon$ as $\bx=\X(s)+r\be_r(s,\theta)$, the map $\Psi$ is given by
 \begin{equation}
 \Psi(\bx) = s\be_z+r\be_r(\theta) = s\be_z+r\cos\theta\,\be_x + r\sin\theta\,\be_y\,,
 \end{equation}
 where $(\be_x,\be_y,\be_z)$ are standard Cartesian basis vectors. We may calculate
 \begin{equation}
 \begin{aligned}
\nabla \Psi^{-1} &= (1-r\wh\kappa)\,\be_{\rm t}\otimes\be_z +\be_\theta(s,\theta)\otimes\be_\theta(\theta) + \be_r(s,\theta)\otimes\be_r(\theta)\\
\nabla \Psi(\nabla \Psi)^{\rm T}\circ \Psi^{-1} &= \frac{1}{(1-r\wh\kappa)^2}\be_z\otimes\be_z + \be_\theta(\theta)\otimes\be_\theta(\theta) + \be_r(\theta)\otimes\be_r(\theta) \,,
\end{aligned} 
\end{equation}
from which we obtain 
\begin{equation}
  \abs{\det \nabla\Psi} = \frac{1}{1-r\wh\kappa}\,.
\end{equation}

\subsection{Pressure estimate}\label{subsec:pressure_estimate}
Here we prove the $\epsilon$-scaling of Lemma \ref{lem:bogovskii}.
Within a vessel with straight centerline, the desired $\epsilon$-dependence follows by the following anisotropic scaling argument. Given a radius function $a(s)$ satisfying Definition \ref{def:rad}, we may define a straight vessel at ``scale 1" by 
\begin{equation}
  \underline{\mc{V}_a} = \big\{\bx= s\be_z + r\be_r(\theta)\,: \, 0\le s\le 1\,, \; 0\le r\le a(s)\,, \; 0\le \theta<2\pi \big\}\,,
\end{equation}
with surfaces $\underline{\Gamma_a}$, $\underline{B_a}$ defined analogously to \eqref{eq:Gameps} and \eqref{eq:Beps}.

Given $f\in L^2(\underline{\mc{V}_a})$ with $|\underline{\mc{V}_a}|^{-1}\int_{\underline{\mc{V}_a}}f\,d\bx=f_{\rm c}$, we consider $\bw$ satisfying 
\begin{equation}\label{eq:scale1_eqn}
\begin{aligned}
  \div\bw &= f \qquad \text{in }\underline{\mc{V}_a} \\
  \bw\big|_{\underline{\Gamma_a}} &= 0\,, \quad  
  \bw\big|_{\underline{B_a}} = -f_{\rm c}\frac{|\underline{\mc{V}_a}|}{|\underline{B_a}|}\be_z\,,
\end{aligned}  
\end{equation}
where we emphasize that the normal vector $\bm{n}=-\be_z$ on the bottom cross section $\underline{B_a}$. By the constructions of \cite{bogovskii1980solutions, galdi2011introduction, fabricius2019pressure}, there exists $\bw\in H^1(\underline{\mc{V}_a})$ satisfying \eqref{eq:scale1_eqn} along with the bound
\begin{equation}\label{eq:scale1_est}
  \norm{\nabla\bw}_{L^2(\underline{\mc{V}_a})} \le C\norm{f}_{L^2(\underline{\mc{V}_a})} + C\abs{f_{\rm c}}
  \le C\norm{f}_{L^2(\underline{\mc{V}_a})}\,.
\end{equation}

Let $\underline{\mc{V}_\epsilon}$, $\underline{\Gamma_\epsilon}$, and $\underline{B_\epsilon}$ denote the images of the respective scale-1 objects under the anisotropic rescaling 
\begin{equation}
  (\overline{x},\overline{y},z)= (\epsilon x,\epsilon y, z)\,.
\end{equation}
For $\bw= (w^x,w^y,w^z)$ satisfying \eqref{eq:scale1_eqn}, \eqref{eq:scale1_est} in $\underline{\mc{V}_a}$, we define $\bw_\epsilon$ by
\begin{equation}
\begin{aligned}
  w^x_\epsilon(\overline{x},\overline{y},z) &= \epsilon \,w^x(\epsilon^{-1}\overline{x},\epsilon^{-1}\overline{y},z)  \\
  w^y_\epsilon(\overline{x},\overline{y},z) &= \epsilon\, w^y(\epsilon^{-1}\overline{x},\epsilon^{-1}\overline{y},z)  \\
  w^z_\epsilon(\overline{x},\overline{y},z) &= w^z(\epsilon^{-1}\overline{x},\epsilon^{-1}\overline{y},z)  \,.
\end{aligned}  
\end{equation}
We may then calculate that $\bw_\epsilon$ satisfies 
\begin{equation}
\begin{aligned}
  \div\bw_\epsilon &= f \quad \text{in }\underline{\mc{V}_\epsilon}\,, \qquad \bw_\epsilon\big|_{\underline{\Gamma_\epsilon}} = 0\,, \qquad \bw_\epsilon\big|_{\underline{B_\epsilon}}= -f_{\rm c}\frac{|\underline{\mc{V}_\epsilon}|}{|\underline{B_\epsilon}|} \be_z\,.
\end{aligned} 
\end{equation}
Here we emphasize that $|\underline{\mc{V}_\epsilon}|$ and $|\underline{B_\epsilon}|$ both scale like $\epsilon^2$. Furthermore, the estimate \eqref{eq:scale1_est} becomes
\begin{equation}\label{eq:scale_eps_est}
  \norm{\nabla\bw_\epsilon}_{L^2(\underline{\mc{V}_\epsilon})} 
  \le C\epsilon^{-1}\norm{f}_{L^2(\underline{\mc{V}_\epsilon})}\,.
\end{equation}

We now consider a curved vessel $\mc{V}_\epsilon$ as in section \ref{sec:setup}. Consider $\bv\in H^1(\mc{V}_\epsilon)$ satisfying
\begin{equation}\label{eq:divvg_eqn}
  \div\bv = g \quad \text{in } \mc{V}_\epsilon \,, \quad \bv\big|_{\Gamma_\epsilon}=0\,, \quad \bv\big|_{B_\epsilon}= g_{\rm c}\frac{|\mc{V}_\epsilon|}{|B_\epsilon|}\bm{n}
\end{equation}
and $\norm{\nabla\bv}_{L^2(\mc{V}_\epsilon)}\le C(\epsilon)\norm{g}_{L^2(\mc{V}_\epsilon)}$, where the $\epsilon$-dependence in the bound is to be determined.

For $\bx\in \mc{V}_\epsilon$, let $\by=\Psi(\bx)$ for $\Psi$ as in section \ref{subsec:straightening}, and denote the pushforward of $\bv$ by
\begin{equation}
  \underline\bv(\by) = \nabla\Psi\,\bv(\Psi^{-1}(\by))\,.
\end{equation}
Similarly, let $\underline g(\by)=g(\Psi^{-1}(\by))$. We may then rewrite \eqref{eq:divvg_eqn} as an equation about a straight vessel $\underline{\mc{V}_\epsilon}$:
\begin{equation}
\begin{aligned}
  \frac{1}{\abs{\det\nabla\Psi}}\div_y\big(\abs{\det\nabla\Psi}\,\underline\bv\big) &= \underline g\,, \qquad \by\in \underline{\mc{V}_\epsilon}\\
  \underline{\bv}\big|_{\underline{\Gamma_\epsilon}} &= 0\,, \quad 
  \underline{\bv}\big|_{\underline{B_\epsilon}}= g_{\rm c}\frac{|\mc{V}_\epsilon|}{|B_\epsilon|}\underline{\bm{n}}\,.
\end{aligned}  
\end{equation}
Using that $\int_{\underline{\mc{V}_\epsilon}}\abs{\det\nabla\Psi}\underline g\,d\by=\int_{\mc{V}_\epsilon} g\,d\bx=g_{\rm c}|\mc{V}_\epsilon|$ and that $\underline{B_\epsilon}=B_\epsilon$ and $\underline{\bm{n}}=-\be_z$, by \eqref{eq:scale_eps_est}, we have the estimate
\begin{equation}
  \norm{\nabla_y(\abs{\det\nabla\Psi}\,\underline\bv)}_{L^2(\underline{\mc{V}_\epsilon})} \le C\epsilon^{-1}\norm{\abs{\det\nabla\Psi}\,\underline g}_{L^2(\underline{\mc{V}_\epsilon})}\,.
\end{equation}
Pulling back to $\mc{V}_\epsilon$, we have
\begin{equation}
  \norm{\nabla\Psi^{-\rm T}\nabla_x(\abs{\det\nabla\Psi}\,\nabla\Psi\,\bv)}_{L^2(\mc{V}_\epsilon)} \le C\epsilon^{-1}\norm{g}_{L^2(\mc{V}_\epsilon)}\,,
\end{equation}
where we have used that $\abs{\det\nabla\Psi}\le C$ independent of $\epsilon$. Similarly, since $C^{-1}\le \abs{\nabla \Psi}\le C$, we have
\begin{equation}
\begin{aligned}
  \norm{\nabla\Psi^{-\rm T}\nabla_x(\abs{\det\nabla\Psi}\,\nabla\Psi\,\bv)}_{L^2(\mc{V}_\epsilon)} &\ge C\norm{\nabla\bv}_{L^2(\mc{V}_\epsilon)}-C\norm{\bv}_{L^2(\mc{V}_\epsilon)} \\
  &\ge C\norm{\nabla\bv}_{L^2(\mc{V}_\epsilon)}-C\epsilon\norm{\nabla\bv}_{L^2(\mc{V}_\epsilon)}\,,
\end{aligned}  
\end{equation}
where the Poincar\'e inequality scaling on $\mc{V}_\epsilon$ may be obtained from Lemma \ref{lem:korn}. Note that $C$ depends on the third derivative of the vessel centerline $\X$. For $\epsilon$ sufficiently small, we obtain Lemma \ref{lem:bogovskii}.
\hfill \qedsymbol

\subsection{Trace inequalities}\label{subsec:trace_app}
Here we prove the series of $\epsilon$-dependent trace inequalities stated in Lemma \ref{lem:trace_ineqs}.

\emph{1. Vessel surface $\Gamma_\epsilon$ into vessel exterior $\Omega_\epsilon$.}
It suffices to show the desired inequality \eqref{eq:trace_Gamma2ext} away from the tip, as we will later show that the more refined bound \eqref{eq:trace_tip2exterior} holds for $s\in[1-\epsilon,1]$.

Recalling the definition \eqref{eq:kappastar} of the maximum centerline curvature $\kappa_\star$, we may use the curvilinear cylindrical coordinates of section \ref{subsec:geometry} to define the region
\begin{equation}
\mc{O}_\star = \big\{ \bx(r,\theta,s)\in \Omega_\epsilon\,:\, 0\le s\le 1-\epsilon\,, \; \epsilon a(s)\le r\le \frac{1}{3\kappa_\star}\,, \;0\le \theta<2\pi\,\big\} \,.
\end{equation}
Within $\mc{O}_\star$, we define a radial cutoff 
\begin{equation}
\psi(r,s) = \begin{cases}
  1 & \text{ if } \epsilon a(s) \le r \le \frac{1}{4\kappa_\star}\,, \\
  0  & \text{ if } r\ge \frac{1}{3\kappa_\star}\,.
\end{cases}
\end{equation}
Here we may take $\abs{\p_r\psi}\le C$ independent of $\epsilon$. 
Using the fundamental theorem of calculus, for $g$ on $\Gamma_\epsilon\cap\{s\in[0,1-\epsilon]\}$, we may write 
\begin{equation}\label{eq:FTCv}
  g\big|_{\Gamma_\epsilon} = -\int_{\epsilon a(s)}^{1/2\kappa_\star}\p_r( \psi(r,s) g)\,dr \,.
\end{equation}
We may then bound 
\begin{equation}\label{eq:v_on_Gamma}
\begin{aligned}
  \abs{g\big|_{\Gamma_\epsilon}} &\le \bigg(\int_{\epsilon a(s)}^{1/2\kappa_\star}\frac{1}{r}\,dr\bigg)^{1/2}\bigg(\int_{\epsilon a(s)}^{1/2\kappa_\star}\abs{\p_r (\psi \,g)}^2\,r\,dr\bigg)^{1/2} \\
  &\le C\abs{\log(\epsilon a(s))}^{1/2}\bigg(\int_{\epsilon a(s)}^{1/2\kappa_\star}\abs{\p_r(\psi\, g)}^2\,r\,dr\bigg)^{1/2}\,.
\end{aligned}
\end{equation}

Multiplying \eqref{eq:v_on_Gamma} by the surface element $\mc{J}_\epsilon(s,\theta)$ and integrating in $s$ and $\theta$, we obtain  
\begin{equation}
\begin{aligned}
  \int_0^{1-\epsilon}\int_0^{2\pi}\abs{g\big|_{\Gamma_\epsilon}}^2\, \mc{J}_\epsilon(s,\theta)\,d\theta ds &\le C\abs{\log\epsilon}\int_0^{1-\epsilon}\int_0^{2\pi} \int_{\epsilon a(s)}^{1/2\kappa_\star}\abs{\p_r(\psi\, g)}^2\,r\,dr \,\mc{J}_\epsilon(s,\theta)\,d\theta ds \\
  &\le C\epsilon\abs{\log\epsilon}\big(\norm{\nabla g}_{L^2(\mc{O}_\star)}^2+ \norm{g}_{L^2(\mc{O}_\star)}^2\big)\,.
\end{aligned}
\end{equation}
Here we have used that $\abs{r\wh\kappa}\le 2r\kappa_\star\le \frac{2}{3}$ within $\mc{O}_\star$, and thus the volume element \eqref{eq:volume_element} satisfies $\abs{\mc{J}}=\abs{r(1-r\wh\kappa)} \sim r$ in $\mc{O}_\star$. We further use that $\abs{\p_r\psi}\le C$ independent of $\epsilon$. Noting that 
\begin{equation}
    \norm{g}_{L^2(\mc{O}_\star)} \le \abs{\mc{O}_\star}^{1/3} \norm{g}_{L^6(\mc{O}_\star)} \le C\norm{\nabla g}_{L^2(\Omega_\epsilon)}
\end{equation}
by the $\epsilon$-independent Sobolev equality in the exterior domain $\Omega_\epsilon$ (cf. \cite[Lemma 2.4]{free_ends}), we obtain \eqref{eq:trace_Gamma2ext} away from the tip.

\emph{2. Flat endpoint into vessel interior.}
Let
\begin{equation}
  \mc{C}_\epsilon = \{ \bx(r, s,\theta) \in \mc{V}_\epsilon\,:\, 0< s < \epsilon \} \,,
\end{equation}
i.e. the $O(\epsilon^3)$-volume bottom portion of the vessel. We consider a uniform rescaling of $\mc{C}_\epsilon$ to $O(1)$ volume. Let $\wt{\bx} = \bx/\epsilon$ and define $\mc{C}_1$ to be the rescaled region $\mc{C}_1=\{\wt{\bx} = \bx/\epsilon\,:\, \bx\in \mc{C}_\epsilon\}$ with similarly rescaled base $B_1$ and sides $\Gamma_1$. Any $\wt v\in H^1(\mc{C}_1)$ satisfies the standard $L^2$ trace inequality
\begin{equation}
  \norm{\wt v}_{L^2(B_1)} \le C\big( \norm{\nabla_{\wt x}\wt v}_{L^2(\mc{C}_1)} + \norm{\wt v}_{L^2(\mc{C}_1)} \big)\,.
\end{equation}
Changing variables back to $\bx=\epsilon\wt\bx$, we obtain the desired $\epsilon$-dependent scaling for $v(\bx)=\wt v(\bx/\epsilon)$.\\

\emph{3. Free end region into vessel exterior.}
As we did for the bottom endpoint, we will rely on a scaling argument at the top end region $\Gamma_\epsilon\cap \{1-\epsilon \le s\le 1 \}$. However, we will also make use of additional smallness from the decaying radius $a(s)$ as $s\to 1$.

We begin by considering the case that the vessel centerline is straight at the tip, i.e. $\X(s) = \X_0+ s\be$ for $1-\epsilon\le s\le 1$ and some fixed unit vector $\be$. In this case, let
\begin{equation}
  \underline{\Gamma_\epsilon}^{\rm tip} = \{\bx= \X_0+ s\be + \epsilon a(s)\be_r(\theta)\,: \, 1-\epsilon \le s\le 1\} 
\end{equation}
denote the straight tip, where $\be_r(\theta)$ is a cylindrical basis vector normal to the fixed vector $\be$.  
We further subdivide $\underline{\Gamma_\epsilon}^{\rm tip}$ into
\begin{equation}
   \underline{\Gamma_\epsilon}^{\rm tip,1} = \{\bx\in \underline{\Gamma_\epsilon}^{\rm tip}\,: \, 1-\epsilon \le s\le 1-\epsilon^2\} \,, \qquad
   \underline{\Gamma_\epsilon}^{\rm tip,2} = \{\bx\in \underline{\Gamma_\epsilon}^{\rm tip}\,: \, 1-\epsilon^2 \le s\le 1\}\,.
\end{equation}

We consider $\underline{\Gamma_\epsilon}^{\rm tip,1}$ first. Define the exterior region
\begin{equation}
  \mc{O}_\epsilon^{1} = \{ \bx \in \Omega_\epsilon\,:\, \bx=\X_0+ s\be + r\be_r(\theta)\,, \; 1-\epsilon< s < 1-\epsilon^2\,,\; \epsilon a(s)<r<\epsilon \}\,. 
\end{equation}
Letting $\wt{\bx}=\bx/\epsilon$, we define the rescaled region $\mc{O}_a^{1}=\{\wt{\bx} = \bx/\epsilon\,: \, \bx\in\mc{O}_\epsilon^{1}\}$, and similarly denote the rescaled surface $\underline{\Gamma_a}^{1}=\{\wt{\bx} = \bx/\epsilon\,: \, \bx\in \underline{\Gamma_\epsilon}^{1}\}$. Within $\mc{O}_a^{1}$, we consider cylindrical coordinates $\sigma = \frac{s-(1-\epsilon)}{\epsilon}$, $\rho=r/\epsilon$, and $\theta$, so that $0\le\sigma\le 1-\epsilon$ on $\underline{\Gamma_a}^{1}$. We define $\wt a(\sigma)=a(\epsilon\sigma+1-\epsilon)$ and note that $\epsilon \le \wt a(\sigma)\le \epsilon^{1/2}$.
We may calculate the surface element on $\underline{\Gamma_a}^{1}$ as $\mc{J}_a(\sigma)=\wt a(\sigma)\sqrt{1+(\p_\sigma\wt a)^2}$, where we note that $\abs{\wt a(\sigma)\p_\sigma\wt a(\sigma)} \le C\epsilon$.

Let $\phi_a(\rho,\sigma)$ be a smooth cutoff in $\rho$ within $\mc{O}_a^{1}$ satisfying
\begin{equation}
  \phi_a(\rho,\sigma) = \begin{cases}
    1\,, & \wt a(\sigma)\le \rho \le \frac{1}{2}\\
    0\,, & \rho \ge 1 \,.
  \end{cases}
\end{equation}
Let $\wt g(\wt\bx)=g(\epsilon\wt\bx)$. On $\underline{\Gamma_a}^{1}$, as in \eqref{eq:FTCv}, we may write 
\begin{equation}
  \wt g\big|_{\underline{\Gamma_a}^{1}} = -\int_{\wt{a}(\sigma)}^1 \p_\rho (\phi_a \,\wt g)\,d\rho\,. 
\end{equation}
We then have
\begin{equation}
\begin{aligned}
  \abs{\wt g\big|_{\underline{\Gamma_a}^{1}}} &\le \bigg(\int_{\wt{a}}^1 \frac{1}{\rho}\,d\rho\bigg)^{1/2}\bigg(\int_{\wt{a}}^1\abs{\p_\rho (\phi_a\, \wt g)}^2\,\rho\,d\rho\bigg)^{1/2} \\
  &\le C\abs{\log(\wt{a}(\sigma))}^{1/2}\bigg(\int_{\wt{a}}^1\abs{\p_\rho(\phi_a\, \wt g)}^2\,\rho\,d\rho\bigg)^{1/2}\,.
\end{aligned}
\end{equation}
Multiplying by the surface element $\mc{J}_a$ on $\underline{\Gamma_a}^{1}$ and integrating in $\sigma$ and $\theta$, we have 
\begin{equation}\label{eq:compute_Teps_trace}
\begin{aligned}
  \int_0^{1-\epsilon}\int_0^{2\pi}\abs{\wt g\big|_{\underline{\Gamma_a}^{1}}}^2\, \mc{J}_a\,d\theta d\sigma &\le C\int_0^{1-\epsilon}\int_0^{2\pi} \abs{\log \wt{a}(\sigma)}\int_{\wt{a}}^1\abs{\p_\rho(\phi_a \wt g)}^2\,\rho\,d\rho \,\mc{J}_a\,d\theta d\sigma \\
  &\le C\sup_{\sigma\in[0,1-\epsilon]}\abs{\wt{a}\log \wt{a}}\big(\norm{\nabla \wt g}_{L^2(\mc{O}_a^{1})}^2+ \norm{\wt g}_{L^2(\mc{O}_a^{1})}^2\big)\\
  &\le C\epsilon^{1/2}\abs{\log\epsilon}\big(\norm{\nabla \wt g}_{L^2(\mc{O}_a^{1})}^2+ |\mc{O}_a^{1}|^{2/3}\norm{\wt g}_{L^6(\mc{O}_a^{1})}^2\big)\\
  &\le C\epsilon^{1/2}\abs{\log\epsilon}\norm{\nabla \wt g}_{L^2(\mc{O}_a^{\rm f,1})}^2
\end{aligned}
\end{equation}
by the Sobolev inequality in $\mc{O}_a^{1}$.
Using the uniform rescaling $\bx=\epsilon\wt{\bx}$ yields 
\begin{equation}\label{eq:tipest1}
  \norm{g}_{L^2(\underline{\Gamma_\epsilon}^{\rm tip,1})} \le C\epsilon^{5/4}\abs{\log\epsilon}^{1/2}\norm{\nabla g}_{L^2(\mc{O}_\epsilon^{1})}\,.
\end{equation}

On $\underline{\Gamma_\epsilon}^{\rm tip,2}$, we may make use of the fact that $a(s)\le \epsilon$ by the spheroidal endpoint condition \eqref{eq:spheroidal}. In particular, we may consider the $O(\epsilon^6)$-volume region 
\begin{equation}\label{eq:Oepsf2}
   \mc{O}_\epsilon^{2} = \{ \bx \in \Omega_\epsilon\,:\, {\rm dist}(\bx,\underline{\Gamma_\epsilon}^{\rm tip,2})\le \epsilon^2  \}\,.
\end{equation} 
Let $\mc{O}_1^{2}$ denote the region \eqref{eq:Oepsf2} under the uniform rescaling $\wh\bx=\bx/\epsilon^2$, and define $\underline{\Gamma_1}^{2}$ analogously. We have that $\wh g(\wh\bx) = g(\epsilon^2\bx)$ satisfies the following trace inequality in $\mc{O}_1^{2}$: 
\begin{equation}
  \norm{\wh g}_{L^2(\underline{\Gamma_1}^{2})} \le C\norm{\nabla \wh g}_{L^2(\mc{O}_1^{2})}\,,
\end{equation}
where, as in \eqref{eq:compute_Teps_trace}, we may use the Sobolev inequality on $\mc{O}_1^{2}$ to absorb the lower part of the norm on the right hand side. The uniform rescaling $\bx=\epsilon^2\wh\bx$ yields 
\begin{equation}\label{eq:tipest2}
  \norm{g}_{L^2(\underline{\Gamma_\epsilon}^{\rm tip,2})} \le C\epsilon^2\norm{\nabla g}_{L^2(\mc{O}_\epsilon^{2})}\,.
\end{equation}

Finally, we need to incorporate the effects of a curved centerline $\X(s)$. Let $\Gamma_\epsilon^{\rm tip}$ denote the tip with a general curved centerline as in section \ref{sec:setup}. For $\bx\in \Gamma_\epsilon^{\rm tip}$, we use the change of variables $\bx=\Psi^{-1}(\by)$ from section \ref{subsec:straightening} to map to $\underline{\Gamma_\epsilon}^{\rm tip}$. 

By \eqref{eq:tipest1} and \eqref{eq:tipest2}, we have
\begin{equation}
\begin{aligned}
  \norm{g}_{L^2(\Gamma_\epsilon^{\rm tip})} 
  &=\norm{g(\Psi^{-1}(\cdot))}_{L^2(\underline{\Gamma_\epsilon}^{\rm tip})} 
  \le C\epsilon^{5/4}\abs{\log\epsilon}^{1/2}\norm{\nabla_y \,g(\Psi^{-1}(\cdot))}_{L^2(\mc{O}_\epsilon^{1}\cup\mc{O}_\epsilon^{2})}\\
  &\le C\epsilon^{5/4}\abs{\log\epsilon}^{1/2}\norm{\nabla\Psi^{-\rm T}\nabla_x \,g(\Psi^{-1}(\cdot))}_{L^2(\Psi^{-1}(\mc{O}_\epsilon^{1}\cup\mc{O}_\epsilon^{2}))}\,.
\end{aligned}  
\end{equation}
Using that $C^{-1}\le \abs{\nabla\Psi}\le C$ for $C$ independent of $\epsilon$, we obtain the bound \eqref{eq:trace_tip2exterior}.
\hfill \qedsymbol

%% file: partII.bbl
\begin{thebibliography}{10}

\bibitem{berrone2023optimization}
S.~Berrone, C.~Giverso, D.~Grappein, L.~Preziosi, and S.~Scial{\`o}.
\newblock An optimization based 3d-1d coupling strategy for tissue perfusion and chemical transport during tumor-induced angiogenesis.
\newblock {\em Computers \& Mathematics with Applications}, 151:252--270, 2023.

\bibitem{bishop1975there}
R.~L. Bishop.
\newblock There is more than one way to frame a curve.
\newblock {\em The American Mathematical Monthly}, 82(3):246--251, 1975.

\bibitem{blanco2007unified}
P.~Blanco, R.~Feij{\'o}o, and S.~Urquiza.
\newblock A unified variational approach for coupling 3d--1d models and its blood flow applications.
\newblock {\em Computer Methods in Applied Mechanics and Engineering}, 196(41-44):4391--4410, 2007.

\bibitem{blanco2009potentialities}
P.~J. Blanco, M.~Pivello, S.~Urquiza, and R.~A. Feij{\'o}o.
\newblock On the potentialities of 3d--1d coupled models in hemodynamics simulations.
\newblock {\em Journal of biomechanics}, 42(7):919--930, 2009.

\bibitem{bogovskii1980solutions}
M.~Bogovski{\i}.
\newblock Solutions of some problems of vector analysis, associated with the operators div and grad.
\newblock {\em Theory of cubature formulas and the application of functional analysis to problems of mathematical physics}, 1980:5--40, 1980.

\bibitem{boulakia2024mathematical}
M.~Boulakia, C.~Grandmont, F.~Lespagnol, and P.~Zunino.
\newblock Mathematical and numerical analysis of reduced order interface conditions and augmented finite elements for mixed dimensional problems.
\newblock {\em Computers \& Mathematics with Applications}, 175:536--569, 2024.

\bibitem{canic2003effective}
S.~Canic and A.~Mikelic.
\newblock Effective equations modeling the flow of a viscous incompressible fluid through a long elastic tube arising in the study of blood flow through small arteries.
\newblock {\em SIAM Journal on Applied Dynamical Systems}, 2(3):431--463, 2003.

\bibitem{cassot2010branching}
F.~Cassot, F.~Lauwers, S.~Lorthois, P.~Puwanarajah, V.~Cances-Lauwers, and H.~Duvernoy.
\newblock Branching patterns for arterioles and venules of the human cerebral cortex.
\newblock {\em Brain research}, 1313:62--78, 2010.

\bibitem{castineira2019rigorous}
G.~Castineira, E.~Maru{\v{s}}i{\'c}-Paloka, I.~Pa{\v{z}}anin, and J.~M. Rodriguez.
\newblock Rigorous justification of the asymptotic model describing a curved-pipe flow in a time-dependent domain, 2019.

\bibitem{d2007multiscale}
C.~D'Angelo.
\newblock {\em Multiscale modelling of metabolism and transport phenomena in living tissues}.
\newblock PhD thesis, EPFL, 2007.

\bibitem{d2012finite}
C.~D'Angelo.
\newblock Finite element approximation of elliptic problems with dirac measure terms in weighted spaces: applications to one-and three-dimensional coupled problems.
\newblock {\em SIAM Journal on Numerical Analysis}, 50(1):194--215, 2012.

\bibitem{d2008coupling}
C.~D'angelo and A.~Quarteroni.
\newblock On the coupling of 1d and 3d diffusion-reaction equations: application to tissue perfusion problems.
\newblock {\em Mathematical Models and Methods in Applied Sciences}, 18(08):1481--1504, 2008.

\bibitem{d2011robust}
C.~D'Angelo and P.~Zunino.
\newblock Robust numerical approximation of coupled stokes' and darcy's flows applied to vascular hemodynamics and biochemical transport.
\newblock {\em ESAIM: Mathematical Modelling and Numerical Analysis}, 45(3):447--476, 2011.

\bibitem{discacciati2009navier}
M.~Discacciati, A.~Quarteroni, et~al.
\newblock Navier-stokes/darcy coupling: modeling, analysis, and numerical approximation.
\newblock {\em Rev. Mat. Complut}, 22(2):315--426, 2009.

\bibitem{fabricius2019pressure}
J.~Fabricius, E.~Miroshnikova, A.~Tsandzana, and P.~Wall.
\newblock Pressure-driven flow in thin domains.
\newblock {\em Asymptotic Analysis}, 116(1):1--26, 2019.

\bibitem{formaggia2006coupling}
L.~Formaggia, A.~Moura, and F.~Nobile.
\newblock Coupling 3d and 1d fluid-structure interaction models for blood flow simulations.
\newblock In {\em PAMM: Proceedings in Applied Mathematics and Mechanics}, volume~6, pages 27--30. Wiley Online Library, 2006.

\bibitem{formaggia2007stability}
L.~Formaggia, A.~Moura, and F.~Nobile.
\newblock On the stability of the coupling of 3d and 1d fluid-structure interaction models for blood flow simulations.
\newblock {\em ESAIM: Mathematical Modelling and Numerical Analysis}, 41(4):743--769, 2007.

\bibitem{galdi2011introduction}
G.~P. Galdi.
\newblock {\em An introduction to the mathematical theory of the {N}avier-{S}tokes equations: Steady-state problems}.
\newblock Springer Science \& Business Media, 2011.

\bibitem{ghosh2021modified}
A.~Ghosh, V.~Kozlov, and S.~Nazarov.
\newblock Modified reynolds equation for steady flow through a curved pipe.
\newblock {\em Journal of Mathematical Fluid Mechanics}, 23:1--22, 2021.

\bibitem{heltai2023reduced}
L.~Heltai and P.~Zunino.
\newblock Reduced lagrange multiplier approach for non-matching coupling of mixed-dimensional domains.
\newblock {\em Mathematical Models and Methods in Applied Sciences}, 33(12):2425--2462, 2023.

\bibitem{hsu1989green}
R.~Hsu and T.~W. Secomb.
\newblock A green's function method for analysis of oxygen delivery to tissue by microvascular networks.
\newblock {\em Mathematical biosciences}, 96(1):61--78, 1989.

\bibitem{hyde2013parameterisation}
E.~R. Hyde, C.~Michler, J.~Lee, A.~N. Cookson, R.~Chabiniok, D.~A. Nordsletten, and N.~P. Smith.
\newblock Parameterisation of multi-scale continuum perfusion models from discrete vascular networks.
\newblock {\em Medical \& biological engineering \& computing}, 51:557--570, 2013.

\bibitem{johnston2009tortuosity}
S.~C. Johnston, D.~K. Wallace, S.~F. Freedman, T.~L. Yanovitch, and Z.~Zhao.
\newblock Tortuosity of arterioles and venules in quantifying plus disease.
\newblock {\em Journal of American Association for Pediatric Ophthalmology and Strabismus}, 13(2):181--185, 2009.

\bibitem{koch2020modeling}
T.~Koch, M.~Schneider, R.~Helmig, and P.~Jenny.
\newblock Modeling tissue perfusion in terms of 1d-3d embedded mixed-dimension coupled problems with distributed sources.
\newblock {\em Journal of Computational Physics}, 410:109370, 2020.

\bibitem{koppl20203d}
T.~K{\"o}ppl, E.~Vidotto, and B.~Wohlmuth.
\newblock A 3d-1d coupled blood flow and oxygen transport model to generate microvascular networks.
\newblock {\em International journal for numerical methods in biomedical engineering}, 36(10):e3386, 2020.

\bibitem{kuchta2021analysis}
M.~Kuchta, F.~Laurino, K.-A. Mardal, and P.~Zunino.
\newblock Analysis and approximation of mixed-dimensional pdes on 3d-1d domains coupled with lagrange multipliers.
\newblock {\em SIAM Journal on Numerical Analysis}, 59(1):558--582, 2021.

\bibitem{laurino2019derivation}
F.~Laurino and P.~Zunino.
\newblock Derivation and analysis of coupled pdes on manifolds with high dimensionality gap arising from topological model reduction.
\newblock {\em ESAIM: Mathematical Modelling and Numerical Analysis}, 53(6):2047--2080, 2019.

\bibitem{lorthois2014tortuosity}
S.~Lorthois, F.~Lauwers, and F.~Cassot.
\newblock Tortuosity and other vessel attributes for arterioles and venules of the human cerebral cortex.
\newblock {\em Microvascular research}, 91:99--109, 2014.

\bibitem{manjate2024mathematical}
S.~R. Manjate.
\newblock Mathematical modelling of pressure-driven flow in thin domains.
\newblock {\em Doctoral thesis, Lule\o{a} University of Technology}, 2024.

\bibitem{maruvsic2000two}
S.~Maru{\v{s}}i{\'c} and E.~Maru{\v{s}}i{\'c}-Paloka.
\newblock Two-scale convergence for thin domains and its applications to some lower-dimensional models in fluid mechanics.
\newblock {\em Asymptotic Analysis}, 23(1):23--57, 2000.

\bibitem{maruvsic2001effects}
E.~Maru{\v{s}}i{\'c}-Paloka.
\newblock The effects of flexion and torsion on a fluid flow through a curved pipe.
\newblock {\em Applied mathematics and optimization}, 44:245--272, 2001.

\bibitem{mironescu2018role}
P.~Mironescu.
\newblock The role of the hardy type inequalities in the theory of function spaces.
\newblock {\em Revue roumaine de math{\'e}matiques pures et appliqu{\'e}es}, 63(4):447--525, 2018.

\bibitem{miroshnikova2020pressure}
E.~Miroshnikova.
\newblock Pressure-driven flow in a thin pipe with rough boundary.
\newblock {\em Zeitschrift f{\"u}r angewandte Mathematik und Physik}, 71(4):138, 2020.

\bibitem{rigid}
Y.~Mori and L.~Ohm.
\newblock An error bound for the slender body approximation of a thin, rigid fiber sedimenting in stokes flow.
\newblock {\em Research in the Mathematical Sciences}, 7(2):1--27, 2020.

\bibitem{inverse}
Y.~Mori and L.~Ohm.
\newblock Accuracy of slender body theory in approximating force exerted by thin fiber on viscous fluid.
\newblock {\em Studies in Applied Mathematics}, 2021.

\bibitem{closed_loop}
Y.~Mori, L.~Ohm, and D.~Spirn.
\newblock Theoretical justification and error analysis for slender body theory.
\newblock {\em Communications on Pure and Applied Mathematics}, 73(6):1245--1314, 2020.

\bibitem{free_ends}
Y.~Mori, L.~Ohm, and D.~Spirn.
\newblock Theoretical justification and error analysis for slender body theory with free ends.
\newblock {\em Archive for Rational Mechanics and Analysis}, 235(3):1905--1978, 2020.

\bibitem{nishimura2007penetrating}
N.~Nishimura, C.~B. Schaffer, B.~Friedman, P.~D. Lyden, and D.~Kleinfeld.
\newblock Penetrating arterioles are a bottleneck in the perfusion of neocortex.
\newblock {\em Proceedings of the National Academy of Sciences}, 104(1):365--370, 2007.

\bibitem{nobile2009coupling}
F.~Nobile.
\newblock Coupling strategies for the numerical simulation of blood flow in deformable arteries by 3d and 1d models.
\newblock {\em Mathematical and Computer Modelling}, 49(11-12):2152--2160, 2009.

\bibitem{notaro2016mixed}
D.~Notaro, L.~Cattaneo, L.~Formaggia, A.~Scotti, and P.~Zunino.
\newblock A mixed finite element method for modeling the fluid exchange between microcirculation and tissue interstitium.
\newblock {\em Advances in discretization methods: discontinuities, virtual elements, fictitious domain methods}, pages 3--25, 2016.

\bibitem{laplace}
L.~Ohm.
\newblock On an angle-averaged {N}eumann-to-{D}irichlet map for thin filaments.
\newblock {\em arXiv preprint arXiv:2308.06592}, 2023.

\bibitem{partI}
L.~Ohm and S.~Strikwerda.
\newblock A hierarchy of blood vessel models, part {I}: 3{D}-1{D} to 1{D}.
\newblock {\em arXiv preprint}, 2025.

\bibitem{pittman1995influence}
R.~N. Pittman.
\newblock Influence of microvascular architecture on oxygen exchange in skeletal muscle.
\newblock {\em Microcirculation}, 2(1):1--18, 1995.

\bibitem{possenti2021mesoscale}
L.~Possenti, A.~Cicchetti, R.~Rosati, D.~Cerroni, M.~L. Costantino, T.~Rancati, and P.~Zunino.
\newblock A mesoscale computational model for microvascular oxygen transfer.
\newblock {\em Annals of Biomedical Engineering}, 49:3356--3373, 2021.

\bibitem{pries2008blood}
A.~R. Pries and T.~W. Secomb.
\newblock Blood flow in microvascular networks.
\newblock In {\em Microcirculation}, pages 3--36. Elsevier, 2008.

\bibitem{qi2024hemodynamic}
Y.~Qi, S.-S. Chang, Y.~Wang, C.~Chen, K.~I. Baek, T.~Hsiai, and M.~Roper.
\newblock Hemodynamic regulation allows stable growth of microvascular networks.
\newblock {\em Proceedings of the National Academy of Sciences}, 121(9):e2310993121, 2024.

\bibitem{qi2021control}
Y.~Qi and M.~Roper.
\newblock Control of low flow regions in the cortical vasculature determines optimal arterio-venous ratios.
\newblock {\em Proceedings of the National Academy of Sciences}, 118(34):e2021840118, 2021.

\bibitem{secomb2004green}
T.~W. Secomb, R.~Hsu, E.~Y. Park, and M.~W. Dewhirst.
\newblock Green's function methods for analysis of oxygen delivery to tissue by microvascular networks.
\newblock {\em Annals of biomedical engineering}, 32:1519--1529, 2004.

\bibitem{shih2015robust}
A.~Y. Shih, C.~R{\"u}hlmann, P.~Blinder, A.~Devor, P.~J. Drew, B.~Friedman, P.~M. Knutsen, P.~D. Lyden, C.~Mateo, L.~Mellander, et~al.
\newblock Robust and fragile aspects of cortical blood flow in relation to the underlying angioarchitecture.
\newblock {\em Microcirculation}, 22(3):204--218, 2015.

\bibitem{shipley2020hybrid}
R.~J. Shipley, A.~F. Smith, P.~W. Sweeney, A.~R. Pries, and T.~W. Secomb.
\newblock A hybrid discrete--continuum approach for modelling microcirculatory blood flow.
\newblock {\em Mathematical medicine and biology: a journal of the IMA}, 37(1):40--57, 2020.

\bibitem{sweeney2024three}
P.~W. Sweeney, C.~Walsh, S.~Walker-Samuel, and R.~J. Shipley.
\newblock A three-dimensional, discrete-continuum model of blood pressure in microvascular networks.
\newblock {\em International Journal for Numerical Methods in Biomedical Engineering}, 40(8):e3832, 2024.

\bibitem{zunino2002mathematical}
P.~Zunino.
\newblock {\em Mathematical and numerical modeling of mass transfer in the vascular system}.
\newblock PhD thesis, Ecole Polytechnique Fédérale de Lausanne, Switzerland, 2002.

\end{thebibliography}
